\documentclass[10pt]{article}

\usepackage{graphics}

\usepackage{amsfonts}

\usepackage{amscd}
\usepackage{amsthm}
\usepackage{indentfirst}

\usepackage[all]{xy}

\usepackage{amssymb, amsmath, amsthm, amsgen, amstext, amsbsy, amsopn}
\usepackage{epic,eepic}

\newtheorem{thm}{Theorem}[section]
\newtheorem{lemma}[thm]{Lemma}
\newtheorem{prop}[thm]{Proposition}
\newtheorem{cor}[thm]{Corollary}

\newtheorem{rem}[thm]{Remark}
\newtheorem{exam}[thm]{Example}

\newtheorem{question}[thm]{Question}

\newcommand{\R}{{\mathbb{R}}}

\newcommand{\A}{{\mathbb{A}}}

\newcommand{\T}{{\mathbb{T}}}
\newcommand{\Z}{{\mathbb{Z}}}
\newcommand{\N}{{\mathbb{N}}}
\newcommand{\C}{{\mathbb{C}}}

\newcommand{\D}{{\mathbb{D}}}

\newcommand{\cD}{{\mathcal{D}}}

\newcommand{\cF}{{\mathcal{F}}}

\newcommand{\cJ}{{\mathcal{J}}}
\newcommand{\cK}{{\mathcal{K}}}
\newcommand{\cL}{{\mathcal{L}}}
\newcommand{\cM}{{\mathcal{M}}}

\newcommand{\cP}{{\mathcal{P}}}

\newcommand{\cT}{{\mathcal{T}}}

\newcommand{\cX}{{\mathcal{X}}}
\newcommand{\cY}{{\mathcal{Y}}}

\newcommand{\s}{{\text{stab}}}

\def\id{{1\hskip-2.5pt{\rm l}}}

\newcommand{\tzeta}{{\widetilde{\zeta}}}

\newcommand{\Symp}{{\hbox{\it Symp} }}

\newcommand{\Qed}{\ \hfill \qedsymbol \bigskip}

\newcommand{\op}{{ \it Op}}

\DeclareMathOperator{\sgrad}{sgrad}

\begin{document}

\title{ Poisson brackets and symplectic invariants}

\renewcommand{\thefootnote}{\alph{footnote}}

\author{\textsc Lev Buhovsky$^{a}$, Michael Entov$^{b}$,\ Leonid
Polterovich$^{c}$ }

\footnotetext[1]{ This author also uses the spelling ``Buhovski"
for his family name.}

\footnotetext[2]{Partially supported by the Israel Science
Foundation grant $\#$ 723/10 and by the Japan Technion Society
Research Fund.}

\footnotetext[3]{ Partially supported by the National Science
Foundation grant DMS-1006610. }

\date{\today}

\maketitle

\begin{abstract}

\noindent We introduce new invariants associated to collections of
compact sub\-sets of a symplectic manifold. They are defined
through an elementa\-ry-looking variational problem involving
Poisson brackets. The proof of the non-triviality of these
invariants involves various flavors of Floer theory, including the
$\mu^3$-operation in Donaldson-Fukaya category. We present
applications to approximation theory on symplec\-tic manifolds and
to Hamiltonian dynamics.

\end{abstract}

\bigskip
\noindent
{\bf MSC classes:} 	53Dxx, 37J05

\bigskip
\noindent
{\bf Keywords:} symplectic manifold, Poisson brackets, Hamiltonian chord, quasi-state, Donaldson-Fukaya category

\tableofcontents

\vfill\eject

\section{Introduction and main results}

\subsection{$C^0$-robustness of the Poisson bracket}
\label{subsec-robustness}

Let $(M^{2n},\omega)$ be a symplectic manifold.  Consider the
space $C^\infty_c (M)$ of smooth compactly supported functions on
$M$ equipped with the uniform norm $\| F\| := \max_{x\in M}
|F(x)|$ and with the Poisson bracket $\{F,G\}$. Most of the action
in the present paper takes place in the space $\cF = C^\infty_c
(M)\times C^\infty_c (M)$. It was established in
\cite{EP-Poisson1} (cf. \cite{Cardin-Viterbo, Z-JMD}) that the
functional
$$ \cF \to [0;+\infty), (F,G) \mapsto ||\{F,G\}||,$$
is lower semi-continuous in the uniform norm, meaning that
\begin{equation}\label{eq-lowsem}
\liminf_{\overline{F},\overline{G}\stackrel{C^0}{\longrightarrow}F,G}
\|\{\overline{F},\overline{G}\}\|= \| \{ F, G\} \| \;\;\forall F,G
\in \cF.
\end{equation}
This result can be considered as a manifestation of symplectic
rigidity in the function space $\cF$. The surprising feature here is
that the Poisson bracket involves first derivatives of functions,
while the convergence in \eqref{eq-lowsem} is only in the
$C^0$-sense.

The main observation of the present paper is that certain
variational problems involving the functional $(F,G) \mapsto
||\{F,G\}||$ give rise to invariants of (collections of) compact
subsets of symplectic manifolds. Even though their definition
involves only elementary calculus, their study is based on a variety
of ``hard" symplectic methods such as Gromov's pseudo-holomorphic
curves, Floer theory, Donaldson-Fukaya category and symplectic field
theory. The applications of these invariants include approximation
theory on symplectic manifolds and Hamiltonian dynamics.

\medskip
\noindent
\subsection{Introducing the Poisson bracket
invariants}\label{subsec-intro-intro}

We introduce the following two versions of the Poisson bracket
invariants.

\begin{figure}
 \begin{center}
\scalebox{0.5}{\includegraphics*{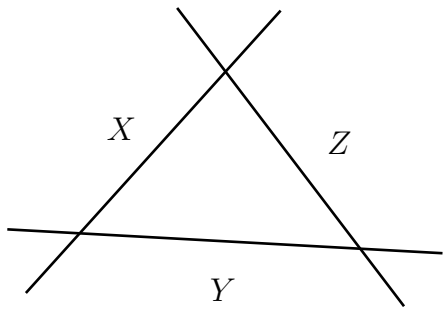}}
 \end{center}
  \caption{$X \cap Y \cap Z=\emptyset$}
   \label{fig1}
\end{figure}

\medskip
\noindent{\sc Invariants of triples:} Let $X,Y,Z \subset M$ be a
triple of compact sets. Put
$$pb_3(X,Y,Z) = \inf ||\{F,G\}||,$$
where the infimum is taken over the class
\begin{equation}\label{eq-F3}
\cF_3 (X,Y,Z) :=\\
\{ (F,G)\ |\  F|_X \leq 0, G|_Y \leq 0, (F+G)|_Z \geq 1\ \}
\end{equation}
of pairs of functions from $\cF$. Note that this class is non-empty
whenever
\begin{equation}\label{eq-inter-3}
X \cap Y \cap Z=\emptyset,
\end{equation}
see Figure \ref{fig1}.
Indeed, it contains any partition of unity subordinated to the
covering $(M \setminus X, M \setminus Y, M \setminus Z)$ of $M$.
If the latter condition is violated, we put $pb_3(X,Y,Z) =
+\infty$.

An easy check shows that $pb_3 (X,Y,Z)$ is symmetric with respect
to $X,Y,Z$.

The next toy example shows that this variational problem is
non-trivial.

\begin{exam}
\label{exam-S2}
{\rm  Consider the sphere
$$S^2 = \{ (x,y,z)\in\R^3\ |\ x^2 +y^2+z^2=1\}$$
with the standard symplectic form. Take three big circles
$$X=\{x=0\},\; Y = \{y=0\},\; Z=\{z=0\}.$$
It turns out that $pb_3(X,Y,Z) > 0$,  see
Example~\ref{exam-S2-superheavy} below. Later on we shall discuss
various generalizations of this example to higher dimensions and
to singular subsets which are not necessarily submanifolds.}
\end{exam}

\begin{figure}
 \begin{center}
\scalebox{0.5}{\includegraphics*{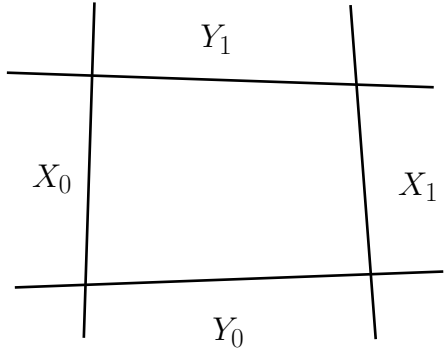}}
 \end{center}
  \caption{$X_0 \cap X_1 = Y_0 \cap Y_1 =\emptyset$}
   \label{fig1a}
\end{figure}

\medskip
\noindent{\sc Invariants of quadruples:} Let $X_0,X_1,Y_0,Y_1
\subset M$ be a quadruple of compact sets. Put
$$pb_4(X_0,X_1,Y_0,Y_1) = \inf ||\{F,G\}||,$$
where the infimum is taken over the class
\begin{multline}\label{eq-F4} \cF_4 (X_0, X_1, Y_0,Y_1) :=  \\ \{ (F,G)\ |\
F|_{X_0} \leq 0,\ F|_{X_1} \geq 1,\ G|_{Y_0} \leq 0,\ G|_{Y_1} \geq
1\ \}\end{multline} of pairs of functions from $\cF$. Note that this
class is non-empty whenever
\begin{equation}\label{eq-inter-4}
X_0 \cap X_1 = Y_0 \cap Y_1 =\emptyset,\end{equation}
see Figure \ref{fig1a}.

If the
latter condition is violated, we put $pb_4(X_0,X_1,Y_0,Y_1) =
+\infty$.

\begin{exam}
\label{exam-pb4m3} {\rm Consider two parallels $X_0$ and $X_1$ and
two meridians $Y_0$ and $Y_1$ on a two-dimensional torus $\T^2$.
They divide $\T^2$ into four squares. Pick three of the squares,
attach a handle to each of them and call the resulting genus-$4$
surface $M$. Equip $M$ with an area form $\omega$. We shall see in
Remark~\ref{rem-comparison-square-4-meridians} (or, alternatively,
in Section~\ref{subsec-intro-pb-lagrangian}) that
$pb_4(X_0,X_1,Y_0,Y_1) > 0$. Furthermore, this example is stable
in the following sense. Consider the product $M \times T^*S^1$
equipped with the symplectic form $\omega + dp \wedge dq$. Let
$K_1,\ldots,K_4$ be any four exact sections\footnote{Here and
further on by an exact section of a cotangent bundle we mean a
graph of the differential of a smooth function on the base.} of
$T^*S^1$. Then $pb_4({X}_0\times K_1,X_1\times K_2,{Y}_0\times
K_3,{Y}_1\times K_4)
>0$ (see Theorem~\ref{thm-lag-main} and Remark~\ref{rm-stb-more}).
Exactly the same conclusion holds true for the quadruple of
circles $X_0,X_1,Y_0,Y_1$ on the torus $\T^2$ (no handles
attached) as well as its stabilization by four exact sections of
$T^*S^1$, see Remark \ref{rem-nohandles}. }
\end{exam}

\medskip
One can easily check that $pb_4(X_0,X_1,Y_0,Y_1)$ does not change
under permutations which switch $X_0$ with $X_1$, $Y_0$ with $Y_1$
or the pair $(X_0,X_1)$ with the pair $(Y_0,Y_1)$.

In what follows we shall often use a slightly different but
equivalent definition of the Poisson bracket invariants. Given a
closed subset $X\subset M$, we denote by $\op(X)$ a sufficiently
small neighborhood of $X$. When we say that $F=0$ on $\op(X)$ we
mean that $F$ vanishes on {\it some} neighborhood of $X$. For a
triple of compact subsets $X,Y,Z \subset M$ satisfying
\eqref{eq-inter-3} define a class

\medskip
\noindent
\begin{multline*}
\cF'_3 (X,Y,Z) := \{ (F,G)\ |\ F \geq 0, G \geq 0, F+G \leq 1,\\
F|_{\op(X)} = 0, G|_{\op(Y)} = 0, (F+G)|_{\op(Z)} = 1\ \}.
\end{multline*}
 Similarly, for a quadruple of compact sets $X_0,X_1,Y_0,Y_1
\subset M$ satisfying \eqref{eq-inter-4} put
\begin{multline*}
\cF'_4 (X_0, X_1, Y_0,Y_1):= \{ (F,G)\ |\ 0\leq F\leq 1,\
F|_{\op(X_0)} = 0,\ F|_{\op(X_1)} =
1,\\
 0\leq G\leq 1,\ G|_{\op(Y_0)} = 0,\ G|_{\op(Y_1)} = 1\ \}.
\end{multline*}

\begin{prop} \label{prop-op}
\begin{equation}\label{eq-op-3}
pb_3(X,Y,Z) = \inf_{(F,G)\in \cF_3' (X,Y,Z)} ||\{F,G\}||
\end{equation}
and
\begin{equation}\label{eq-op-4}
pb_4(X_0,X_1,Y_0,Y_1)  = \inf_{(F,G)\in \cF'_4 (X_0,X_1,Y_0,Y_1)}
||\{F,G\}||.
\end{equation}
\end{prop}

\medskip
The proof will be given in Section~\ref{sec-pfs-basic-props-pb3-pb4}.

\subsection{An application to symplectic
approximation}
\label{subsec-intro-approx}

Non-vanishing of the Poisson bracket invariants can be interpreted
in terms of geometry in the space $\cF$ equipped with the uniform
distance $$d((F,G),(H,K)) = ||F-H|| + ||G-K||\;$$ as follows.
Consider the family of subsets $\cK_s \subset \cF$, $s \geq 0$,
given by $$\cK_s =\{(H,K)\in\cF\;:\; ||\{H,K\}|| \leq s \}.$$
Define {\it the profile function} $\rho_{F,G}: [0;+\infty)\to\R$
associated with a pair $(F,G) \in \cF$ (cf. \cite{EPR}) as
$$\rho_{F,G} (s) := d((F,G), \cK_s).$$
Obviously, $\rho_{F,G} (||\{F,G\}||)=0$ and the function
$\rho_{F,G} (s)$ is non-increasing and non-negative. The value
$\rho_{F,G}(0)$ is responsible for the optimal uniform
approximation of $(F,G)$ by a pair of Poisson-commuting functions.
Many results of the function theory on symplectic manifolds can be
expressed in terms of profile functions. For instance, the lower
semi-continuity of the functional $(F,G) \mapsto ||\{F,G\}||$
discussed in the beginning of this paper means that for any
$F,G\in C^\infty_c (M)$ we have that $\rho_{F,G} (s)>0$ for any
$s\in [0;||\{ F,G\}||)$. The study of the modulus of the lower
semi-continuity of this functional performed in \cite{Buh}, cf.
\cite{EP-Poisson2}, yields a sharp estimate on the convergence
rate of $\rho_{F,G} (s)$ to zero as $s\to ||\{ F,G\}||$. Below we
focus on behavior of profile functions at and near $s=0$.

Consider a triple $(X,Y,Z)$ or a quadruple $(X_0,X_1,Y_0,Y_1)$ of
compact subsets of $M$ satisfying intersection conditions
\eqref{eq-inter-3} and \eqref{eq-inter-4} respectively. In both
cases denote by $p$ the Poisson bracket invariant $pb_3(X,Y,Z)$
or, respectively, $pb_4(X_0,X_1,Y_0,Y_1)$. Define subclasses
$$\cF_3^\flat(X,Y,Z) \subset \cF_3(X,Y,Z),\;
\cF_4^\flat(X_0,X_1,Y_0,Y_1) \subset \cF_4(X_0,X_1,Y_0,Y_1)$$
consisting of all pairs $(F,G)$ such that {\it at least one of the
functions $F$,$G$} has its range in $[0;1]$. We shall often
abbreviate these classes as $\cF_3^\flat$ and $\cF_4^\flat$.

The main result of this section shows that the profile functions
associated to pairs from $\cF^{\flat}_k$ exhibit quite different
patterns of behavior depending on whether $p=0$ or $p>0$.
Furthermore, when $p
> 0$, there is a difference between the cases $k=3$ and $k=4$.

\begin{thm}
\label{thm-dichotomy}{\bf [Dichotomy]}
\begin{itemize}
\item[{(i)}] Assume that $p=0$. In this case for every $s>0$ there
exists $(F,G) \in \cF_k^{\flat}$ (where $k=3,4$)   with
$\rho_{F,G}(s) =0$.
\item[{(ii)}] Assume that $p >0$. Then for every $(F,G) \in
\cF_k^{\flat}$ (where $k=3,4$) the profile function $\rho_{F,G}$ is
continuous, $\rho_{F,G}(0) = 1/2$ and
\begin{equation}\label{eq-rho-up-trivial}
\frac{1}{2}- \frac{1}{2||\{F,G\}||}\cdot s \geq \rho_{F,G}(s) \
\forall s \in [0; ||\{F,G\}||].
\end{equation}
Furthermore,
\begin{equation}\label{eq-rho-new}
\rho_{F,G}(s) \geq \frac{1}{2}- \frac{1}{2\sqrt{p}}\cdot \sqrt{s}
\;\;\;\;\forall (F,G) \in \cF_3^{\flat} \ \forall s\geq 0,
\end{equation}
and
\begin{equation}\label{eq-rho-new-pb4}
\rho_{F,G}(s) \geq \frac{1}{2}- \frac{1}{2{p}}\cdot
{s}\;\;\;\;\forall (F,G) \in \cF_4^\flat\ \forall s\geq 0.
\end{equation}
\end{itemize}
\end{thm}

\medskip
This result, whose proof is given in
Section~\ref{subsec-sympl-approx-pfs-pb}, deserves a discussion. The
appearance of the class $\cF_k^\flat$ in our story is quite natural:
it follows from Proposition~\ref{prop-op}    that $p = \inf
||\{F,G\}||$, where the infimum is taken over all $(F,G) \in
\cF_k^\flat$. This immediately yields part (i) of the dichotomy.

A comparison between estimates \eqref{eq-rho-up-trivial} and
\eqref{eq-rho-new-pb4} shows that for $(F,G) \in \cF_4^\flat$ and
$p>0$ $$1/2-\rho_{F,G}(s) \sim s$$ for small $s$, and thus we have
captured a sharp rate, in terms of the power of $s$, of the
profile function near $0$. (Here and below we write $a(s) \sim
b(s)$ whenever for all sufficiently small $s>0$ the ratio
$a(s)/b(s)$ of non-negative functions $a$ and $b$ is bounded away
from $0$ and $+\infty$.)

In contrast to this, when $(F,G) \in \cF_3^\flat$, there is a
discrepancy in the powers of $s$ in upper bound
\eqref{eq-rho-up-trivial} and lower bound \eqref{eq-rho-new}.
Interestingly enough, for a certain triple of closed subsets
$X,Y,Z$ with a positive Poisson bracket invariant $pb_3$, both
rates $1/2-\rho_{F,G}(s) \sim s$ and $1/2-\rho_{F,G}(s) \sim
\sqrt{s}$ can be achieved by suitable pairs $(F,G) \in
\cF_3^\flat(X,Y,Z)$.

Indeed, consider the sphere
$$S^2 = \{ (x,y,z)\in\R^3\ |\ x^2 +y^2+z^2=1\}$$
with the standard symplectic form $\omega$ on it. Define $ F,G : S^2
\rightarrow \R $ by $F (x,y,z)= x^2$ and $G (x,y,z) = y^2$. These
functions lie in $\cF_3^\flat (X,Y,Z)$, where $X$, $Y$ and $Z$ are the
big circles $\{x=0\}$, $\{y=0\}$ and $\{z=0\}$ respectively. We have
seen in Example~\ref{exam-S2} that $p:=pb_3(X,Y,Z) >0$.

\begin{thm} \label{thm-S2-estimate-exact} For the functions $F,G:
S^2\to\R$ as above one has
\begin{equation}
\label{eqn-S2-estimate-exact} \rho_{F,G} (s) \leq \rho_{F,G} (0) -
C\sqrt{s}
\end{equation}
for some $C>0$.
\end{thm}

\medskip
In particular, by \eqref{eq-rho-new} we get that
$1/2-\rho_{F,G}(s) \sim \sqrt{s}$. The proof will be given in
Section~\ref{subsec-sympl-approx-pfs-ex}.

Further, cover the circle $Z$ by two open subsets, $U$ and $V$ so
that $U \cap X = V \cap Y = \emptyset$. Take any pair of non-negative functions
$F,G$ from $\cF_4^\flat(X, U \cap Z, Y, V \cap Z)$. Observe that
$(F,G)$ automatically lies in $\cF_3^\flat(X,Y,Z)$. By inequality
\eqref{eq-pb3-pb4} below, $$p:= pb_4(X, U \cap Z, Y, V \cap Z)
\geq pb_3(X,Y,Z) >0.$$ Thus by \eqref{eq-rho-new-pb4}
$$\rho_{F,G}(s) \geq \frac{1}{2}- \frac{1}{2{p}}\cdot
{s}\;,$$ and hence, by \eqref{eq-rho-up-trivial}, we get that
$1/2-\rho_{F,G}(s) \sim {s}$.

It would be interesting to explore further the rates of
$1/2-\rho_{F,G}(s)$ as $s \to 0$ for $(F,G) \in
\cF_3^\flat(X,Y,Z)$: Are there intermediate rates between $\sim s$
and $\sim \sqrt{s}$? Is there a generic rate, and if yes, what is
it?

Let us continue the discussion on the Dichotomy Theorem. The
continuity of $\rho_{F,G}(s)$ for $s>0$ holds, in fact, for {\it
any} pair $(F,G) \in \cF$ (which does not necessarily  lie in
$\cF_k^\flat$):

\begin{prop}
\label{prop-profile-fcn-Lipschitz}
 For every $(F,G) \in \cF$, the profile
function $s \rho_{F,G}$ is Lipschitz on $ (0;+\infty) $ with the
Lipschitz constant $3 \min(||F||,||G||)$.
\end{prop}

\medskip
In particular, $\rho_{F,G}$ is continuous on $ (0;+\infty)$. Let us
mention also that the Lipschitz constant of $s \rho_{F,G}(s) $ is
uniformly bounded by $3$ for all $(F,G) \in \cF_k^\flat$. The
proposition is proved in Section~\ref{subsec-sympl-approx-pfs-lip}.

The Dichotomy Theorem leaves unanswered the following natural and
closely related questions on the behavior of profile functions at
$s=0$ which, in general, are currently out of reach. The first one
deals with part (i) of the Dichotomy Theorem:

\medskip
\begin{question}\label{question-p0} Assume that the Poisson
bracket invariant $p$ vanishes. Is it true that $$\inf_{(F,G) \in
\cF_k^\flat} \rho_{F,G}(0)=0\;?$$ Or, even stronger, does there
exist a pair $(F,G)$ in $\cF_k^\flat$ or in its closure in $C(M)$
with $||\{F,G\}||=0$? In the last question we assume for simplicity that $M$ is compact,
and we define $||\{F,G\}||$ for continuous $F$ and $G$ by formula \eqref{eq-lowsem}.
\end{question}

\medskip
\noindent The second question is as follows:

\medskip
\begin{question}
\label{Q:rho-continuous-zero}
Is the function $ \rho_{F,G} $ continuous at $0$ for {\bf any} pair
of functions $(F,G) \in \cF$?
\end{question}

\medskip
\noindent It turns out that for closed manifolds of dimension two
the answers to both questions are affirmative. This readily
follows from a recent result by Zapolsky \cite{Z} which states
that every pair of functions $F,G$ on a surface with $||\{F,G\}||
\sim s$ lies at the distance $\sim \sqrt{s}$ from a
Poisson-commuting pair. In fact this yields the following more
detailed answer to Question~\ref{Q:rho-continuous-zero}, compare
with inequality \eqref{eq-rho-new}:

\begin{prop}
\label{prop-2-dim-case-almost-commut}
Suppose $(M,\omega)$ is a closed connected $2$-dimensional
symplectic manifold. For any $(F,G)$ the profile function $
\rho_{F,G} $ satisfies the inequality
\begin{equation}
\label{eqn-profile-fcn-2-dim-case} \rho_{F,G}(s) \geqslant
\rho_{F,G}(0) - C\sqrt{s},
\end{equation}
 for some constant $ C=C(M,\omega) > 0 $. In particular, $\rho_{F,G}$
 is continuous at $0$.
\end{prop}

\medskip
We refer to  Section~\ref{subsec-sympl-approx-pfs-2D} for the
proofs and further discussion.

\subsection{An application to dynamics: Hamiltonian chords}
\label{subsec-intro-chords}

\noindent
\begin{thm}
\label{thm-chords1}

Let $X_0,X_1,Y_0,Y_1 \subset M$ be a quadruple of compact sets
with $X_0 \cap X_1 = Y_0 \cap Y_1 =\emptyset$ and
$pb_4(X_0,X_1,Y_0,Y_1) = p > 0$. Let $G \in C^\infty_c (M)$ be a
Hamiltonian function with $G|_{Y_0} \leq 0$ and $G|_{Y_1} \geq 1$
generating a Hamiltonian flow $g_t$. Then $g_T x \in X_1$ for some
point $x \in X_0$ and some time moment $T \in [-1/p; 1/p]$.

\end{thm}

\medskip
We refer to the curve $\{g_t x\}_{t \in [0;T]}$ as to {\it a
Hamiltonian chord} of $g_t$ (or, for brevity, of the Hamiltonian
$G$) of time-length $|T|$ connecting $X_0$ and $X_1$.

Hamiltonian chords joining two disjoint subsets (notably, Lagrangian
submanifolds) of a symplectic manifold arise in several interesting
contexts such as Arnold diffusion (see e.g. \cite{KL},\cite[Question
0.1]{Bernard}) or optimal control (see e.g. \cite{Pontryagin et al},
\cite[Ch.12]{Nijmeijer-vdSchaft}, \cite{Ivan-applied}). Furthermore,
Hamiltonian chords had been studied on various occasions in
symplectic topology (see e.g. \cite{AF,Merry}).

Theorem~\ref{thm-chords1} has a flavor of the following well-known
phenomenon in symplectic dynamics: For a suitably chosen pair of
subsets $Y_0$ and $Y_1$ of a symplectic manifold the condition
$\min_{Y_1} F - \max_{Y_0} F \geq C$ yields existence of a periodic
orbit of the Hamiltonian flow of $F$ with some interesting
properties provided $C$ is large enough, see \cite{Hofer-Zehnder,
Gatien-Lalonde, BPS}. Theorem~\ref{thm-chords1} extends this
phenomenon to the case of non-closed orbits, i.e. Hamiltonian
chords.

It turns out that the bound on the time-length of a Hamiltonian
chord given in Theorem~\ref{thm-chords1} is sharp in the following
sense. Given two disjoint compact subsets $X_0,X_1$ of $M$ and a
Hamiltonian $G\in C^\infty_c(M)$, denote by $T(X_0,X_1;G)$ the
minimal time-length of a Hamiltonian chord of $G$ which connects
$X_0$ and $X_1$. (Here we set $\inf \emptyset := +\infty$.) Put
$$T(X_0,X_1,Y_0,Y_1)= \sup\, \{\ T(X_0,X_1;G) : G\in C^\infty_c (M),\ G|_{Y_0} \leq 0,\
G|_{Y_1} \geq 1\ \}.$$

\begin{thm}\label{thm-eq-pb-T}
$$pb_4(X_0,X_1,Y_0,Y_1) = T(X_0,X_1,Y_0,Y_1)^{-1}.$$
\end{thm}

\medskip
The proof will be given in Section~\ref{sec-chords-proofs}. This
result can be considered as a dynamical interpretation of the
invariant $pb_4$. It immediately yields Theorem~\ref{thm-chords1}.

Let us pass to the case of Hamiltonian chords for non-autonomous
flows. We shall need the following notion.

\medskip
\noindent{\sc Stabilization:} Identify the cotangent bundle
$T^*S^1$ with the cylinder $\R \times S^1$ equipped with the
coordinates $r$ and $\theta\, ({\rm mod}\, 1)$ and the standard
symplectic form $dr \wedge d\theta$. Denote by $\A_R \subset
T^*S^1$, $0< R \leq \infty $, the annulus $\{|r| < R\}$. Given a
compact subset $X$ of a symplectic manifold $(M,\omega)$, define
its {\it $R$-stabilization}
$$\s _R X:= X \times S^1 \subset (M \times \A_R, \omega + dr \wedge
d\theta),$$ where $S^1$ is identified with the zero section
$\{r=0\}$. We shall abbreviate $\s X$ for $\s_\infty X$.

\begin{thm}\label{thm-non-autonomous}
 Let $X_0,X_1,Y_0,Y_1 \subset M$ be a quadruple of compact sets
with $X_0 \cap X_1 = Y_0 \cap Y_1 =\emptyset$ and $$pb_4(\s_R
X_0,\s_R X_1,\s_R Y_0,\s_R Y_1) = p > 0$$ for some $R \in
(1;+\infty]$. Let $G \in C^\infty_c(M \times S^1)$ be a
(non-autonomous) $1$-periodic Hamiltonian with $G_t|_{Y_0} \leq
0$, $G_t|_{Y_1} \geq 1$ for all $t\in S^1$ and
\begin{equation}\label{eq-oscillation}
\max G - \min G < R\;
\end{equation}
generating a Hamiltonian flow $g_t$. Then there exists a point $x
\in M$ and time moments $t_0,t_1 \in \R$ with $|t_0-t_1| \leq 1/p$
such that $g_{t_0} x \in X_0$ and $g_{t_1} x \in X_1$.
\end{thm}

\medskip
The proof will be given in Section~\ref{sec-chords-proofs}. Exactly
as in the autonomous case, the curve $\{g_t x\}$, $t \in [t_0;t_1]$,
is called a {\it Hamiltonian chord} passing through $X_0$ and $X_1$.
We refer to Remark \ref{rem-comparison-time} below for a comparison
of the  bounds on the time-length of Hamiltonian chords given by
Theorem~\ref{thm-chords1} (the autonomous case) and
Theorem~\ref{thm-non-autonomous} (the non-autonomous case).

Here is a sample application of our theory. Consider a compact
domain $V \subset T^*\T^n$ whose interior contains the zero
section $\T^n$. Fix a pair of distinct points $q_0,q_1 \in \T^n$
and put $D_i = T^*_{q_i}\T^n \cap V$.

\begin{thm}\label{thm-torus-chords} Let $G: V \times S^1 \to \R$
be a Hamiltonian which vanishes near $\partial V \times S^1$ and
which is $\geq 1$ on $\T^n \times S^1$. Then there exists a
Hamiltonian chord of the Hamiltonian flow of $G$ passing through
$D_0$ and $D_1$. \end{thm}

\medskip
The proof is given in Section~\ref{subseq-pb-qs} below. As an
illustration, assume that the torus $\T^n$ is equipped with a
Riemannian metric and $V=\{|p| \leq 1\}$, where $(p,q)$ are
canonical coordinates on $T^*\T^n$. Suppose that the Hamiltonian
$G$ has the form $u(|p|)$, where $u(s)$ vanishes for $s$ close to
$1$ and $u(0) =1$. Then the projection of the Hamiltonian chord
provided by Theorem~\ref{thm-torus-chords} is a (reparameterized)
Riemannian geodesic segment joining the points $q_0$ and $q_1$.
Theorem~\ref{thm-torus-chords} resembles the one of \cite{BPS}
where under similar assumptions the authors proved the existence
of closed trajectories imitating closed geodesics on the torus.
The approach of \cite{BPS} was based on relative symplectic
homology. It would be interesting to find its footprints in our
context. It would be also interesting to compare our approach with
the one of Merry \cite{Merry} who detects Hamiltonian chords by
using a Lagrangian version of Rabinowitz Floer homology.

\subsection{Poisson bracket invariants and symplectic
quasi-states}\label{subseq-pb-qs}

Now we turn to a discussion of methods for establishing lower
bounds (and, in particular, the positivity) for the Poisson
bracket invariants of certain triples and quadruples of compact
subsets of a symplectic manifold. The first method is based on the
theory of symplectic quasi-states and quasi-measures.

\noindent
{\bf In this
section we assume that $(M^{2n},\omega)$ is a closed connected
symplectic manifold.}

Denote by $C(M)$ the space of the continuous
functions on $M$.  A {\it symplectic quasi-state} \cite{EP-qst} is
a functional $\zeta: C(M) \to \R$ which satisfies the following
axioms:

\medskip
\noindent {\bf (Normalization)} $\zeta(1)=1 $;

\medskip
\noindent {\bf (Positivity)} $\zeta(F) \geq 0$ provided $F \geq 0$;

\medskip
\noindent {\bf (Quasi-linearity)} $\zeta$ is linear on every
Poisson-commutative subspace of $C(M)$.

\medskip
\noindent Here we say that two continuous functions $F,G \in C(M)$
Poisson-commute if there exist sequences of smooth functions
$\{F_i\}$ and $\{G_i\}$ which uniformly converge to $F$ and $G$
respectively so that $||\{F_i,G_i\}||\to 0$ as $i \to +\infty$.
This notion is well-defined due to the $C^0$-robustness of the
Poisson bracket, see Section~\ref{subsec-robustness} above.

Recall that a {\it quasi-measure} associated to a quasi-state
$\zeta$ is a set-function whose value on a closed subset $X$
equals, roughly speaking, $\zeta(\chi_X)$, where $\chi_X$ is the
indicator function of $X$ (see e.g. \cite{EP-qst}). A closed
subset $X \subset M$ is called {\it superheavy} with respect to
$\zeta$ if $\tau(X)=1$. Equivalently, $X$ is superheavy whenever
$\zeta(F) \geq c$ for any $F$ with $F|_X \geq c$, and hence
automatically $\zeta(F)\leq c$ for any $F$ with $F|_X \leq c$,
see \cite{EP-rigid}.

We say that a symplectic quasi-state $\zeta$ satisfies the {\it
PB-inequality} (with ``PB" standing for the ``Poisson brackets"),
if there exists $K=K(M,\omega)> 0$ so that
\begin{equation}\label{eq-frol}
|\zeta(F+G)-\zeta(F)-\zeta(G)| \leq \sqrt{K||\{F,G\}||} \;\;\forall
F,G \in C(M).
\end{equation}
Here $||\{F,G\}||$ for continuous functions $F,G \in C(M)$ is
understood in the sense of \eqref{eq-lowsem}.

At present we know a variety of examples of symplectic manifolds
admitting symplectic quasi-states which satisfy the PB-inequality,
as well as plenty of examples of superheavy subsets \cite{EP-rigid},
\cite{EPZ}, \cite{Ostrover}, \cite{Usher1}, \cite{Usher2}.

\medskip
\begin{exam}\label{exam-superheavy} {\rm The complex projective space $\C
P^n$ equipped with the standard Fubini-Study symplectic form
admits a symplectic quasi-state satisfying PB-inequality. Its
superheavy subsets include certain monotone Lagrangian
submanifolds such as the Clifford torus and the real projective
space $\R P^n$, as well as certain singular subsets such as a
codimension-1 skeleton of a sufficiently fine triangulation. Any
product of $\C P^n$'s with the split symplectic form also admits
such a quasi-state, and the product of superheavy sets is again
superheavy.}
\end{exam}

\medskip
For $M$ of dimension higher than $2$ the only currently known construction of such quasi-states is
based on the Hamiltonian Floer theory and works under the
assumption that the quantum homology algebra $QH_* (M)$ of $M$
splits as an algebra into a direct sum so that one of the summands
is a field (see \cite{EP-toric}, \cite{Usher2}). Such quasi-states
automatically satisfy the PB-inequality (see \cite{EPZ}).

\begin{thm}\label{thm-qstates-pb4-dynamics}
Assume that a closed symplectic manifold $(M,\omega)$ admits a
symplectic quasi-state $\zeta$ which satisfies PB-inequality
\eqref{eq-frol} with a constant $K$.
\begin{itemize}
\item[{(i)}] Let $X,Y,Z\subset M$ be a triple of superheavy
closed sets with $X\cap Y\cap Z=\emptyset$. Then
\begin{equation}
\label{eq-Delta} pb_3 (X,Y,Z) \geq \frac{1}{K}.
\end{equation}
\item[{(ii)}] Let $X_0,X_1,Y_0,Y_1\subset M$ be a quadruple of
closed subsets such that
$$X_0 \cap X_1 = Y_0 \cap Y_1 =\emptyset.$$

If $X_0 \cup Y_0, Y_0 \cup X_1, X_1 \cup Y_1, Y_1 \cup X_0$ are all
superheavy, then
\begin{equation}
\label{eq-Delta-1} pb_4(X_0,X_1,Y_0,Y_1) \geq \frac{1}{4K}.
\end{equation}

If $X_0 \cup Y_0, Y_0 \cup X_1, Y_1$ are all
superheavy (this condition is stronger than the previous one), then
\begin{equation}
\label{eq-Delta-2} pb_4(X_0,X_1,Y_0,Y_1) \geq \frac{1}{K}.
\end{equation}

\end{itemize}
\end{thm}

\medskip\noindent{\bf Proof.} The theorem follows immediately from the
formalism described above (cf. \cite{EPZ}, Theorem 1.7): To prove
(i), assume that $F|_X \leq 0$, $G|_Y \leq 0$, $F+G|_Z \geq 1$. By
the superheaviness, $\zeta(F) \leq 0$, $\zeta(G) \leq 0$ and
$\zeta(F+G) \geq 1$. Applying PB-inequality \eqref{eq-frol} we get
\eqref{eq-Delta}.

Let us pass to the proof of (ii).
Assume
$X_0 \cup Y_0, Y_0 \cup X_1, X_1 \cup Y_1, Y_1 \cup X_0$ are all
superheavy.
By Proposition~\ref{prop-op}(ii), it suffices to find a lower bound
on $a:=||\{F,G\}||$ for pairs $(F,G)\in \cF'_4 (X_0,X_1,Y_0,Y_1)$. Put
$$u_1=FG,\; u_2 = G(1-F),\; u_3=(1-F)(1-G),\;u_4 = F(1-G).$$
These functions vanish on the superheavy sets $X_0 \cup Y_0, Y_0 \cup
X_1, X_1 \cup Y_1, Y_1 \cup X_0$ respectively and hence
$\zeta(u_i)=0$ for all $i$. Also note that $\sum_i u_i=1$.
On the other hand,
%as one can easily check,
%$$\max_{i,j} ||\{u_i,u_j\}|| \leq 2a.$$
$$||\{u_2,u_3\}|| = ||\{G(1-F),(1-G)(1-F)\}|| =||(1-F)\{F,G\}||\leq a,$$
$$||\{u_2+u_3,u_4\}|| = ||\{ 1-F,F(1-G)\}|| = ||F\{ F,G\}||\leq a.$$
Together with PB-inequality \eqref{eq-frol} this yields
$$|\zeta(u_2+u_3) |=|\zeta (u_2 +u_3) - \zeta(u_2) -\zeta(u_3)|\leq \sqrt{K}\sqrt{||\{u_2,u_3\}||}\leq \sqrt{aK},$$
$$|\zeta(u_2+u_3+u_4) - \zeta(u_2+u_3) - \zeta (u_4)|\leq \sqrt{K} \sqrt{||\{u_2+u_3, u_4\}||} \leq \sqrt{aK}.$$
Using the equality $\{u_1,u_2+u_3+u_4\}= \{u_1,1-u_1\}=0$ we get
$$1 = |\zeta (u_1+(u_2+u_3+u_4))| =  |\zeta(u_2+u_3+u_4)|\leq $$
$$\leq  |\zeta (u_2+u_3)|  + \sqrt{K||\{u_2+u_3,u_4\}||}\leq 2\sqrt{aK}.$$
Thus $a\geq 1/(4K)$ which proves \eqref{eq-Delta-1}.

Now assume
$X_0 \cup Y_0, Y_0 \cup X_1, Y_1$ are all
superheavy. Then part (i), together with inequality \eqref{eqn-pb3-pb4-max} below comparing $pb_3$ and $pb_4$, imply
$$pb_4(X_0,X_1,Y_0,Y_1)\geq  pb_3 (X_0\cup Y_0, X_1\cup Y_0, Y_1)\geq 1/K,$$
that is  \eqref{eq-Delta-2}.
\Qed

%and the equalities $\{u_1, u_2+u_3+u_4\} = 0$ and $\zeta (u_i)=0$ this yields
%$$|\zeta(u_2+u_3+u_4) - \zeta(u_2+u_3) - \zeta (u_4)|\leq \sqrt{Ka}
%$$1 = |\zeta (\sum u_i) -\sum \zeta (u_i)| = |\zeta(u_2+u_3+u_4) - \zeta (u_2) - \zeta ( \leq  4 \cdot \sqrt{Ka}$$
%(see e.g. \cite{EPZ}, p. 1044). Thus $||\{F,G\}|| = a/3 \geq 0.01
%\cdot K^{-1}$ as required. \Qed

\medskip
\begin{exam}\label{exam-S2-superheavy}{\rm A big circle of $S^2$
(or, in other words, a Clifford torus of $\C P^1$) is superheavy.
This yields the positivity of $pb_3$ in Example~\ref{exam-S2} above.}
\end{exam}

\medskip
Let us discuss some applications of
Theorem~\ref{thm-qstates-pb4-dynamics} to the existence of
Hamiltonian chords. In order to formulate them we need the following
notion. Consider the sphere $S^2$ equipped with an area form
$\sigma$ of total area $1$. Denote by $E$ the equator of $S^2$. Let
$\zeta$ be a symplectic quasi-state on $M$ satisfying PB-inequality.
We say that $\zeta$ is {\it $S^2$-stable} if for every $c > 0$ the
symplectic manifold $(M \times S^2, \omega + c \sigma)$ admits a
symplectic quasi-state $\tzeta_c$ which satisfies PB-inequality and
such that $Z \times E$ is $\tzeta_c$-superheavy for every superheavy
subset $Z \subset M$. The quasi-states associated to field factors
of quantum homology are known to be $S^2$-stable \cite{EP-rigid}.
{\bf In part (ii) of the next corollary superheaviness is considered
with respect to a $S^2$-stable quasi-state on $(M,\omega)$.}

\begin{cor}\label{cor-non-autonomous-superheavy}
 Let $X_0,X_1,Y_0,Y_1 \subset M$ be a quadruple of compact sets
 such that
$X_0 \cap X_1 = Y_0 \cap Y_1 =\emptyset$ and the sets $X_0 \cup Y_0,
Y_0 \cup X_1, X_1 \cup Y_1, Y_1 \cup X_0$ are all superheavy. Let $G
\in C^\infty_c(M \times S^1)$ be a $1$-periodic
Hamiltonian with $G_t|_{Y_0} \leq 0$, $G_t|_{Y_1} \geq 1$ for all
$t\in S^1$. Then there exists a point $x \in M$ and time moments
$t_0,t_1 \in \R$ so that $g_{t_0} x \in X_0$ and $g_{t_1} x \in
X_1$. Furthermore,
\begin{itemize}
\item[{(i)}] If $G$ is autonomous, $|t_0-t_1| \leq 4K$.
If in addition $Y_1$ is super-heavy, $|t_0-t_1| \leq K$.
\item[{(ii)}] If $G$ is non-autonomous,
$|t_0 -t_1| \leq C$,  where the constant $C$
depends only on the symplectic quasi-state $\zeta$ and
on the oscillation $\max G - \min G$ of the Hamiltonian $G$.
\end{itemize}
\end{cor}

Part (i) is an immediate consequence of
Theorem~\ref{thm-qstates-pb4-dynamics}(ii) combined  with
Theorem~\ref{thm-chords1}. Part (ii) can be deduced from
Theorems~\ref{thm-qstates-pb4-dynamics}(ii) and
\ref{thm-non-autonomous}, see Section~\ref{sec-chords-proofs}
below for the proof and for more information on $C$.
%
%\begin{cor}\label{cor-non-autonomous-superheavy}
% Let $X_0,X_1,Y_0,Y_1 \subset M$ be a quadruple of compact sets
% such that
%$X_0 \cap X_1 = Y_0 \cap Y_1 =\emptyset$.
% Let $G
%\in C^\infty_c(M \times S^1)$ be a $1$-periodic
%Hamiltonian with $G_t|_{Y_0} \leq 0$, $G_t|_{Y_1} \geq 1$ for all
%$t\in S^1$.
%
%\bigskip
%A. Assume the sets $X_0 \cup Y_0,
%Y_0 \cup X_1, X_1 \cup Y_1, Y_1 \cup X_0$ are all superheavy. Then there exists a point $x \in M$ and time moments
%$t_0,t_1 \in \R$ so that $g_{t_0} x \in X_0$ and $g_{t_1} x \in
%X_1$. Furthermore,
%\begin{itemize}
%\item[{(A1)}] If $G$ is autonomous, $|t_0-t_1| \leq 4K$;
%\item[{(A2)}] If $G$ is non-autonomous,
%$|t_0 -t_1| \leq C$,  where the constant $C$
%depends only on the symplectic quasi-state $\zeta$ and
%on the oscillation $\max G - \min G$ of the Hamiltonian $G$.
%\end{itemize}
%
%\bigskip
%B. Assume  $X_0 \cup Y_0,
%X_1 \cup Y_0, Y_1$ are superheavy (this condition is stronger than the one in part A). Then the same results as in part A are true with the constants $4K$ and $C$
%being replaced, respectively, by $K$ and $C/4$ .
%
%\end{cor}
%
%Part (A1) is an immediate consequence of
%Theorem~\ref{thm-qstates-pb4-dynamics}(ii) combined  with
%Theorem~\ref{thm-chords1}. Part (A2) can be deduced from
%Theorems~\ref{thm-qstates-pb4-dynamics}(ii) and
%\ref{thm-non-autonomous}, while Part B follows from Part A and inequality \eqref{eqn-pb3-pb4-max} comparing $pb_3$ and $pb_4$ -- see Section~\ref{sec-chords-proofs}
%below for the proofs and for more information on $C$.

\medskip
\begin{exam}\label{exam-S2xS2-qstates} {\rm Let $M= S^2 \times \ldots \times S^2$ be the
product of $n$ copies of $S^2$ equipped with the split symplectic
structure $\omega= \sigma \oplus\ldots\oplus \sigma$, where $\int_{S^2}\sigma
= 1$. Denote by $(x_i,y_i,z_i)$ the Euclidean coordinates and by
$(z_i, \phi_i)$ the cylindrical coordinates on the $i$-th copy of
the sphere ($i=1,2$), where $\phi_i$ is the polar angle in the
$(x_i,y_i)$-plane. Define the following subsets in the $i$-th
factor:  Fix $a \in (0;1/2)$ so that the $\sigma$-area of the set
$B_i = \{|z_i| \geq a\}$ is greater than $1/2$. Define an annulus
$A_i = \{|z_i| \leq a\}$. Write $E_i$ for the equator $\{z_i=0\}$
and $C_i^\theta$ for the segment
$$\{\phi_i = \theta\} \cap A_i.$$
Define the following subsets of $M$:
$$Y_0 = M \setminus \prod \text{Interior}(A_i),\;\;\;Y_1= \prod E_i .$$  For $ v = (\theta_1,\ldots,\theta_n) \in
\R^n/2\pi\Z^n$ denote $C^v = \prod C_i^{\theta_i}\;$. Put  $X_0 =
C^v$, $X_1 = C^w$, where $v,w$ are two distinct points in
$\R^n/2\pi\Z^n$.

We claim that the quadruple $X_0,Y_0,X_1,Y_1$ satisfies the
assumptions of Theorem~\ref{thm-qstates-pb4-dynamics}. The
argument uses basic criteria of superheaviness for which we refer
to \cite{EP-rigid}. The set $Y_1$ is the Clifford torus in $M$ and
thus superheavy.
Let us check that $X_0 \cup Y_0$ is superheavy. Note that
$$\prod (B_i \cup C_i^{\theta_i}) \subset Y_0 \cup X_0.$$
But $B_i \cup C_i^\theta$ is the complement to an open disc of area
$< 1/2$ and hence superheavy in $S^2$. Since the product of
superheavy sets is again superheavy, we conclude that $X_0 \cup Y_0$
is superheavy. Analogously, $X_1 \cup Y_0$ is superheavy. Thus, $X_0\cup Y_0$, $X_1\cup Y_0$ and $Y_1$ are superheavy and therefore,
by Theorem~\ref{thm-qstates-pb4-dynamics}(ii),
$pb_4(X_0,X_1,Y_0,Y_1)\geq  1/K$.
}

\end{exam}

\medskip
As we shall see right now, Theorem~\ref{thm-torus-chords} can be
easily reduced to the situation analyzed in the previous example.

\medskip
\noindent {\bf Proof of Theorem~\ref{thm-torus-chords}.} We use
the notations of Example~\ref{exam-S2xS2-qstates}. Identify the
interior of $\prod A_i$ with a neighborhood $W$ of the zero
section in $T^*\T^n$ so that the zero section corresponds to the
Lagrangian torus $Y_1$ and every cotangent fiber intersects $W$
along the cube $C^u \setminus \partial C^u$ for some $u \in \T^n$.
Making, if necessary, the rescaling $(p,q) \to (\mu p,q)$ with a
sufficiently small $\mu >0$, we can assume that the domain $V$
from the formulation of the theorem is contained in $W$. Then the
sets $D_0$ and $D_1$ are identified with $X_0 \cap V$ and $X_1
\cap V$ respectively.

Let $Z$ be the closure of $M \setminus V$. Observe that $Z \cup D_i$
contains $Y_0 \cup X_i$ and hence is superheavy.

Take any function $G: V \times S^1 \to \R$ which vanishes near
$\partial V \times S^1$ and is $\geq 1$ on $Y_1 \times S^1$. Extend
it by zero to the whole $M \times S^1$. By
Corollary~\ref{cor-non-autonomous-superheavy} applied to the
quadruple $(D_0,D_1,Z,Y_1)$, the Hamiltonian flow $g_t$ of $G$ has a
chord passing through $D_0$ and $D_1$. Since $g_t$ is the identity
outside $V$, this chord is entirely contained in $V$. The
time-length of this chord admits an upper bound provided by
Corollary~\ref{cor-non-autonomous-superheavy}. \Qed

\subsection{Poisson bracket and deformations of the symplectic
form}\label{subsec-deformation}

In this section we present yet another approach to the positivity
of the Poisson bracket invariants which is applicable to certain
triples and quadruples of (sometimes singular) Lagrangian
submanifolds. Our method is based on a special deformation of the
symplectic form on $M^{2n}$ combined with the study of
``persistent" pseudo-holomorphic curves with Lagrangian boundary
conditions (cf. \cite{Akveld-Salamon}).

\subsubsection{A lower bound}\label{subsubsec-def-lowerbd}

Given two functions $F,G\in C^\infty_c (M)$, consider the family of
forms
$$\omega_s := \omega -sdF\wedge dG.$$
Note that $$dF \wedge dG \wedge \omega^{n-1} =
\frac{1}{n}\{F,G\}\cdot \omega^n.$$ Thus
$$\omega_s^n =(1-s\{F,G\})\omega^n.$$ Therefore
the form $\omega_s$ is symplectic for all
$$s \in I:= [0; 1/||\{F,G\}||).$$
(We set $1/||\{F,G\}|| = +\infty$ if $\{ F,G\}\equiv 0$.)

Recall that an almost complex structure $J$ on $M$ is said to be
compatible with $\omega$ if $\omega (\xi, J \eta)$ is a Riemannian
metric on $M$. Choose a generic family of almost complex
structures $J_s$, $s \in I$, compatible with $\omega_s$.

The next elementary proposition allows to relate Poisson brackets to
pseudo-holomorphic curves:

\begin{prop}
\label{prop-psh-curves-general} Let $F,G\in C^\infty_c (M)$.
Assume that there exist
\begin{itemize}
\item{} a family of almost complex structures $J_s$, $s\in I$,
such
that each $J_s$ is compatible with the symplectic form $\omega_s =
\omega -sdF\wedge dG$,

\item{} a family of $J_s$-holomorphic maps $u_s: \Sigma_s\to M$,
$s\in I$, where each $\Sigma_s$ is a compact Riemann surface with
boundary and possibly with corners,

\item{} positive constants $C_1, C_2$,

\end{itemize}
so that for all $s\in I$
\begin{equation}
\label{eqn-int-omega-psh-curve} \int_{\Sigma_s} u_s^* \omega \leq
C_1
\end{equation}
and
\begin{equation}
\label{eqn-int-fdg-boundary} \int_{\partial \Sigma_s} u_s^* (FdG)
\geq C_2.
\end{equation}
Then $||\{ F,G\}|| \geq C_2/C_1$.

\end{prop}

\medskip\noindent{\bf Proof.}
Applying the Stokes theorem together with
\eqref{eqn-int-omega-psh-curve} and \eqref{eqn-int-fdg-boundary} we
get
$$0\leq \int_{\Sigma_s} u_s^* \omega_s = \int_{\Sigma_s} u_s^* \omega
- s\int_{\partial\Sigma_s} u_s^* (FdG) \leq C_1 -
s\int_{\partial\Sigma_s} u_s^* (FdG).$$ Hence
$$C_2 s\leq s\int_{\partial\Sigma_s} u_s^* (FdG) \leq C_1$$
and thus $C_2 s \leq C_1$ for any $s\in I=[0; 1/||\{F,G\}||)$. Note
that $||\{F,G\}||\neq 0$ (since $C_2$ is assumed to be positive) and
therefore $C_2/||\{F,G\}||\leq C_1$ and thus $C_2/C_1\leq
||\{F,G\}||$. \Qed

\medskip
\noindent We always apply
Proposition~\ref{prop-psh-curves-general} in the following
situation. First, assume the pair of functions $(F,G)$ lies in
$\cF'_3(X,Y,Z)$ (respectively in $\cF'_4(X_0,X_1,Y_0,Y_1)$),
where $X \cap Y \cap Z=\emptyset$ (respectively $X_0 \cap X_1 =
Y_0 \cap Y_1 = \emptyset$.) Put $W= X \cup Y \cup Z$ (resp. $W=
X_0 \cup X_1 \cup Y_0 \cup Y_1$). Observe that the $1$-form $FdG$
is necessarily closed in a sufficiently small neighborhood $U$ of
$W$. Moreover, the image of $[FdG]$ in $H^1(W,\R)$ under the
natural morphism $H^1(U,\R) \to H^1(W,\R)$ does not depend on the
specific choice of $(F,G)$. Second, assume that the boundaries of
the curves $u(\Sigma_s)$ lie on $W$. In view of this discussion,
$\int_{\partial \Sigma_s} u_s^* (FdG)$ is fully determined by the
homology class of $u_s(\partial \Sigma_s)$ in $H_1(W,\Z)$.
Similarly, $\int_{\Sigma_s} u_s^* \omega$ is determined by the
relative homology class of $u_s(\Sigma_s)$ in $H_2(M,W,\Z)$. The
conclusion of this discussion is that under these two assumptions
inequalities \eqref{eqn-int-omega-psh-curve} and
\eqref{eqn-int-fdg-boundary} have purely topological nature, and
hence can be easily verified.

We will now discuss various specific cases where
Proposition~\ref{prop-psh-curves-general} can be applied.
Before moving further let us illustrate our main idea by the following
elementary example which does not involve any advanced machinery.

\subsubsection{Case study: quadrilaterals on
surfaces}\label{exam-pb3-surfaces}

Let $(M,\omega)$ be a  symplectic surface of area $B \in
(0;+\infty]$.    Consider a curvilinear quadrilateral $\Pi\subset
M$ of area $A$ with sides denoted in the cyclic order by $a_1,
a_2, a_3, a_4$  -- that is $\Pi$ is a topological disc bounded by
the union of four smooth embedded curves $a_1, a_2, a_3, a_4$
connecting four distinct points in $M$ in the cyclic order as
listed here and (transversally) intersecting each other only at
their common end-points. Our objective is to calculate/estimate
the value of $pb_4(a_1,a_3,a_2,a_4)$. Recall from
Section~\ref{subsec-intro-intro} that
\begin{equation}
\label{eq-or} pb_4(a_1,a_3,a_2,a_4)= pb_4(a_1,a_3,a_4,a_2).
\end{equation}
Thus without loss of generality we can assume that the orientation
of $\partial \Pi$ induced by the cyclic order of $a_i$'s coincides
with the boundary orientation.

\medskip
\noindent
\begin{thm}\label{thm-surface-1}
$pb_4(a_1,a_3,a_2,a_4)= \max(1/A, 1/(B-A))$.
\end{thm}

\medskip
\noindent {\bf Proof.} $\;$

\medskip
\noindent{\sc Lower bound:} Pick any $(F,G)\in \cF'_4(a_1,a_3,
a_2, a_4)$. Note that the quadrilateral $\Pi$ is $J$-holomorphic
for {\it any} (almost) complex structure $J$ on $M$ compatible
with the orientation. Also note that, by a direct calculation,
$\int_{\partial \Pi} FdG = 1$. Thus one can apply
Proposition~\ref{prop-psh-curves-general} with $\Sigma=\Pi$ and
get that
$$||\{ F,G\}||\geq 1/A.$$
Since this is true for any $(F,G)\in \cF'_4(a_1,a_3, a_2, a_4)$, we
get that
\begin{equation}\label{eq-A-bound}
pb_4(a_1,a_3, a_2, a_4) \geq 1/A.
\end{equation}

Further, if $M$ is a closed surface apply
Proposition~\ref{prop-psh-curves-general} with $\Sigma=\overline{M
\setminus \Pi}$. We get that (mind the order of sides)
\begin{equation}\label{eq-B-A}
pb_4(a_1,a_3,a_4,a_2) \geq 1/(B-A).
\end{equation}

If $M$ is open, the surface $\Sigma$ is not compact. However,
since $\Sigma$ is properly embedded and the functions $F$ and $G$
are compactly supported, Proposition~\ref{prop-psh-curves-general}
is still applicable (after an obvious modification) and yields
inequality \eqref{eq-B-A}.

Combining inequalities \eqref{eq-A-bound} and \eqref{eq-B-A} with
\eqref{eq-or} we get that
\begin{equation}\label{eq-quadr-low}
pb_4(a_1,a_3,a_2,a_4) \geq \max(1/A, 1/(B-A)).
\end{equation}

\bigskip \noindent{\sc Upper bound:} Put $\alpha = \sqrt{A}$ and
choose any $\beta \in (\alpha;\sqrt{B})$. By Moser's theorem
\cite{Moser}, we can assume that for $\epsilon >0$ small enough
$M$ contains a square $K= [-\epsilon;\beta+\epsilon]^2$ equipped
with coordinates $(p,q)$ so that the symplectic form $\omega$ is
given by $dp \wedge dq$ and the quadrilateral $\Pi$ is given by
$[0;\alpha]^2$. Define a piece-wise linear function $u(t)$ so that
$u(t)=0$ for $t <0$ and $t
> \beta$, $u(t) = t/\alpha$ for $t \in [0;\alpha]$ and $u(t) =
(\beta-t)/(\beta-\alpha)$ for $t \in [\alpha;\beta]$. For $\delta
>0$ denote by $u_\delta$ a smoothing of $u$ with $u_\delta=0$ outside $(0;\beta)$,
$u_\delta(\alpha)=1$ and
$$|u'_\delta(t)| \leq \gamma:= \max(1/{\alpha},
{1}/({\beta-\alpha}))+\delta.$$ Take any cut-off function $v$ on
$K$ which is supported in the interior of $K$ and equals $1$ on
$[0;\beta]^2$. Consider the functions $F:= v(p,q)u_\delta (p)$ and
$G=v(p,q)u_\delta (q)$ which we extend by $0$ to the whole $M$.
Note that (after an appropriate labelling of the sides of $\Pi$)
$(F,G) \in \cF_4(a_1,a_3,a_2,a_4)$ and a straightforward
calculation shows that $||\{F,G\}|| \leq \gamma^2$. Since such $F$
and $G$ exist for all $\beta, \delta$, we get that
$pb_4(a_1,a_3,a_2,a_4) \leq \max(1/A, 1/(B-A))$. Together with
\eqref{eq-quadr-low}, this yields the theorem. \qed

\medskip

Let us discuss now what happens with the $pb_4$-invariant for
stabilizations of the sets $a_1,a_2,a_3,a_4$. Interestingly enough,
the situation is quite subtle. Suppose that $M \neq S^2$.  Let $K$
be any exact section of $T^*S^1$. We claim that
\begin{equation}\label{eq-stab-surf-1}
pb_4(a_1 \times K, a_3 \times K, a_2 \times K,  a_4 \times K) \geq
1/A.
\end{equation}

Indeed, after a $C^0$-perturbation $a_1^\epsilon, a_2^\epsilon,
a_3^\epsilon, a_4^\epsilon$ of $a_1, a_2, a_3, a_4$ we can assume
that
$$ L = a_1^\epsilon \cup a_2^\epsilon \cup a_3^\epsilon \cup
a_4^\epsilon$$ is a smooth embedded circle in $M$ enclosing area
$A^\epsilon$. Take a split complex structure $J$ on $M \times
T^*S^1$. Observe that $\widehat{L}= L \times K$ is a Lagrangian
torus in $M \times T^*S^1$. As the deformation parameter $s$
changes, $\widehat{L}$ remains Lagrangian for the deformed
symplectic structure $\omega_s = \omega-sdF \wedge dG$, provided
$(F,G) \in \cF'_4(a_1^\epsilon, a_3^\epsilon, a_2^\epsilon,
a_4^\epsilon)$, but its symplectic area class alters. (By Moser's
theorem \cite{Moser}, an equivalent viewpoint is that the
symplectic form on $M \times T^*S^1$ is fixed, but $\widehat{L}$
undergoes the process of a non-exact Lagrangian isotopy.) The
class $\alpha:= [\Pi \times \text{point}]$ is the generator of
$\pi_2(M \times T^*S^1, L \times K)$. Thus the standard Gromov's
theory \cite{Gromov-pshc} yields that for a generic deformation
$J_s$ of $J$ as in Proposition~\ref{prop-psh-curves-general} there
exists a pseudo-holomorphic disc $\Sigma_s$ in the class $\alpha$
(this argument breaks down for $M=S^2$ due to possible bubbling).
Therefore Proposition~\ref{prop-psh-curves-general} with $C_1 =
A^\epsilon$ and $C_2 =1$ (the latter readily follows from the
Stokes theorem) yields
$$pb_4(a_1^\epsilon \times K, a_3^\epsilon
\times K, a_2^\epsilon \times K,  a_4^\epsilon \times K) \geq
1/A^\epsilon.$$ Passing to the limit as the size of perturbation
$\epsilon$ goes to $0$ (this procedure is justified in
Proposition~\ref{prop-Hausdorff} below) we get inequality
\eqref{eq-stab-surf-1}.

Let us emphasize that the Gromov-theoretical argument as above
does not work for surfaces  $\Sigma$ other than discs, and, in
particular, it is not applicable to $\overline{M \setminus \Pi}$.
Thus we are unable to find the lower bound for $pb_4$ in terms of
$1/(B-A)$ as it was done in the proof of
Theorem~\ref{thm-surface-1} in the two-dimensional case. Therefore
in general we do not know the exact value of $pb_4$ in this
situation. However, we have the following partial result.

\medskip
\noindent \begin{prop}\label{thm-surfaces-2} Assume that $M \neq
S^2$ and $2A \leq B$. Then
$$pb_4(a_1 \times K, a_3 \times K, a_2 \times K,  a_4 \times K) =
1/A.$$
\end{prop}

\medskip
\noindent {\bf Proof.} This follows from
$$1/A= pb_4(a_1, a_3, a_2, a_4)\geq
pb_4(a_1 \times K, a_3 \times K, a_2 \times K, a_4 \times K) \geq
1/A.$$ The equality on the left is guaranteed by
Theorem~\ref{thm-surface-1} and the inequality on the right
follows from \eqref{eq-stab-surf-1}. For the inequality in the
middle which deals with the behavior of the Poisson bracket
invariants under stabilizations we refer to \eqref{eq-pb-products}
below. \qed

\medskip
The previous argument does not work for $M=S^2$. However, in this
case we have a stronger result:

\medskip
\noindent \begin{prop}\label{prop-surfaces-sphere} Assume that $M =
S^2$. Then
$$pb_4(a_1 \times K, a_3 \times K, a_2 \times K,  a_4 \times K) =
\max(1/A, 1/(B-A)).$$
\end{prop}

\medskip
The proof based on a method of symplectic field theory is given in
Section \ref{sec-SFT}  below.

\bigskip

In contrast to the previous situation,  if one stabilizes each of
the sets $a_1$, $a_2$, $a_3$, $a_4$ by {\it its own} exact
Lagrangian section of $T^*S^1$, a transition from rigidity to
flexibility takes place:

\begin{prop}\label{prop-baby-example}
Let $K_1,\ldots,K_4$ be a generic quadruple of sections of $T^*S^1$.
Then
$$pb_4(a_1 \times K_1, a_3 \times K_3, a_2 \times K_2,  a_4 \times
K_4) = 0.$$
%(This is true for any closed symplectic surface $M$ without any restrictions
%on its genus).
\end{prop}

\medskip
The proof will be given in Section~\ref{sec-vanishing}.

\medskip
\begin{rem}
\label{rem-comparison-square-4-meridians} {\rm Let us compare the
situations considered in Proposition~\ref{prop-baby-example} and
Example~\ref{exam-pb4m3}. Namely, just as in
Example~\ref{exam-pb4m3}, consider two parallels $X_0$ and $X_1$
and two meridians $Y_0$ and $Y_1$ on a two-dimensional torus
$\T^2$. They divide $\T^2$ into four squares. Pick three of the
squares, attach a handle to each of them and call the obtained
genus-$4$ surface $M$. Denote the remaining $4$-th square by
$\Pi$, its area by $A$ and its sides by $a_1\subset X_0, a_2
\subset Y_0, a_3\subset X_1, a_4\subset Y_1$ (similarly to
Example~\ref{exam-pb3-surfaces}). Equip $M$ with an area form
$\omega$.  Theorem \ref{thm-surface-1} and monotonicity of $pb_4$
with respect to inclusions of sets (see \eqref{eqn-inclusion1}
below) yield
$$pb_4(X_0,X_1,Y_0,Y_1)\geq pb_4 (a_1,a_3,a_2,a_4)\geq 1/A>0.$$
Furthermore, if $K$ is an exact section of $T^* S^1$ then by
\eqref{eq-stab-surf-1}
$$pb_4(X_0\times K, X_1\times K, Y_0\times K, Y_1\times K)
\geq pb_4 (a_1\times K, a_3\times K, a_2\times K, a_4\times K)\geq
1/A.$$ However, if $K_1,\ldots,K_4$ a generic quadruple of exact
sections of $T^*S^1$, then, according to
Proposition~\ref{prop-baby-example},
$$pb_4 (a_1\times K_1,a_3\times K_3,a_2\times K_2,a_4\times K_4) =
0,$$ while, on the other hand,
$$pb_4({X}_0\times K_1,X_1\times
K_3,{Y}_0\times K_2,{Y}_1\times K_4)
>0.$$ The latter claim follows from  Remark~\ref{rm-stb-more} below
which will be proved by means of ``per\-sis\-tent"
pseudo-holomorphic curves coming from operations in
Donald\-son-Fukaya category. We present this technique right away
in the next section.
 }
\end{rem}
%%%%%%%%%%%%%%%%%%%%%%%%%%%%%%%%%%%%%%%%%%%%%%%%%%%%%%%%%%%%%%%%%%%%%%%%%%%%%

\subsection{Poisson bracket invariants and Lagrangian Floer
homology}\label{subsec-intro-pb-lagrangian}

Recall that a Lagrangian submanifold $L \subset M$ is called {\it
monotone} if there exists a positive {\it monotonicity} constant
$\kappa >0$ so that $\omega(A) = \kappa \cdot m_L(A)$ for every $A
\in \pi_2(M,L)$. Here $m_L: \pi_2(M,L) \to \Z$ is the Maslov class
of $L$.

Suppose for a moment that $\omega(B) \neq 0$ for some $B \in
\pi_2(M)$. Write $B'$ for the image of $B$ in $\pi_2(M,L)$. Then
$$\omega(B) = \omega(B') = \kappa m_L(B') = 2\kappa c_1 (B)\;,$$
where $c_1$ is the first Chern class of $TM$. Thus $\kappa =
\omega(B)/(2c_1(B))$. In particular, in this case the monotonicity
constant $\kappa$ does not depend on the monotone Lagrangian
submanifold $L$.

When $c_1$ and $\omega$ vanish on $\pi_2(M)$ but $m_L$ and $\omega$
do not vanish on $\pi_2(M,L)$ and are positively proportional, the
monotonicity constant $\kappa$ of $L$ may depend on $L$ (think of
circles of different radii in the plane).

Finally, if both $m_L$ and $\omega$ vanish on $\pi_2(M,L)$, the
monotonicity constant of $L$ is not defined uniquely: every $\kappa
>0$ does the job (think of a meridian of the two-torus).

In what follows we deal with monotone Lagrangian submanifolds of a
symplectic manifold $(M,\omega)$. We shall study collections $\cL
=(L_0,L_1,\ldots, L_{k-1})$ of Lagrangian submanifolds in $M$ in general position
satisfying the following topological condition. Consider the set
$\cT_k$ of homotopy classes of $k$-gons in $M$ whose sides (in the
natural cyclic order) lie, respectively, in $L_0,L_1,\ldots,
L_{k-1}$. For every class $\alpha \in \cT_k$ denote by $m(\alpha)$
its Maslov index and by $\omega(\alpha)$ its symplectic area. We say
that $\cL$ is {\it of finite type} if for every $N \in \Z$
\begin{equation}\label{eq-fintype}
A(L_0,\ldots,L_{k-1};N):= \sup \{\omega(\alpha) : \alpha \in
\cT_k, m(\alpha) = N\} < +\infty.
\end{equation}
Here we put $\sup \emptyset := -\infty$.

\medskip
\begin{exam}\label{exam-pb4m3-1} {\rm Four curves on a genus-$4$
surface $M$ as in Example~\ref{exam-pb4m3} form a collection of
finite type (cf. \cite{DeSilva}).   Indeed, recall that $M$ was
obtained by attaching three handles to the torus $\T^2$. Passing
to the abelian cover $\widetilde{M}$ of $M$ associated to the
universal cover of $\T^2$, we see that up to the action of the
group of deck transformations $\Z^2$ and, up to a change of the
orientation, there is a {\it unique} homotopy class of
quadrilaterals in $\widetilde{M}$ with boundaries on the lifts of
our curves. This yields the finite type condition for the curves.
At the same time the quadruple of circles on the {\it torus} (see
Example~\ref{exam-pb4m3}) is not of finite type: passing to the
universal cover $\R^2$ of $\T^2$ we see that there exist index-$2$
squares of arbitrarily large area with boundaries  on the lifts of
our curves (see Remark~\ref{rem-nohandles} for a further
discussion). }
\end{exam}

\medskip
Another class of examples is as follows.

\begin{prop}\label{prop-finite-type}
Assume that all $L_i$'s have the same monotonicity constant and
the morphism $\pi_1(L_i) \to \pi_1(M)$ has a finite image for
every $i$. Then the collection $\cL$ is of finite type.
\end{prop}

\medskip The proof is given in Section~\ref{subsec-finitetype}. In
what follows we deal only with collections of finite type.

The main results of this section involve Floer theory of monotone
Lagrangian submanifolds. We write $HF(L_0,L_1)$ for the Lagrangian
Floer homology and denote by $\mu^k$ the operations in (Donaldson-)Fukaya
category. We refer to Section~\ref{sec-pb-pseudo-holo} for
preliminaries.

\begin{thm}\label{thm-triang} Let $L_0,L_1,L_2$, $L_0 \cap
L_1 \cap L_2 = \emptyset$, be a finite type collection of closed
Lagrangian submanifolds of $M^{2n}$. Assume that the product
$$\mu^2: HF(L_0,L_1) \otimes HF(L_1,L_2) \to HF(L_0,L_2)$$
is well-defined, invariant under exact Lagrangian isotopies and
does not vanish. Then $$pb_3 (L_0, L_1, L_2)
>\frac{1}{2A(L_0,L_1,L_2;2n)}.$$
\end{thm}

The proof will be given further in this section.

\medskip
\begin{exam} {\rm Let $L$ be a closed connected manifold with a
finite fundamental group and let $M:=T^*L$ be equipped with the
standard symplectic structure. Identify $L$ with the zero section
of $T^* L$. The group $HF(L,L)$ and the product $\mu^2$ on it are
non-trivial: $HF(L,L)$ is isomorphic to the singular homology of
$L$ \cite{Floer-unreg}. Under this isomorphism the product in the
Floer homology corresponds to the classical intersection product
in the singular homology of $L$ \cite{FO}. Thus $pb_3(L_0,L_1,L_2)
>0$ for three exact sections of $T^*L$ in general position.}
\end{exam}

\medskip
\noindent \begin{exam} {\rm Let $M:=S^2\times S^2$ be equipped with
the symplectic structure $\omega\oplus\omega$, where $\omega$ is an
area form on $S^2$. Let $L:= \{ (x,-x)\in S^2\times S^2\}$ be the
anti-diagonal. It is a Lagrangian sphere. The group $HF(L,L)$ and
the product on it are non-trivial -- one should just recall that
$(S^2\times S^2, \omega\oplus\omega)$ is symplectomorphic to a
quadric in $\C P^3$ with the symplectic structure induced by the
Fubini-Study form and apply \cite{Biran-Cornea-GT}, Theorem 2.3.4
and Remark 2.2.1.  Thus $pb_3(L_0,L_1,L_2) >0$ for three generic
images of $L$ under Hamiltonian isotopies. }\end{exam}

\medskip
A more sophisticated example where $L_i$'s are Lagrangian spheres
and the triangle product is non-trivial is given in
Section~\ref{subsec-triangle-example} below.

Assume now that we have a finite type collection $L_0, L_1, L_2,
L_3\subset M^{2n}$ of Lagrangian submanifolds such that
\begin{equation}
\label{eq-lag-intersect} L_0\cap L_2=L_1\cap L_3=\emptyset.
\end{equation}
Assuming that the Lagrangian Floer homology groups $HF(L_i,L_j)$
are well-defined, one can define the $\mu^3$-operation in the
Donaldson-Fukaya category:
$$\mu^3: HF(L_0,L_1)\otimes HF(L_1,L_2)\otimes HF(L_2,L_3)\to
HF(L_0,L_3),$$ provided it is well-defined on the chain level.
%defined on the chain level, descends to an
%operation on Floer homology. We shall still denote it by $\mu^3$.
For such a collection of Lagrangian submanifolds we have the following result.

\begin{thm} \label{thm-quadrilat} Assume that the
operation $\mu^3$ is well-defined, invariant under exact
Lagrangian isotopies preserving the intersection condition
\eqref{eq-lag-intersect} and does not vanish. Then
\begin{equation} \label{eq-pb4-Lagra}
pb_4 (L_0, L_2, L_1, L_3) \geq 1/A(L_0, L_1, L_2, L_3;3n-1).
\end{equation}
\end{thm}

\medskip
For the proof see Section~\ref{sec-mainlagthm}. This theorem is
applicable, for instance, to the quadruple of curves on the
genus-$4$ surface from Example~\ref{exam-pb4m3} and their
stabilizations. More sophisticated examples in which $\mu^3$ does
not vanish were found by Smith in \cite{Smith}.

In Section~\ref{sec-mainlagthm} below we discuss an extension of
the lower bounds on the Poisson bracket invariants provided by
Theorems~\ref{thm-triang} and \ref{thm-quadrilat} to
stabilizations of collections of Lagrangian submanifolds. In view
of Theorem~\ref{thm-non-autonomous}, such non-trivial lower bounds
on $pb_4$ for the stabilized Lagrangian submanifolds yield the
existence of Hamiltonian chords for non-autonomous Hamiltonian
flows.

Let us prove Theorem~\ref{thm-triang}  skipping some
technicalities and preliminaries on the operations in Lagrangian
Floer homology which will be given in Section~\ref{sec-mainlagthm}
below.

\medskip
\noindent {\bf Proof of Theorem~\ref{thm-triang}.} We follow the
strategy described in Section~\ref{subsec-deformation} above: Take
a pair of functions $(F,G) \in \cF'_3(L_0,L_1,L_2)$ and consider
the deformation of the symplectic form $\omega$ given by
$$\omega_s := \omega -sdF\wedge dG,$$
where $s \in I:= [0;1/||\{F,G\}||)$. Observe that $\omega_s$ is
cohomologous to $\omega$ and, moreover, $\omega_s$ and $\omega$
represent the same relative cohomology classes in $H^2(M,L_i)$,
$i=0,1,2$. Thus, by Moser's theorem \cite{Moser}, there exists an
ambient isotopy
 $f_s : M \to M$ with $f_s^*\omega_s =
\omega$. Furthermore\footnote{Warning: in general there is no
ambient Hamiltonian isotopy of $M$ taking $L_i$ to $L_i^s$ for all
$i$ simultaneously!}, $L_i^s:= f_s^{-1}(L_i)$ is an exact isotopy of
$L^0_i = L_i$. Thus the product in the Lagrangian Floer homology does
not change with $s$.

Choose a generic family of almost complex structures $J_s$, $s \in
I$, compatible with $\omega_s$. The non-vanishing of the product
in the Lagrangian Floer homology guarantees that for every $s \in
I$ there exists a $J_s$-holomorphic triangle, say $\Sigma$, whose
$i$-th side lies on $L_i$ for $i=0,1,2$. The dimension of the
moduli space of such triangles equals $m(\Sigma)-2n = 0$ (see
\eqref{exp-dim} below) and thus the finite type condition
\eqref{eq-fintype} guarantees that $\omega(\Sigma) \leq A$ with $A
= A(L_0, L_1, L_2;2n)$. Observe that, by the Stokes formula,
$\int_{\partial \Sigma} FdG = 1/2$. Hence
Proposition~\ref{prop-psh-curves-general}
 yields
 $$||\{F,G\}|| \geq \frac{1}{2A}$$ and therefore
 $$pb_3(L_0,L_1,L_3) \geq \frac{1}{2A}.$$ \Qed

\bigskip
\noindent {\bf Organization of the paper.} In
Section~\ref{sec-pb-3-pb4-defs} we discuss  basic properties of the
Poisson bracket invariants.

In Section~\ref{sec-sympl-approx} we prove the results on
symplectic approximation stated in
Section~\ref{subsec-intro-approx} above and discuss a
generalization of Theorem~\ref{thm-dichotomy} (ii) to the case of
iterated Poisson brackets.

In Section~\ref{sec-chords-proofs} we establish the existence of
Hamiltonian chords (see Section~\ref{subsec-intro-chords})  and
discuss more examples and applications.

In Section~\ref{sec-pb-pseudo-holo} we give preliminaries on
Lagrangian Floer homology and operations in Donaldson-Fukaya
category. We use them for the proof of the results   stated in
Section~\ref{subsec-intro-pb-lagrangian} above.

In Section~\ref{sec-SFT} we apply symplectic field theory to
calculation of the Poisson bracket invariant for a stabilized
quadrilateral on the two-sphere.

In Section~\ref{sec-vanishing} we present a sufficient condition for
the vanishing of the Poisson bracket invariants.

In Section~\ref{sec-discussion} we formulate various open problems
and outline directions of further study. We present connections to
control theory, speculate on an extension of the Poisson bracket
invariants to $k$-tuples of sets for $k>4$ and continue the
discussion on vanishing of $pb_3$ and $pb_4$.

\section{Preliminaries on Poisson bracket invariants}
\label{sec-pb-3-pb4-defs}

\subsection{Definitions and notations}\label{subsec-defnot2}

Let $(M^{2n},\omega)$ be a connected symplectic manifold (either open or closed).
We use the following sign conventions in the definitions of a
Hamiltonian vector field and the Poisson bracket on $M$: the
Hamiltonian vector field $\sgrad F$ of a Hamiltonian $F$ is
defined by
$$i_{\sgrad {F}}
\omega = -dF$$ and the Poisson bracket of two Hamiltonians $F$,
$G$ is given by
$$\{F,G\} := \omega(\sgrad G,\sgrad F) = dF (\sgrad G) = - dG
(\sgrad F) =$$
$$= L_{\sgrad G} F = - L_{\sgrad F} G.$$

Let $\cX=(X_1,\ldots,X_k)$ be an ordered collection of $k$ compact
subsets of a symplectic manifold $(M,\omega)$. In what follows
$pb(\cX)$ stands for $pb_3(X_1,X_2,X_3)$ if $k=3$ and for
$pb_4(X_1,X_2,X_3,X_4)$ if $k=4$. Furthermore, we write
$pb(\cX;M)$ whenever we wish to emphasize dependence of the
Poisson bracket invariants on the ambient symplectic manifold $M$.

Let $\cX$ and $\cY$ be two collections as above (with the same $k$).
We say that $\cX \subset \cY$ if $X_i \subset Y_i$ for all $i$.
Given a compact subset $Y$ of a manifold $N$, we put
$$\cX \times Y := (X_1 \times Y,\ldots,X_k \times Y).$$

Let us say that a sequence of subsets $X^{(j)}$ of $M$ {\it
converges} to a limit set $Y$ if every open neighborhood of $Y$
contains all but a finite number of sets from the sequence. This
is denoted by $X^{(j)} \to Y$. Given  collections $\cX^{(j)}$ and
$\cY$, we write $\cX^{(j)} \to \cY$ if $X^{(j)}_i \to Y_i$ for all
$i=1,\ldots,k$.

\subsection{Basic properties of Poisson bracket invariants}

All the properties listed below (except the last one) readily follow
from the definitions and Proposition~\ref{prop-op}.

\medskip
\noindent{\sc Semi-continuity:}
\begin{prop}\label{prop-Hausdorff}
\label{cor-Hausdorff} Suppose that a sequence of collections
$\cX^{(j)}$ of $k=3$ or $4$ ordered subsets of a symplectic manifold
converges to a collection $ \cY$. Then
$$\limsup_{j\to +\infty} pb(\cX^{(j)}) \leq pb(\cY).$$
\end{prop}

\medskip
\noindent{\sc Behavior under symplectic embeddings:} Assume that
$(M,\omega)$ and $(N,\sigma)$ are symplectic manifolds of the same
dimension. Let $\phi: M \to N$ be a symplectic embedding. Let $\cX$
and $\cY$ be collections of $k$ ordered subsets of $M$ and $N$
respectively with $\phi(\cX) \supset \cY$. Then
\begin{equation}
\label{eqn-inclusion} pb(\cX;M) \geq pb(\cY;N).
\end{equation}
In particular, if $\cX,\cY$ are collections of $k$ ordered subsets
of $M$, then
\begin{equation}
\label{eqn-inclusion1} \cX\supset \cY\ \Longrightarrow\ pb(\cX) \geq
pb(\cY).
\end{equation}

\medskip
\noindent{\sc Behavior under products:}
 Suppose that $M$ and $N$ are connected
symplectic manifolds. Equip $M\times N$ with the product symplectic
form. Let $A\subset N$ be a compact subset. Then for every
collection $\cX$ of $k=3$ or $4$ compact subsets of $M$
\begin{equation}\label{eq-pb-products}
pb(\cX,M) \geq pb(\cX \times A, M \times N).
\end{equation}

\medskip
\noindent{\sc Comparing $pb_3$ and $pb_4$:} The invariants $pb_3$
and $pb_4$ are related by the following inequality.

\begin{prop}
\label{prop-pb-3-4-relation} Let $X_0,X_1,X_2,X_3 \subset M$ be
compact subsets such that
$$X_0 \cap X_2 = X_1 \cap X_3 =\emptyset.$$
Then
\begin{equation}
\label{eqn-pb3-pb4-max} pb_4 (X_0,X_2,X_1,X_3)\geq \max_{i=0,1,2,3}
pb_3 (X_i\cup X_{i+1}, X_{i+1}\cup X_{i+2}, X_{i+3}), \end{equation}
where all the indices are taken modulo $4$.
\end{prop}

\medskip
Combining inequality \eqref{eqn-pb3-pb4-max} with monotonicity
property \eqref{eqn-inclusion1} we get that
\begin{equation}\label{eq-pb3-pb4}
pb_4(X_0,X_2,X_1,X_3)\geq pb_3(X_0,X_1,X_2 \cup X_3).
\end{equation}

\medskip
\noindent{\sc Expansion property:}

\begin{prop}\label{prop-expan}
Consider a quadruple of compact subsets $X_0,X_1,Y_0,Y_1 \subset M$
such that $X_0 \cap X_1 = Y_0 \cap Y_1 =\emptyset.$ Let $A \subset
Y_0$ be a compact subset disjoint from $X_1$. Then
\begin{equation}\label{eq-expan}
pb_4(X_0 \cup A,X_1,Y_0,Y_1) = pb_4(X_0,X_1,Y_0,Y_1).
\end{equation}
\end{prop}

\medskip
The proof is given below in
Section~\ref{sec-pfs-basic-props-pb3-pb4}.

\subsection{Proofs of the basic properties of $pb_3$ and $pb_4$}
\label{sec-pfs-basic-props-pb3-pb4}

\noindent {\bf Proof of Proposition~\ref{prop-op}(i).}

\noindent \begin{lemma}\label{lem-cutoff}  Denote by $\Delta \subset
\R^2(s,t)$ the triangle $s\geq 0,t\geq 0, s+t \leq 1$. Then for
every $\kappa > 0$ there exists a smooth map $T= (T_1,T_2):
\R^2(s,t) \to \Delta$ and $\delta=\delta (\kappa)>0$ so that

\begin{itemize}

\item{} $\lim_{\kappa\to 0} \delta(\kappa) = 0$,

\item{} $T$ maps $\{s \leq \delta\}$ to $\{s=0\}$, $\{t \leq
\delta\}$ to $\{t=0\}$, $\{s+t \geq 1-\delta \}$ to $\{s+t=1\} $,

\item{} $||\{T_1,T_2\}|| \leq 1+\kappa$. (The Poisson bracket $\{
T_1, T_2\}$ is taken with respect to the standard area form on
$\R^2$).

\end{itemize}

\end{lemma}

\begin{figure}
 \begin{center}
\scalebox{0.5}{\includegraphics*{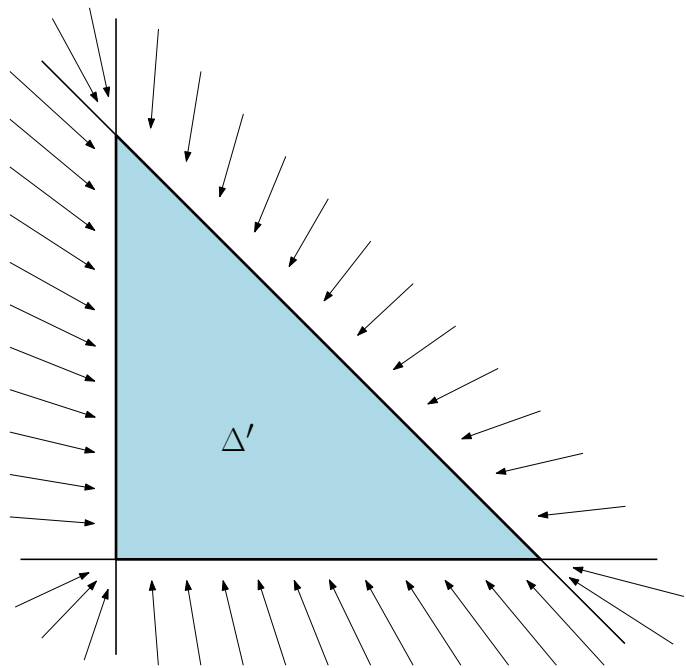}}
 \end{center}
  \caption{A piece-wise projective approximation to $T$}
   \label{fig4}
\end{figure}

\medskip\noindent{\bf Proof.} Take $K > 1$. Set $\delta = 1 - \frac{1}{K}$. Let
$\Delta'$ be the triangle bounded by the lines $ l_1' = \{s =
\delta\} $, $ l_2' = \{t = \delta\} $ and $ l_3' = \{ s+t =
1-\delta\} $. The desired map $T$ will be obtained as a perturbation of
a piece-wise projective map $\R^2 \to \Delta'$ presented on Figure \ref{fig4}.
Denote by $ a_1',a_2',a_3' $ the sides of $ \Delta'
$ lying on $ l_1',l_2',l_3' $ respectively, and by $
v_1',v_2',v_3' $ the opposite vertices. The lines $ l_1',l_2',l_3'
$ divide the plane into 7 closed domains: $\Delta'$, 3 exterior
angles $A(v_i')$ corresponding to the vertices $v_i'$ and 3
unbounded domains $D(a_i')$, so that $D(a_i')$ has the side $a_i'$
as a part of its boundary.

Pick a vertex $ v_i' $. Introduce polar coordinates $ (r,\theta) $
on the plane so that $ v_i' $ is the center of the coordinate
system. When a point $ x $ runs through the straight line $ l_i'
$, the value of $ \theta(x) $ runs through an open interval in
$S^1$ -- denote it by $ J_i \subset S^1 $. For any $ \theta \in
J_i $ denote by $ R_i(\theta) $ the distance between $ v_i' $ and
the intersection point of $ l_i' $ with the ray from $ v_i' $
having the angle $ \theta $ (note that $R_i (\theta)\to +\infty$
as $\theta$ approaches an end-point of $J_i$). Let $ \psi : (0;
+\infty) \rightarrow \R $ be a smooth function so that $ 0 \leq
\psi'(t) < K $ for $ t \in (0;+\infty) $, $ \psi(t) = t $ for $ t
\in (0;1/2) $, $ \psi(t) = 1 $ for $ t \in [1;+\infty) $. In
particular, $ \psi(t) < K t $. Define $ \phi_i: (0;+\infty) \times
S^1 \rightarrow (0;+\infty) $ by
\begin{equation*}
 \phi_i(r,\theta) =
 \left\{
\begin{array}{lllll}
R_i(\theta) \psi(\frac{r}{R_i(\theta)}) , \text{ if } \theta \in J_i, \\[0.05in]
r, \text{ if } \theta \notin J_i. \\[0.05in]
\end{array} \right.
\end{equation*}
An easy check shows that $ \phi_i $ is a smooth function. We have
$ \frac{\partial \phi_i}{\partial r}(r,\theta) < K $ and $
\phi_i(r,\theta) < K r $ for $ (r,\theta) \in (0;+\infty) \times
S^1 $. Define $ \Phi_i : \R^2 \rightarrow \R^2 $ by $
\Phi_i(r,\theta) = (\phi_i(r,\theta),\theta) $ for $ (r,\theta)
\in (0;+\infty) \times S^1 $, and $ \Phi_i(v_i') = v_i' $. Then at
the point $ (r,\theta) $ we have
$$ \Phi_i^* \omega = \frac{\phi_i(r,\theta)}{r} \frac{\partial \phi_i}{\partial r}(r,\theta) \omega $$
and hence $$ | d \Phi_i | = \frac{\phi_i(r,\theta)}{r} \frac{\partial \phi_i}{\partial r}(r,\theta) < K^2 ,$$
since we know that $ \frac{\partial \phi_i}{\partial r}(r,\theta) < K $ and $
\phi_i(r,\theta) < K r $ for $ (r,\theta) \in (0;+\infty) \times
S^1 $. The map $ \Phi_i $ maps
the region $ D(a_i') $ onto $ a_i' $. Moreover, it maps the region
$ A(v_{i+1}') $ (we use the cyclic numbering of vertices modulo
$3$) onto the ray starting at $v_{i+1}'$  and going outwards from
$\Delta'$ along $l_i'$. Similarly, $ \Phi_i $ maps the region $
A(v_{i+2}') $ onto the ray starting at $v_{i+2}'$ and going
outwards from $\Delta'$ along $l_i'$.

Consider an affine map $ \Psi : \R^2 \rightarrow \R^2 $ defined by
$$ \Psi(s,t) =
\bigg(\frac{s-\delta}{1-\delta}, \frac{t-\delta}{1-\delta}\bigg)
$$ in the standard coordinates $ (s,t) $. Then $ \|d\Psi\|= K^2 $.
Now define $ T: \R^2 \rightarrow \R^2 $ by $ T = \Psi \circ \Phi_1
\circ \Phi_2 \circ \Phi_3 $. We have $ | d T | \leqslant K^8 $, and
it is easy to see that $T$ maps $\{s \leq \delta\}$ to $\{s=0\}$,
$\{t \leq \delta\}$ to $\{t=0\}$, $\{s+t \geq 1-\delta \}$ to
$\{s+t=1\} $. Finally it remains to take $ K = \sqrt[8]{1+\kappa} $.
\Qed

\bigskip Put
$$pb'_3 := \inf_{(F',G')\in \cF'_3 (X,Y,Z)} ||\{F',G'\}||. $$
Clearly, $pb_3 \leq pb'_3$. Thus it is enough to show that $pb_3\geq
pb'_3$.

Indeed, let $(F,G)\in \cF_3 (X,Y,Z)$, that is
$$F|_X \leq 0, G|_Y \leq 0, (F+G)|_Z \geq 1.$$
Take $T$ from Lemma~\ref{lem-cutoff} (with a small enough
$\kappa>0$) and put
$$F' = T_1(F,G), G' = T_2(F,G).$$
An immediate check shows that $(F',G')\in \cF'_3 (X,Y,Z)$ and
$$||\{F',G'\}|| \leq ||\{F,G\}||(1+\kappa).$$
Choosing $\kappa$ arbitrarily small and taking the infimums over
$(F',G')$ and $(F,G)$ in both sides of the inequality we get that
$$pb'_3\leq pb_3,$$
and hence $pb_3=pb'_3$, as required. \Qed

\medskip
\noindent {\bf Proof of Proposition~\ref{prop-op}(ii).}
Set
$$pb_4 := pb_4(X_0,X_1,Y_0,Y_1),$$
$$pb'_4:= \inf_{(F,G)\in \cF'_4 (X_0,X_1,Y_0,Y_1)} ||\{F,G\}||.$$
Clearly,
$$pb_4\leq  pb'_4,$$
so it suffices to prove that
\begin{equation} \label{eq-equality-vsp}
pb_4 \geq pb'_4.
\end{equation}

Fix $\epsilon > 0$.  Choose $(F,G)\in \cF_4 (X_0, X_1, Y_0,Y_1)$ so
that
$$pb_4 \geq ||\{F,G\}||-\epsilon.$$
Fix a small enough $\delta>0$ and choose a smooth non-decreasing
function $u: \R \to [0;1]$ so that $u(s) = 0$ for $s \leq \delta$,
$u(s) =1$ for $s \geq 1-\delta$ and $u'(s) \leq 1+2\delta$. Put
$\phi = u \circ F$ and $\psi = u \circ G$. An immediate check shows
that $(\phi,\psi)\in \cF'_4 (X_0, X_1, Y_0,Y_1)$. Now note that
$$
pb'_4 \leq ||\{\phi,\psi\}|| \leq
$$
$$\leq
(1+2\delta)^2 \cdot ||\{F,G\}|| \leq (1+2\delta)^2 (pb_4
+\epsilon).$$ Choosing $\epsilon$ and $ \delta$ arbitrarily small,
we get \eqref{eq-equality-vsp} which completes the proof.\Qed

\medskip
\noindent {\bf Proof of Proposition~\ref{prop-expan}.} By
monotonicity,
$$pb_4(X_0 \cup A,X_1,Y_0,Y_1) \geq pb_4(X_0,X_1,Y_0,Y_1).$$
Let us prove the reverse inequality. Fix $\epsilon >0$. By
Proposition~\ref{prop-op}(ii), there exist functions $F,G\in
C^\infty (M)$ with $F = 0$ on $\op(X_0)$, $F=1$ on $\op(X_1)$, $G=0$
on $\op(Y_0)$, $G=1$ on $\op(Y_1)$, and
$$pb_4(X_0,X_1,Y_0,Y_1) \geq ||\{F,G\}||-\epsilon.$$ Let $U
\subset \op(Y_0) $ be a neighborhood of $A$, $U \cap X_1
=\emptyset$. Choose a smooth cut-off function $u: M \to [0;1]$
which vanishes on $A$ and equals $1$ outside $U$. Put $F' = uF$.
Note that
$$\{F',G\} = u\{F,G\} + F\{u,G\}.$$
If $x \in U$, the function $G$ is constant near $x$ and hence
$\{F',G\} = \{F,G\} =0$. If $x \notin U$, we have $u=1$ near $x$
and hence again $\{F',G\} = \{F,G\}$. Thus $\{F',G\} = \{F,G\}$
everywhere and therefore, since $F' = 0$ on $X_0 \cup A$ and
$F'=1$ on $X_1$, we get that
$$pb_4(X_0 \cup A,X_1,Y_0,Y_1) \leq ||\{F',G\}||=||\{F,G\}|| \leq pb_4(X_0,X_1,Y_0,Y_1)+\epsilon.$$
Since this holds for every $\epsilon >0$, we get that
$$pb_4(X_0 \cup A,X_1,Y_0,Y_1) \leq  pb_4 (X_0,X_1,Y_0,Y_1)$$
which completes the proof.\Qed

\section{Symplectic approximation}
\label{sec-sympl-approx}

In this section we prove the results on symplectic approximation
stated in Section~\ref{subsec-intro-approx}.

\subsection{The Lipschitz property of the function $s\rho_{F,G}$}\label{subsec-sympl-approx-pfs-lip}

\noindent {\bf Proof of Proposition
\ref{prop-profile-fcn-Lipschitz}:} Fix $ c > 0 $. Then for any $
\epsilon > 0 $ there exist $ F',G' $ such that $ \| \{ F', G' \}
\| \leqslant c $ and $ \|F-F'\| + \|G-G'\| < \rho_{F,G}(c) +
\epsilon $. Take any $ 0 < \lambda < 1 $ and denote $ F_1 =
\lambda F' $, $ G_1 = G' $. Then $ (F_1,G_1) \in K_{\lambda c} $,
and we have $$ \|F-F_1 \| + \|G-G_1 \| = \|F- \lambda F'\| +
\|G-G'\| \leqslant $$ $$ \leqslant \|F- F'\| + (1-\lambda) \|F' \|
+ \|G-G'\| \leqslant \rho_{F,G}(c) + \epsilon + (1-\lambda) \|F'
\| \leqslant $$ $$ \leqslant \rho_{F,G}(c) + \epsilon +
(1-\lambda) \|F' - F \| + (1-\lambda) \| F \| \leqslant
$$ $$ \leqslant \rho_{F,G}(c) + \epsilon + (1-\lambda)
\rho_{F,G}(c) + (1-\lambda)\epsilon + (1-\lambda) \| F \| \leqslant
$$ $$ \leqslant \rho_{F,G}(c) + \epsilon + (1-\lambda) \| F \| +
(1-\lambda)\epsilon + (1-\lambda) \| F \| = $$ $$ = \rho_{F,G}(c) +
2 \| F \| (1 - \lambda) + (2-\lambda)\epsilon .$$ Therefore $
\rho_{F,G}(\lambda c) \leqslant \rho_{F,G}(c) + 2 \| F \| (1 -
\lambda) + (2-\lambda)\epsilon $. This is true for any $ \epsilon >
0 $, hence $ \rho_{F,G}(\lambda c) \leqslant \rho_{F,G}(c) + 2 \| F
\| (1 - \lambda) $. Therefore $$ | \lambda c \rho_{F,G}(\lambda c) -
c \rho_{F,G}(c) | \leqslant | \lambda c \rho_{F,G}(\lambda c) - c
\rho_{F,G}(\lambda c) | + | c \rho_{F,G}(\lambda c) - c
\rho_{F,G}(c) | = $$ $$ = \rho_{F,G}(\lambda c) | \lambda c - c | +
c ( \rho_{F,G}(\lambda c) - \rho_{F,G}(c) ) \leqslant $$ $$
\leqslant \|F\| ( c - \lambda c) + 2 c \| F \| (1 - \lambda) = $$ $$
= 3 \|F\| ( c - \lambda c ) .$$ This is true for any $ c > 0 $, $ 1
> \lambda > 0 $. Therefore the function $ s \mapsto s \rho_{F,G}(s)
$ is Lipschitz on $ (0; +\infty) $ with the Lipschitz constant $ 3
\| F \| $. Since $\rho_{F,G}=\rho_{G,F}$, the same property is
true with the Lipschitz constant $3 ||G||$ which finishes the
proof. \Qed

%%%%%%%%%%%%%%%%%%%%%%%%%%%%%%%%%%%%%%%%%%%%%%%%%%

\subsection{The profile function and Poisson bracket invariants}\label{subsec-sympl-approx-pfs-pb}

\noindent {\bf Proof of Theorem~\ref{thm-dichotomy}(ii).} Assume,
without loss of generality, that $0 \leq F \leq 1$. Choose $H \in
C^{\infty}_c (M)$ such that $H = 1/2$ on the union of the supports of
$F$ and $G$ and $0 \leq H \leq 1/2$. Then $\{H,G\} = 0$ and
$$\rho_{F,G}(0) \leq d((F,G), (H,G)) = ||F-H|| \leq 1/2.$$
Therefore, as soon as we prove \eqref{eq-rho-new} and
\eqref{eq-rho-new-pb4},  we would get $\rho_{F,G}(0)=1/2$.
Furthermore,
for any $t\in [0,1]$
we have $||\{ tF+(1-t)H, G\}|| = ||\{tF,G\}||$ and therefore
the pair $(tF+(1-t)H,G)$ lies in $\cK_{t||\{ F,G\}||}$.
Thus $\rho_{F,G} (t||\{ F,G\}||)\leq d((F,G), (tF+(1-t)H,G))=(1-t)||F-H||\leq 1/2 - t/2$.
Setting $s:=t||\{ F,G\}||$ we get  $\rho_{F,G} (s)\leq 1/2-s/(2||\{ F,G\}||)$, that is \eqref{eq-rho-up-trivial}.

Thus it remains to prove inequalities \eqref{eq-rho-new} and
\eqref{eq-rho-new-pb4}.

\medskip \noindent {\sc The case} $(F,G) \in \cF_3^\flat(X,Y,Z)$:

\noindent Fix $s\in [0;p)$. Suppose that $\rho_{F,G}(s) < 1/2$
(otherwise \eqref{eq-rho-new} follows automatically). Take any
$\delta\in (\rho_{F,G} (s), 1/2)$.  Take $(H,K)\in \cF=C^\infty_c
(M)\times C^\infty_c (M)$ with  $d((F,G),(H,K)) \leq \delta$ and
$$||\{H,K\}|| \leq s.$$
Put
$$\alpha := ||F - H||, \beta := ||G -K||.$$
Thus $\alpha +\beta\leq \delta$. Furthermore,
$$H|_X \leq \alpha, K|_Y \leq \beta, (H +K)|_Z
\geq 1 -\delta.$$ Set $$H_1: = \frac{H-\alpha}{1-2\delta},\
K_1:=\frac{K-\beta}{1-2\delta}.$$ Then
$$H_1|_X \leq 0,\, K_1|_Y \leq 0,\, (H_1+K_1)|_Z \geq 1.$$
Note that, unlike $H$ and $K$, the functions $H_1, K_1$ are not
necessarily compactly supported but are constant outside a compact
set (if $M$ is not closed). Pick a smooth compactly supported
function $u: M\to [0;1]$ so that $u\equiv 1$ on an open
neighborhood of $X\cup Y\cup Z\cup {\rm supp}\, H\cup {\rm supp}\,
K$. Set $H_2:=uH_1$, $K_2:=uK_1$. An easy check shows that
$$||\{ H_2, K_2\}||=||\{ H_1, K_1\}||.$$
On the other hand, $(H_2,K_2)\in\cF$ and
$$H_2|_X \leq 0,\, K_2|_Y \leq 0,\, (H_2+K_2)|_Z \geq 1.$$
Therefore, by the definition of $pb_3$, we have
$$||\{H_2,K_2\}|| \geq p=pb_3 (X,Y,Z).$$ Note that
$$||\{ H_2, K_2\}||=||\{ H_1, K_1\}|| = \frac{1}{(1 -2\delta)^2}\cdot ||\{ H,K\}||
\leq \frac{s}{(1 -2\delta)^2}.$$ Thus
$$\frac{s}{(1 -2\delta)^2}\geq p$$
and therefore
$$\delta\geq \frac{1}{2}-\frac{1}{2\sqrt{p}}\cdot\sqrt{s}.$$
Since this is true for every $\delta \in (\rho_{F,G}(s),1/2)$ we get
that
$$\rho_{F,G}(s) \geq \frac{1}{2}-\frac{1}{2\sqrt{p}}\cdot\sqrt{s},$$
as required.

\medskip \noindent {\sc The case} $(F,G) \in \cF_4^\flat(X_0,X_1,Y_0,Y_1)$:

\noindent Fix $s\in [0;p)$. Suppose that $\rho_{F,G}(s) < 1/2$
(otherwise \eqref{eq-rho-new-pb4} follows automatically). Take any
$\delta\in (\rho_{F,G} (s), 1/2)$.  Take $(H,K)\in \cF=C^\infty_c
(M)\times C^\infty_c (M)$ with  $d((F,G),(H,K)) \leq \delta$ and
$$||\{H,K\}|| \leq s.$$
Put
$$\alpha := ||F - H||, \beta := ||G -K||.$$
Thus $\alpha +\beta\leq \delta$ and, in particular, $0\leq
\alpha,\beta \leq \delta < 1/2$. Furthermore,
$$H|_{X_0} \leq \alpha, H|_{X_1} \geq 1-\alpha, K|_{Y_0} \leq \beta,
K|_{Y_1} \geq 1 -\beta.$$ Set $$H_1: =
\frac{H-\alpha}{1-2\alpha},\ K_1:=\frac{K-\beta}{1-2\beta}.$$ Then
$$H_1|_{X_0} \leq 0, H_1|_{X_1} \geq 1, K_1|_{Y_0} \leq 0,
K_1|_{Y_1} \geq 1.$$ Note that, unlike $H$ and $K$, the functions
$H_1, K_1$ are not necessarily compactly supported but are
constant outside a compact set (if $M$ is not closed). Pick a
smooth compactly supported function $u: M\to [0;1]$ so that
$u\equiv 1$ on an open neighborhood of $X_0\cup X_1\cup Y_0\cup
Y_1\cup {\rm supp}\, H\cup {\rm supp}\, K$. Set $H_2:=uH_1$,
$K_2:=uK_1$. An easy check shows that $$||\{ H_2, K_2\}||=||\{
H_1, K_1\}||.$$ On the other hand, $(H_2,K_2)\in\cF$ and
$$H_2|_{X_0} \leq 0, H_2|_{X_1} \geq 1, K_2|_{Y_0} \leq 0,
K_2|_{Y_1} \geq 1.$$ Therefore, by the definition of $pb_4$, we
have
$$||\{H_1,K_1\}|| = ||\{H_2,K_2\}|| \geq p=pb_4(X_0, X_1, Y_0, Y_1).$$
Note that
$$||\{ H_1, K_1\}|| = \frac{1}{(1 -2\alpha)(1-2\beta)}\cdot ||\{ H,K\}||.$$
Since $\alpha,\beta \geq 0$ and $\alpha+\beta \leq \delta$,
$$(1 -2\alpha)(1-2\beta) \geq 1-2(\alpha+\beta)+4\alpha\beta \geq  1-2\delta.$$
Thus $s/(1-2\delta) \geq p$ and therefore $\delta\geq 1/2-s/(2p)$.
Since this is true for every $\delta \in (\rho_{F,G}(s),1/2)$, we
get that $\rho_{F,G}(s) \geq 1/2-s/(2p)$ as required. \Qed

%%%%%%%%%%%%%%%%REMARK ADDED%%%%%%%%%%%%%%%%%%%%%%%%%
\begin{rem}

{\rm Theorem~\ref{thm-dichotomy}(ii) has the following generalization concerning
{\it iterated} Poisson brackets of two functions.

Namely, denote by $\cP_N$, $N\geq 2$, the set of Lie monomials in
two variables of degree $N$ (i.e. if the Lie brackets are denoted by
$\{\cdot,\cdot\}$, the set $\cP_2$ consists
 of $\{A,B\}$, $\cP_3$ of $\{\{A,B\},A\}$ and $\{\{A,B\},B\}$, and so on).
 For $F,G \in (C^\infty_c (M), \{\cdot,\cdot\})$ set
$$Q_N(F,G) := \sum_{p \in \cP_N} || p(F,G)||,$$
$$\cK_s^{(N)}:= \{(F,G) \in \cF : Q_N(F,G) \leq
s\}.$$ In particular, for $N=2$ we get the sets $\cK_s$ defined in
Section 1.2: $\cK_s=\cK_s^{(2)}$. The sets $\cK_s^{(N)}$ can be
viewed as ``tubular neighborhoods" of the set of Poisson-commuting
pairs of functions on $M$: indeed, a symplectic version of the
Landau-Hadamard-Kol\-mo\-go\-rov inequality (see \cite{EP-Poisson2},
\cite{EPR}) implies that $\cK_0^{(N)}=\cK_0$ for any $N$. Now,
similarly to $\rho_{F,G}$, define a new profile function (cf.
\cite{EPR}):
$$\rho_{F,G}^{(N)} (s):= d((F,G),\cK_s^{(N)}).$$
In particular, for $N=2$ we get the profile function $\rho_{F,G}$
studied above: $\rho_{F,G}^{(2)}=\rho_{F,G}$.

It turns out that, similarly to Theorem~\ref{thm-dichotomy}(ii),
for certain $(F,G)$ one can estimate $\rho_{F,G}^{(N)} (s)$ from
below for small $s$ using an analogue of $pb_3$ for iterated
Poisson brackets. Namely, given a triple $(X,Y,Z)$ of compact
subsets of $M$ with $X \cap Y \cap Z = \emptyset$, define
$$pb_3^{(N)} (X,Y,Z):= \inf_{(F,G)} Q_N (F,G),$$
where the infimum is taken over $\cF_3 (X,Y,Z)$.

Then the proof of Theorem~\ref{thm-dichotomy}(ii)  can be
carried over directly to the case of iterated Poisson brackets
yielding the following claim:

\medskip
{\it Put $p_N= pb_3^{(N)} (X,Y,Z)$  and
let $\cF_3^{\flat}$ be defined as in Section 1.3. Assume that $p_N
>0$.

Then for every $(F,G) \in \cF_3^{\flat}$ the profile function
$\rho_{F,G}^{(N)}$ is continuous. It satisfies
 $\rho_{F,G}^{(N)}(0) = 1/2$ and
$$\rho_{F,G}^{(N)}(s) \geq \frac{1}{2}- \frac{C(N)}{ p_N^{1/N}} s^{1/N},$$
for all $s\in [0; p)$, where $C(N)>0$ is a positive constant
depending only on $N$. }
\medskip

Let us note that a similar result for another class of pairs $(F,G)$
(defined by means of a symplectic quasi-state) has been proved in
\cite{EPR}. It would be interesting to find out whether such a lower
bound on the profile function is (asymptotically) exact.

It follows from \cite{EPR} that $pb_3^{(N)}(X,Y,Z) >0$ for all $N \geq 2$
provided the sets $X,Y,Z$ are superheavy
(see Section \ref{subseq-pb-qs} above).

One can similarly define  the natural analogue $pb_4^{(N)}$ of the $pb_4$-invariant in the context of iterated Poisson brackets, and repeat the proof of Theorem~\ref{thm-dichotomy}(ii)  to get a lower bound for the generalized profile function $\rho_{F,G}^{(N)}$. However, at the moment we have no tools for
proving the positivity of $pb_4^{(N)}$ in any example.
 }
\end{rem}

%%%%%%%%%%%%%%%%%%%%%%%%%%%%%%%%%%%%%%%%%%%%%%%%%%%%%%%%

\subsection{The two-dimensional case}\label{subsec-sympl-approx-pfs-2D}

In the two-dimensional case, the continuity of the profile function
at $0$ readily follows from the following result by Zapolsky.

\begin{prop}[\cite{Z}] \label{Prop:Zapolsky}
Let $(M,\omega)$ be a closed connected $2$-dimensional symplectic
manifold.  Let $ (F,G) \in \cF $ be a pair of functions with $ \|
\{ F,G \} \| \leqslant s$. Then there exist a pair of
Poisson-commuting functions $ (F',G') \in \cF$ with $ \| F-F'\|+
\|G-G'\| \leqslant C\sqrt{s}  $, where the constant $C$ depends
only on $(M,\omega)$.
\end{prop}

\medskip In other words, every almost commuting pair of functions
is nearly commuting, that is it can be approximated by a commuting
pair. Similar statements for various types of matrices and, more
generally, elements of $C^*$-algebras have been extensively
studied -- see e.g. \cite{PS,Hastings} and the references therein.
However, no analogue of Proposition~\ref{Prop:Zapolsky} is known
for higher-dimensional symplectic manifolds and there might be a
counterexample.

\medskip
\noindent {\bf Proof of
Proposition~\ref{prop-2-dim-case-almost-commut}.} Fix $\delta >0$
small enough. Take $(F_1,G_1) \in \cK_s$ with
$$d((F,G),(F_1,G_1)) \leq \rho_{F,G}(s) + \delta.$$
By Proposition~\ref{Prop:Zapolsky}, there exist Poisson-commuting
functions $F_2$ and $G_2$ with
$$d((F_1,G_1),(F_2,G_2)) \leq C\sqrt{s}.$$
By the triangle inequality,
$$\rho_{F,G}(0) \leq d((F,G),(F_2,G_2)) \leq  \rho_{F,G}(s) + C\sqrt{s}+
\delta$$ for every $\delta >0$. This yields inequality
\eqref{eqn-profile-fcn-2-dim-case}. \Qed

\subsection{Sharpness of the convergence rate: an example}
\label{subsec-sympl-approx-pfs-ex}

\noindent {\bf Proof of Theorem~\ref{thm-S2-estimate-exact}:}

We need to show that
$$\rho_{F,G} (\varepsilon)\leq \rho_{F,G} (0) -
C\sqrt{\varepsilon}$$ for some $C>0$ and any sufficiently small
$\varepsilon$ (since $\rho_{F,G}$ is non-increasing, by choosing a
smaller $C$ we can get the inequality for any $\epsilon$).

The standard symplectic form on the upper hemi-sphere can be
expressed as $ \omega = \frac{ dx \wedge dy }{ \sqrt{ 1-x^2-y^2 }
}$, while on the lower hemi-sphere we have $ \omega = - \frac{ dx
\wedge dy }{ \sqrt{ 1-x^2-y^2 } } $. Therefore, for a given pair
of functions $ f,g : S^2 \rightarrow \R $, on the upper
hemi-sphere we have
$$ \{ f(x,y),g(x,y) \}_{S^2} = \sqrt{1-x^2-y^2} \{ f(x,y),g(x,y)
\}_{\R^2 (x,y)},$$ while on the lower hemi-sphere we have
$$  \{ f(x,y),g(x,y) \}_{S^2} = - \sqrt{1-x^2-y^2} \{ f(x,y),g(x,y)
\}_{\R^2 (x,y)}.$$ In any case we have $$ | \{ f(x,y),g(x,y)
\}_{S^2} | = \sqrt{1-x^2-y^2}\ | \{ f(x,y),g(x,y) \}_{\R^2 (x,y)}
|.$$

\begin{figure}
 \begin{center}
\scalebox{0.75}{\includegraphics*{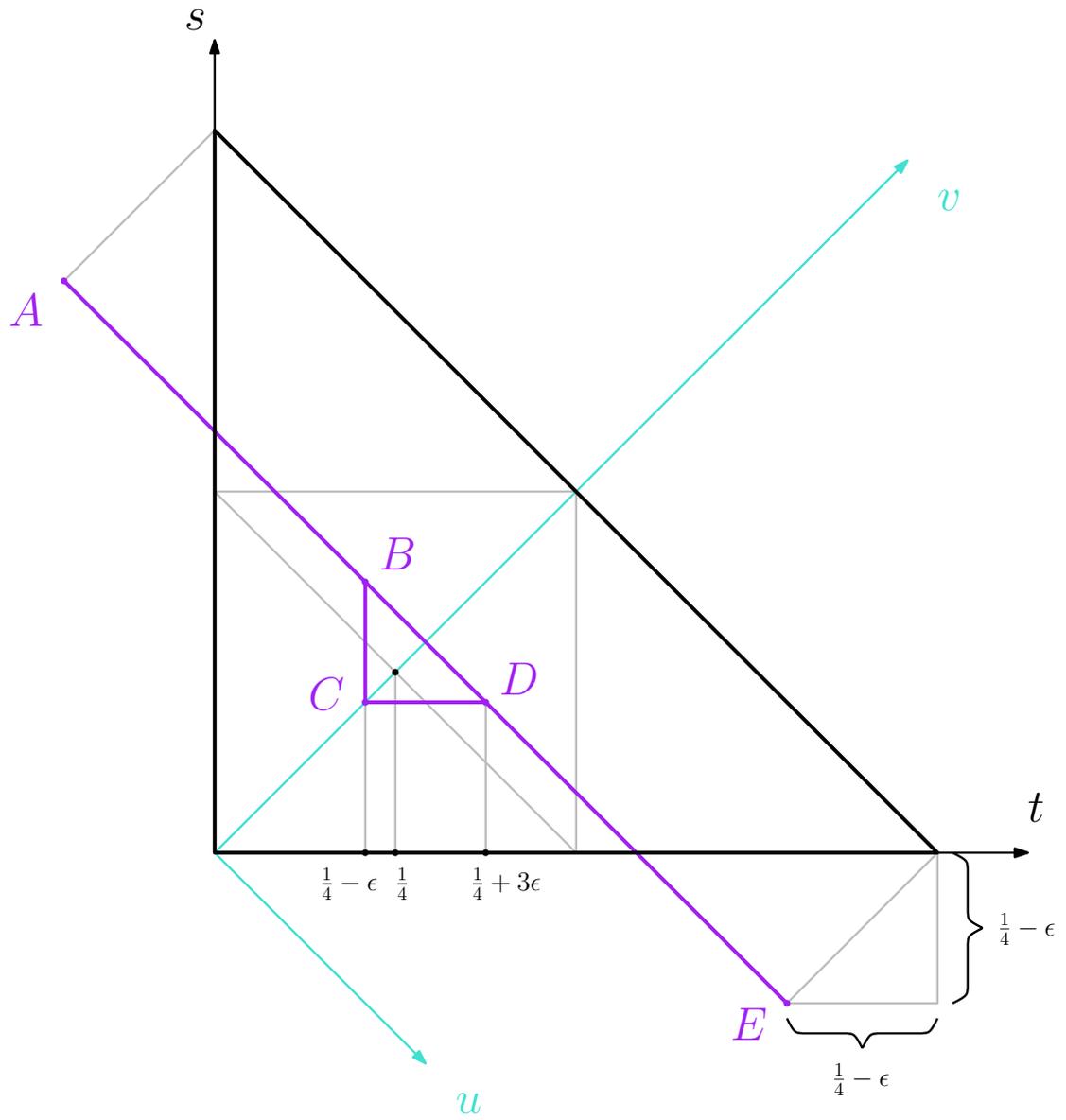}}
 \end{center}
  \caption{The image of $T$: $[ A ; E ] \cup \triangle BCD$}
   \label{fig5}
\end{figure}

Our purpose is to find smooth functions $ F_1 , G_1 : S^2
\rightarrow \R $, depending on a small parameter $ \varepsilon
> 0 $, such that $ \| F_1 - F \| + \| G_1 - G \| \leqslant 1/2 -
O(\varepsilon) $, while $ \| \{ F_1, G_1 \} \| \leqslant
O(\varepsilon^2) $. We will search for functions $ F_1,G_1 $ of
the form $ F_1 = f(x^2,y^2) $, $ G_1 = g(x^2,y^2) $,
where $ f,g : \bigtriangleup = \{ (t,s) | t,s \geqslant 0, t+s
\leqslant 1 \} \rightarrow \R $ are smooth functions. Further on
we use the notation $ t = x^2 $, $ s = y^2 $. We have
$$ \| F_1 - F \| + \| G_1 - G \| =  \| f(t,s) - t
\|_{\bigtriangleup} + \| g(t,s) - s \|_{\bigtriangleup}  ,$$
while
$$ |\{ F_1 , G_1 \}_{S^2}| = |\{ f(x^2,y^2),g(x^2,y^2) \}_{S^2}| $$ $$
= \sqrt{1-x^4-y^4}\ |\{ f(x^2,y^2),g(x^2,y^2) \}_{\R^2 (x,y)}|
$$
$$ = 4 \sqrt{1-x^4-y^4}\ |xy| |\{ f(t,s),g(t,s) \}_{\R^2 (t,s)}| $$
$$ = 4 \sqrt{(1-t^2-s^2)ts}\  |\{ f(t,s),g(t,s) \}_{\R^2 (t,s)}| .$$
For our purposes it is enough to find smooth $ f,g : \bigtriangleup
\rightarrow \R $ that satisfy $$ \| f(t,s) - t \|_{\bigtriangleup} +
\| g(t,s) - s \|_{\bigtriangleup} \leqslant 1/2 - O(\varepsilon) ,$$
$$ \| \{ f(t,s),g(t,s) \}_{\R^2 (t,s)} \| \leqslant O(\varepsilon^2)
.$$ Consider new coordinates $ u = t-s $, $ v = t+s $. In these
coordinates we have $ \bigtriangleup = \{ (u,v) | 0 \leqslant v
\leqslant 1 , -v \leqslant u \leqslant v \} $. We take the
functions $ f,g $ to be of the form
$$  f(u,v) = \phi(v) + u \psi(v) ,$$ $$ g(u,v) = \phi(v) - u \psi(v),$$
for some $\phi,\psi: [0,1]\to\R$,
 or, in regular coordinates $ (t,s)$, $$ f(t,s) = \phi(t+s) +
(t-s) \psi(t+s) ,$$
$$ g(t,s) = \phi(t+s) - (t-s)
\psi(t+s) .$$ We have
$$ \{ f,g \}_{\R^{2}(t,s)} = 2 \{ f,g
\}_{\R^{2}(u,v)} = 2 \{ \phi(v) + u \psi(v), \phi(v) - u \psi(v)
\}_{\R^{2}(u,v)} $$ $$ = - 4 \{ \phi(v) , u \psi(v)
\}_{\R^{2}(u,v)} = - 4 \psi(v) \{ \phi(v) , u \}_{\R^{2}(u,v)} = 4
\psi(v) \phi'(v) .$$ First, let us find a pair of {\em continuous}
functions $ \phi, \psi $, such that
\begin{equation}\label{eq-figure5} \| f(t,s) - t \|_{\bigtriangleup} , \| g(t,s) -
s \|_{\bigtriangleup} \leqslant 1/4 - \varepsilon.\end{equation}
The image of the corresponding map $T: (s,t) \mapsto
(f(s,t),g(s,t))$ consists of the union of a segment and a triangle
attached to it, see Figure \ref{fig5}.
Because of the symmetry, in order to verify \eqref{eq-figure5} it
is enough to check only that $ \| f(t,s) - t \|_{\bigtriangleup}
\leqslant 1/4 - \varepsilon $. We have $ f(t,s) - t = f(u,v) -
(u+v)/2 = \phi(v) + u\psi(v) - u/2 - v/2$. For a fixed $ v $ this
is a linear function of $ u $. Recall that $ \bigtriangleup = \{
(u,v) | 0 \leqslant v \leqslant 1 , -v \leqslant u \leqslant v \}
$. As a conclusion, it is enough to check the inequality $ |f(u,v)
- (u+v)/2| = | \phi(v) + u\psi(v) - u/2 - v/2 | \leqslant 1/4 -
\varepsilon $ only for the cases $ u = v $ and $ u = -v $ while $
0 \leqslant v \leqslant 1 $. Substituting $ u=v $, $ u=-v $ we see
that it is enough to check that $$ | \phi(v) + v\psi(v) - v |
\leqslant 1/4 - \varepsilon , $$ $$ | \phi(v) - v\psi(v) |
\leqslant 1/4 - \varepsilon $$ for $ 0 \leqslant v \leqslant 1 $.
We define our continuous $ \phi, \psi $ to be
$$ \phi(v) = (1/4 - \varepsilon) + 4\varepsilon v \text{ for
} v \in [0;1/2] ,$$
$$ \phi(v) = 1/4 + \varepsilon \text{ for } v \in [1/2;1] ,$$
and
$$ \psi(v) = 4\varepsilon \text{ for } v \in [0;1/2] ,$$
$$ \psi(v) = (1 - 4\varepsilon) + \frac{1}{v} ( -1/2 + 4\varepsilon ) \text{ for } v \in [1/2;1] .$$
Because of our choice of the functions $ \phi,\psi $, the
functions $ \phi(v), v\psi(v) $ are linear on each one of
intervals $ [0;1/2] $ and $ [1/2;1] $, and hence the functions  $
\phi(v) + v\psi(v) - v ,\phi(v) - v\psi(v) $ are linear on the
intervals $ [0;1/2] $ and $ [1/2;1] $ as well. Therefore it is
enough to check that $$ | \phi(v) + v\psi(v) - v | \leqslant 1/4 -
\varepsilon$$ and
$$ | \phi(v) - v\psi(v) | \leqslant 1/4 -
\varepsilon $$ only for $ v=0, 1/2,1 $. We have
$$ \phi(0) + 0 \cdot \psi(0) - 0 = \phi(0) = 1/4 - \varepsilon ,$$
$$ \phi(1/2) + 1/2 \psi(1/2) - 1/2 = (1/4 + \varepsilon) + 4\varepsilon/2 -
 1/2 = - 1/4 + 3\varepsilon,$$
$$ \phi(1) + 1 \cdot \psi(1) - 1 = (1/4 + \varepsilon) + 1/2 - 1 = - 1/4 + \varepsilon,$$
$$ \phi(0) - 0 \cdot \psi(0) = \phi(0) = 1/4 - \varepsilon ,$$
$$ \phi(1/2) - \psi(1/2)/2 = (1/4 + \varepsilon) - 4\varepsilon/2   = 1/4 - \varepsilon,$$
$$ \phi(1) - 1 \cdot \psi(1) = (1/4 + \varepsilon) - 1/2 = - 1/4 + \varepsilon .$$
In all the cases the absolute value of the result is not bigger
than $1/4 - \varepsilon $.

Hence we have found continuous $ \phi, \psi : [0;1] \rightarrow \R
$ for which
$$ \| f(t,s) - t \|_{\bigtriangleup}, \|
g(t,s) - s \|_{\bigtriangleup} \leqslant 1/4 - \varepsilon .$$ One
can easily make the functions $ \phi, \psi $ smooth by a
sufficiently $C^0$-small perturbation so that we will still have $
\| f(t,s) - t \|_{\bigtriangleup}, \| g(t,s) - s
\|_{\bigtriangleup} \leqslant 1/4 - \varepsilon/2$, and, moreover,
$$ | \phi'(v) | \leqslant 4 \varepsilon \text{ for } v \in [0;1/2] ,$$
$$ | \phi'(v) | \leqslant 32 \varepsilon^2 \text{ for } v \in [1/2;1] ,$$
$$ | \psi(v) | \leqslant 4 \varepsilon \text{ for } v \in [0;1/2] ,$$
$$ | \psi(v) | \leqslant 1/2 \text{ for } v \in [1/2;1] .$$
Then for any $ v \in [0;1] $ we have $ | 4 \psi(v) \phi'(v) |
\leqslant 64 \varepsilon^2 $. As a conclusion, we obtain
$$ \| f(t,s) -
t \|_{\bigtriangleup} + \| g(t,s) - s \|_{\bigtriangleup} \leqslant
(1/4 - \varepsilon/2) +  (1/4 - \varepsilon/2) = 1/2 -
\varepsilon,$$ and $$ | \{ f,g \}_{\R^{2}(t,s)} | = | 4 \psi(v)
\phi'(v) | \leqslant 64 \varepsilon^2 $$ at any point $ (t,s) \in
\bigtriangleup $. \Qed

\medskip
\begin{rem} {\rm At the moment we are unable to decide
whether the example constructed above has a counterpart in the
context of matrix algebras (for instance, for $su(n)$).}
\end{rem}

%%%%%%%%%%%%%%%%%%%%%%%%%%%%%%%%%%%%%%%%%%%%%%%%%%

\section{Detecting Hamiltonian chords}
\label{sec-chords-proofs}

\subsection{Proofs of the results about Hamiltonian chords}

In this section we prove
Theorems~\ref{thm-eq-pb-T},\ref{thm-non-autonomous} and part (ii) of
Corollary~\ref{cor-non-autonomous-superheavy}.

\medskip
\noindent {\bf Proof of Theorem~\ref{thm-eq-pb-T}.}
Let $X_0,X_1 \subset M$ be disjoint compact subsets, and let $G$ be
a function from $C^{\infty}_c(M)$. Set
$$Pb(X_0,X_1;G) := \inf ||\{F,G\}||,$$
where the infimum is taken over all functions $F\in C^\infty_c (M)$
with $F|_{X_0} \leq 0, F|_{X_1} \geq 1$, and
$$Pb' (X_0,X_1;G) := \inf ||\{F,G\}||,$$
where the infimum is taken over all functions $F\in C^\infty_c (M)$
with $F|_{X_0} = 0, F|_{X_1} = 1$.

 Put
$T=T(X_0,X_1;G)$ and let $F$ be any smooth compactly supported
function with
$$F|_{X_0} \leq 0, F|_{X_1} \geq 1.$$ There exist $i\in \{0;1\}$
and $x \in X_i$ so that $g_T x \in X_{1-i}$ (recall that $g_t$ is
the Hamiltonian flow of $G$). Thus $|F(g_Tx)-F(x)| \geq 1$, which
yields $||\{F,G\}|| \geq T^{-1}$. Therefore
$$T(X_0,X_1;G)\geq \frac{1}{Pb(X_0,X_1;G)}.$$
Since, obviously, $Pb(X_0,X_1;G)\leq Pb'(X_0,X_1;G)$, it remains
to prove that
$$T(X_0,X_1;G)\leq \frac{1}{Pb'(X_0,X_1;G)}.$$
We shall need the following lemma.

\begin{lemma}\label{lem-one-flow-exist-fcn} Let $v$ be a smooth
compactly supported vector field on $M$ and $X_0,X_1$ be a pair of
disjoint compact subsets of $M$. Denote by $g_t$ the flow of $v$.
Assume that $g_t X_0 \cap X_1 = \emptyset$ for all $t \in [-a;a]$
for some $a>0$. Then there exists a smooth compactly supported
function $F: M \to [0;1]$ such that $F|_{X_0} =0, F|_{X_1} =1$ and
$||L_v F|| < 1/a$.
\end{lemma}

\medskip
\noindent{\bf Proof of Lemma~\ref{lem-one-flow-exist-fcn}.} Choose
$b>a$, sufficiently close to $a$, so that $g_t X_0 \cap X_1 =
\emptyset$ for all $t \in [-b;b]$. The sets
$$\check{X}_0 := \bigcup_{t \in [0;b]} g_t(X_0)$$
and
$$\check{X}_1 := \bigcup_{t \in [0;b]} g_t(X_1)$$
do not intersect. Take any  smooth compactly supported function $H:
M \to [0;1]$ so that $H=0$ on $\check{X}_0$ and $H=1$ on
$\check{X}_1$. Put
$$F:= \frac{1}{b}\int_0^b H \circ g_t \; dt.$$
Clearly, $F$ has values in $[0;1]$, is compactly supported, $F=0$ on
$X_0$, $F=1$ on $X_1$ and
$$L_v F = \frac{1}{b} \int_0^b \frac{d}{dt} H \circ g_t \; dt =
\frac{1}{b}(H \circ g_b -H).\;$$ It follows that $||L_v F||\leq 1/b
< 1/a$, and we are done. \Qed

\medskip
Let us return to the proof of the theorem. Put $v =\sgrad G$.
Assume on the contrary that
$$T(X_0,X_1;G) > \frac{1}{Pb'(X_0,X_1;G)}.$$
Thus there exists
$$a > \frac{1}{Pb'(X_0,X_1;G)}$$
so that $g_t X_0 \cap X_1 = \emptyset$ for all $t \in [-a;a]$. By
Lemma~\ref{lem-one-flow-exist-fcn}, there exists a  smooth
compactly supported function $F:M \to [0;1]$ so that $F|_{X_0} =0,
F|_{X_1} =1$ and $||L_v F|| < 1/a$. But $L_v F = \{F,G\}$  and we
conclude that $Pb'(X_0,X_1,Y_0,Y_1)< 1/a$, which means a
contradiction. This completes the proof. \Qed

\medskip\noindent{\bf Proof of Theorem~\ref{thm-non-autonomous}.}
Choose $a >0$ so that
\begin{equation}
\label{eq-oscG} \max G - \min G < R -a. \end{equation} Let $u: \R
\to [0;+\infty)$ be a cut-off function which is equal to $1$ on
the interval $[-(R-a);R-a]$ and whose support lies in $(-R;R)$.
Consider a new autonomous compactly supported Hamiltonian
$$H: M \times \A_R \to \R,\; (x,r,\theta) \to u(r)(G(x,\theta)+r)\;$$
generating a Hamiltonian flow $h_t$. Since $H \leq 0$ on $\s_R Y_0$
and $H \geq 1$ on $\s_R Y_1$, Theorem~\ref{thm-chords1} guarantees
existence of a point $z=(y,0,\theta_0) \in \s_R X_0$ and $T \in
[-1/p;1/p]$ so that $h_{T} z \in \s_R X_1$.

We claim that the piece of trajectory $ z(t) = \{h_t z\}$, $t \in
[0;T]$ is entirely contained in the domain $V = \{|r| < R-a\}
\subset M \times \A_R$. Indeed, assuming the contrary, choose $\tau
\in [0;T]$ so that $h_{\tau}z \in \partial V$. Write $z(t) =
(x(t),r(t),\theta(t))$. We have that  $r(0)=0$ and $r(\tau) = \pm
(R-a)$. By the energy conservation law, $H(z(0)) = H(z(\tau))$ and
hence
$$G(z(0),\theta(0)) = G(z(\tau),\theta(\tau)) \pm (R-a).$$
This contradicts assumption \eqref{eq-oscG} and the claim follows.

It follows that $u(r(t))=1$ for all $t \in [0;T]$. Hence the
projection of $z(t)$ to $M$ is a curve $\alpha$ of the form
$\{g_{\theta_0+t} g_{\theta_0}^{-1}y\}$, $t \in [0;T]$. Put $x=
g_{\theta_0}^{-1}y$, $t_0 = \theta_0$ and $t_1=\theta_0 +T$. We see
that $g_{t_0}x =y \in X_0$ and $g_{t_1}x = x(T) \in X_1$. Thus
$\alpha$ is a required Hamiltonian chord. \Qed

\medskip
\noindent{\bf Proof of
Corollary~\ref{cor-non-autonomous-superheavy}, part (ii).} Choose $R >
\max G -\min G$. Identify the annulus $A_R$ with the sphere $S^2$
of area $2R$ with punctured the North and the South Poles. Under
this identification the zero section $\{r=0\}$ corresponds to the
equator, say $E$, of the sphere. Thus we consider $M \times \A_R$
as a domain in $M \times S^2$. The latter manifold is equipped
with the symplectic form $\omega + 2R\sigma$, where $\sigma$ is
the standard area form on $S^2$ of the total area $1$.
The $S^2$-stability of the quasi-state $\zeta$ on $(M,\omega)$ yields a quasi-state $\tzeta_{2R}$ on $(M \times S^2, \omega
+2R\sigma)$. Denote by $K_{2R}$ the constant from the
PB-inequality for $\tzeta_{2R}$.

Assume that the sets $X_0\cup Y_0, Y_0\cup X_1, X_1\cup Y_1, Y_1\cup X_0$ are superheavy.
Due to the
$S^2$-stability of $\zeta$, the sets
$$(X_0 \cup Y_0) \times E , (Y_0 \cup X_1)\times E, (X_1 \cup
Y_1)\times E, (Y_1 \cup X_0)\times E$$ are all superheavy with
respect to the quasi-state $\tzeta_{2R}$. By inequality
\eqref{eqn-inclusion} and Theorem~\ref{thm-qstates-pb4-dynamics},
\begin{multline*}
pb_4 (\s_R X_0, \s_R X_1,\s_R Y_0, \s_R Y_1) \\ \geq pb_4 (X_0
\times E, X_1 \times E,Y_0\times E, Y_1\times E) \geq
\frac{1}{4 K_{2R}}.\end{multline*}

Finally, the existence of the required
Hamiltonian chord follows now from Theorem~\ref{thm-non-autonomous}.
This finishes the proof.
\Qed

\begin{rem}
{\rm
Note that if we assume that $Y_1$ is superheavy (and so are $X_0\cup Y_0$ and $X_1\cup Y_0$), then the constant appearing in the previous proof can be improved
from ${1}/({4 K_{2R}})$ to ${1}/{K_{2R}}$.
Indeed, by  the
$S^2$-stability of the quasi-state on $(M,\omega)$, the sets $(X_0\cup Y_0)\times E, (X_1\cup Y_0)\times E, Y_1\times E$
are superheavy with respect to $\tzeta_{2R}$.
Then, by \eqref{eq-Delta-2},
\[
pb_4 (\s_R X_0, \s_R X_1,\s_R Y_0, \s_R Y_1) \geq  \frac{1}{K_{2R}}.
\]

}
\end{rem}

\subsection{Miscellaneous remarks}

Let us make a few more remarks on the interplay between
superheaviness and Hamiltonian chords for autonomous Hamiltonians.
{\bf In this section we assume that $M$ is closed.}

\medskip
\noindent
\begin{rem} [Recurrence of Hamiltonian chords] {\rm
Let $F$ be a smooth function on $M$. Denote its Hamiltonian flow by
$f_t$. Put $Y_0 = \{F \leq 0\}$ and $Y_1 = \{F \geq 1\}$. A subset
$X$ is called a {\it ballast} if $X \cup Y_i$ is superheavy for
$i=0,1$. For instance, in Example~\ref{exam-S2xS2-qstates} above the
role of ballasts is played by the Lagrangian discs $C^v$.

Given two ballasts $X_0,X_1$, denote by $P \subset \R$ the set of
all $\tau$ such that $f_\tau X_0 \cap X_1 \neq \emptyset$. We
claim that Hamiltonian chords between $X_0$ and $X_1$ exhibit a
recurrent behavior in the following sense: The set $P$ intersects
every interval of time-length $8K$. Indeed, since $Y_i$ are
invariant under $f_t$, the image of a ballast under $f_t$ is again
a ballast. Take any $s \notin P$ so that $f_sX_0 \cap X_1
=\emptyset$. Thus the quadruple $(f_sX_0,X_1,Y_0,Y_1)$ satisfies
the assumptions of Corollary~\ref{cor-non-autonomous-superheavy}.
Hence there exists $t \in [-4K;4K]$ so that $f_{t+s}X_0 =
f_tf_s X_0$ intersects $X_1$, and the claim follows. } \end{rem}

\medskip\noindent \begin{rem}[Energy control] {\rm
Let us follow the notations and the set-up of the previous
example. Fix an interval $I=[a;b]$ with $0 \leq a < b \leq 1$ and
put $X_i^I = X_i \cap F^{-1}(I)$, $i=0,1$, where $X_i$, $i=0,1$,
are disjoint ballasts.

We claim there exists a Hamiltonian chord of $f_t$ of time-length
$8K/(b-a)$ which touches both $X_0^I$ and $X_1^I$.

Interestingly enough, this statement has a flavor of time-energy
uncertainty: we have to pay for the precision of our knowledge of
the energy level carrying a chord by an uncertainty in our
knowledge of the time interval on which the chord is defined.

To prove the claim put $Y_0' = \{F \leq a\}, Y_1' = \{F \geq b\}$.
One can deduce from Proposition~\ref{prop-expan} that
$$pb_4(X_0^I, X_1^I, Y_0',Y_1') = pb_4(X_0,X_1,Y_0',Y_1') \geq \frac{1}{4K}.$$ Put $F' = (F-a)/(b-a)$. Then $F' \leq 0$ on $Y_0'$ and
$F' \geq 1$ on $Y_1'$. Therefore, by Theorem~\ref{thm-chords1},
the Hamiltonian flow $f'_t$ of $F'$ admits a chord of time-length
at most $8K$ touching both $X_0^I$ and $X_1^I$. The claim
follows from the fact that $f_t = f'_{t(b-a)}$.}\end{rem}

\medskip\noindent \begin{rem}[Producing rigid subsets from flexible ones] {\rm
Let $X_0$, $Y_0$, $Y_1$ be subsets of $M$ so that $Y_0$ and $Y_1$
are disjoint, and $X_0 \cup Y_0, Y_1 \cup X_0$ are superheavy.
Take any Hamiltonian $G$ such that $G|_{Y_0} \leq 0,G|_{Y_1} \geq
1$ and denote by $g_t$ its Hamiltonian flow. Put
$$Z:= \bigcup_{t \in [-4K;4K]} g_t X_0.$$
Theorem~\ref{thm-qstates-pb4-dynamics} implies that $Z$ intersects
{\it every} superheavy subset $X_1 \subset M$ and hence exhibits a
``symplectically rigid" behavior. To illustrate this, assume in
addition that the quasi-state $\zeta$ is invariant under the
identity component $Symp_0$ of the symplectomorphism group of
$(M,\omega)$: this happens in all known higher-dimensional
examples. Let $r_0(M,\omega):=\sup r(B)$, where $r(B)$ is the
radius of a symplectically embedded open ball $B\subset M$ and the
supremum is taken over all balls $B$ whose complement contains a
superheavy subset. It follows that $Z$ cannot be mapped into any
symplectically embedded ball $B \subset M$ of radius $r <
r_0(M,\omega)$ by a diffeomorphism from $\Symp_0$. A somewhat
paradoxical point here is that $X_0$ itself could be absolutely
``flexible", e.g. a closed Lagrangian disc. Of course, the
Hamiltonian function $G$ as above is quite special, hence there is
no contradiction.}\end{rem}

\medskip
\noindent \begin{exam} {\rm Here we present a construction of
subsets $X_0,X_1,Y_0,Y_1$ satisfying the assumptions of
Theorem~\ref{thm-qstates-pb4-dynamics} (ii). Let $A_i$,
$i=1,2,3,4$, be four closed superheavy subsets such that no three
of them have a common point. Put $A_{ij} :=A_i \cap A_j$. Present
each $A_i$ as a union of closed subsets, $A_i = B_i \cup C_i$, so
that
$$B_i \cap C_i \cap A_{ij} = \emptyset \;\forall i,j,$$ and
\begin{itemize}
\item $A_{12}\cup A_{13} \subset B_1$, $A_{14} \subset C_1$;
\item $A_{23} \cup A_{24} \subset B_2$, $A_{21} \subset C_2$;
\item $A_{34} \subset B_3$, $A_{31} \cup A_{32} \subset C_3$;
\item $A_{41} \subset B_4$, $A_{42} \cup A_{43} \subset C_4$.
\end{itemize}
Put
$$X_0 := B_1 \cup C_2, Y_0 := B_2 \cup C_3, X_1 := B_3 \cup C_4, Y_1 :=
B_4 \cup C_1.$$ Obviously, the sets $X_0 \cup Y_0$, $Y_0 \cup
X_1$, $X_1 \cup Y_1$, $Y_1 \cup X_0$ contain superheavy sets $A_2$,
$A_3$, $A_4$, $A_1$ respectively. At the same time it is
straightforward to check that
$$X_0 \cap X_1 = Y_0 \cap Y_1 =\emptyset,$$ as required.

One can also construct a quadruple of sets satisfying the
assumptions of Theorem~\ref{thm-qstates-pb4-dynamics} from a
triple of superheavy sets. Namely, assume $A, B, C\subset M$ are
closed superheavy sets with $A\cap B\cap C=\emptyset$. Let $U$ be
an open neighborhood of $A\cap B$ such that $\overline{U}$ is
disjoint from $C$. Set $X_0:=\overline{U}\cap (A\cup B)$,
$X_1:=C$, $Y_0:=A\setminus  U $, $Y_1:=B\setminus U $. Then $X_0,
X_1, Y_0, Y_1$ satisfy the assumptions of
Theorem~\ref{thm-qstates-pb4-dynamics}. By Proposition
\ref{prop-expan}, formula \eqref{eq-pb3-pb4} and monotonicity,
$pb_4(X_0,X_1,Y_0,Y_1) \geq pb_3(A,B,C)$. Note that if $A,B,C$
are, for instance, superheavy Lagrangian submanifolds intersecting
transversally and $U$ is the complement of a sufficiently small
closed tubular neighborhood of $C$, the sets $Y_0$ and $Y_1$ given
by this construction are finite unions of small Lagrangian discs.}
\end{exam}

\medskip \noindent
\begin{rem}\label{rem-comparison-time}
{\rm Let us compare the bounds on the time-length of Hamiltonian
chords given by Theorem~\ref{thm-chords1}
 (the autonomous case) and Theorem~\ref{thm-non-autonomous} (the
non-autonomous case). We will compare the bounds for the case of
an autonomous Hamiltonian and for $R=+\infty$ (i.e. when both
estimates are applicable and there are no restrictions on the
oscillation of the Hamiltonian). We have seen in
\eqref{eq-pb-products} above that
\begin{equation}
\label{eqn-pb4-stab-comparison} 1/pb_4(X_0,X_1,Y_0,Y_1) \leq
1/pb_4(\s X_0,\s X_1,\s Y_0,\s Y_1).
\end{equation}
Thus for Hamiltonian chords of autonomous Hamiltonians the
``autonomous" bound from Theorem~\ref{thm-chords1} is a priori
better than the ``non-autonomous" one from
Theorem~\ref{thm-non-autonomous}. As it was mentioned above, the
``autonomous" bound is sharp and therefore whenever one has the
equality in \eqref{eqn-pb4-stab-comparison} the ``non-autonomous"
bound is sharp as well (see Theorem \ref{thm-surface-1} and
Proposition \ref{thm-surfaces-2} above for an example where the
equality in \eqref{eqn-pb4-stab-comparison} is actually reached). It
would be interesting to find out whether the bound on the
time-length of the Hamiltonian chord given by
Theorem~\ref{thm-non-autonomous} is {\it always} sharp. In other words, the question is whether
for any compact $X_0,X_1,Y_0,Y_1\subset M$, $X_0\cap X_1=Y_0\cap
Y_1=\emptyset$, one can find time-dependent Hamiltonians as in Theorem~\ref{thm-non-autonomous}
admitting Hamiltonian chords that connect $X_0$
and $X_1$ and have time-lengths arbitrarily close to $1/pb_4(\s
X_0,\s X_1,\s Y_0,\s Y_1)$. }
\end{rem}

\section{Poisson brackets and pseudo-holomorphic po\-ly\-gons}\label{sec-pb-pseudo-holo}

\subsection{Defining polygons}

Let $\D \subset \C$ be the unit disc. Take $k \geq 2$ pairwise
distinct points $z_0,\ldots,z_{k-1}$ on the unit circle in
$\partial \D$ in the counter-clockwise cyclic order (thus further
on we use the convention $(k-1)+1=0$ for the indices). They divide
the circle into $k$ arcs
$$a_0 = [z_{k-1};z_0], a_1 = [z_0;z_1],\ldots, a_{k-1} = [z_{k-2};z_{k-1}].$$
Let $\cL = (L_0,\ldots,L_{k-1})$ be a collection of Lagrangian
submanifolds in a symplectic manifold $(M,\omega)$. A {\it
parameterized $k$-gon} with boundary on $\cL$ is a smooth map
$\phi: \D \to M$ such that $\phi(a_i) \subset L_i$ for all $i$.
For the sake of brevity we shall often refer to the image
$\phi(\D)$ as to a $k$-gon with boundary on $\cL$ with edges
$\phi(a_i)$ and (cyclically oriented) vertices $\phi(z_i)$. The
$k$-gons are called triangles for $k=3$ and quadrilaterals for
$k=4$.

Denote by $\D(z_0,\ldots,z_{k-1})$ the unit disc in $\C$ with $k$
counter-clockwise cyclically ordered marked points
$z_0,\ldots,z_{k-1}$ on the boundary. The space $\cP^k$ of all such
discs is naturally identified with a subset of $(\partial \D)^k$.
The group $PU(1,1)$ acts on $\cD$ by holomorphic automorphisms, and
hence acts on $\cP^k$. Given an almost complex structure $J$ on
$(M,\omega)$ consider the set of all pairs $(z,\phi)$ where $z=
(z_0,\ldots,z_{k-1}) \in \cP^k$ and $\phi: \D(z_0,\ldots,z_{k-1})
\to M$ is a $J$-holomorphic parameterized $k$-gon with boundary on
$\cL$. Its quotient by the natural action of the group $PU(1,1)$ is
called the moduli space of $J$-holomorphic $k$-gons with boundary on
$\cL$ and is denoted by $\cM$ (with some extra decorations which
will be introduced later).

\subsection{A reminder on the Maslov
class}\label{subsec-pseudo-Maslov}

Let $\ell_1, \ell_2$ be a pair of Lagrangian subspaces in a
symplectic vector space $V$. Pick any compatible almost complex
structure $J$ on $V$ with $J\ell_1 =\ell_2$. Denote by
$\gamma_J(L_1,L_2)$ the path $e^{tJ}\ell_1$, $t \in [0;\pi/2]$, of
Lagrangian subspaces joining $\ell_1$ with $\ell_2$.

Let now $L_0,\ldots,L_{k-1}$ be a collection of Lagrangian
submanifolds of a symplectic manifold $(M^{2n},\omega)$ in general
position: every pair from this collection intersects transversally
and there are no triple intersections. Let $P $ be a $k$-gon whose
edges $e_i$ lie on $L_i$. Choose a parametrization $e_i(t)$ of the
edges yielding the cyclic orientation of the boundary of the
polygon. Denote by $v_{i,i+1}$ the vertex lying on $L_i \cap
L_{i+1}$, where the indices are taken modulo $k$.

Let $\Lambda M \to M$ be a canonical fibration whose fiber over a
point $x \in M$ is the Lagrangian Grassmannian $\Lambda_n (T_xM)$.
For every edge $e_i$ consider its canonical lift $\hat{e}_i(t) =
T_{e_i(t)} L_i$ to $\Lambda M$. Fix an $\omega$-compatible almost
complex structure $J$ on $M$. The curves
$$\hat{e}_0, \gamma_J(T_{v_{0,1}}L_0
, T_{v_{0,1}}L_1), \hat{e}_1, \gamma_J(T_{v_{1,2}}L_1,
T_{v_{1,2}}L_2),\ldots,$$
$$\ldots \hat{e}_{k-1},
\gamma_J(T_{v_{k-1,0}}L_{k-1}, T_{v_{k-1,0}}L_0)$$ form a loop,
say $\theta$, in $\Lambda M$.

Take a symplectic trivialization of the tangent bundle $TM$ over
$P$ so that the restriction of $\Lambda M$ to $P$ splits as $P
\times \Lambda_n$, where $\Lambda_n$ is the Lagrangian
Grassmannian (the space of Lagrangian planes in the symplectic
vector space $\R^{2n}$). Write $\theta'$ for the projection of
$\theta$ to $\Lambda_n$.

Recall that $\cT_k$ is the set of homotopy classes of $k$-gons in
$M$ whose sides (in the natural cyclic order) lie, respectively,
in $L_0,L_1,\ldots, L_{k-1}$. Let $\alpha =[P] \in \cT_k$ be the
homotopy class of a $k$-gon $P$. By definition, the Maslov index
$m(\alpha)$ is the Maslov index of $\theta'$ in $\Lambda_n$. This
definition is independent of the choices of $J$, the specific
polygon $P$ inside the homotopy class $\alpha$ and the symplectic
trivialization. We refer to \cite{FOOO, Viterbo} for the details.

\subsection{Gluing polygons}\label{subsec-gluing}

Let $\cL = (L_0,\ldots,L_{k-1})$ be a collection of Lagrangian
submanifolds of a symplectic manifold $(M^{2n},\omega)$ in general
position. Given a homotopy class $\alpha$ of polygons with
boundary on $\cL$, we can perform two operations on it:
\begin{itemize}
\item Take a representative $P$ of $\alpha$ and attach a disc with
the boundary on some $L_i$ at a point lying on the $i$-th edge of
$P$; \item Attach a sphere at a point of $P$.
\end{itemize}
We say that two homotopy classes $\alpha$ and $\beta$ of the
polygons are {\it equivalent} if $\beta$ can be obtained from
$\alpha$ by a sequence of such operations. For brevity we shall
write
$$\beta = \alpha + \sum_{i=0}^{k-1} D_i + S,$$
where $D_i \in \pi_2(M,L_i)$ and $S \in \pi_2(M)$. Observe that
this representation is not unique: for instance, $S \in
\pi_2(M,L_i)$ for all $i$. The Maslov indices of $\alpha$ and
$\beta$ are related by the standard formula (cf. \cite{Viterbo})
\begin{equation}\label{eq-Maslov-sum}
m(\beta) = m(\alpha) + \sum m_{L_i}(D_i) + 2c_1(S),
\end{equation}
where $m_{L_i}$ is the Maslov class of $L_1$ and $c_1$ is the first
Chern class of $(M,\omega)$.

We shall need also another gluing operation. Let $P$ be a $k$-gon
with vertices $p_j \in L_j \cap L_{j+1}$, $j=0,\ldots,k-1$, with
boundary on $\cL$ and let $\alpha$ be a digon with vertices
$v,p_i$ and boundaries on $(L_i,L_{i+1})$ (the marked points $z_0,
z_1\in\partial \D$ are mapped, respectively, to $v$ and $p_i$;
accordingly, the arcs $a_0, a_1\subset \partial D$ are mapped,
respectively, into $L_{i+1}$ and $L_i$). Attaching $\alpha$ to $P$
along $p_i$ we get in a natural way a new $k$-gon $P'$ with
boundary on $\cL$ and the vertices
$$p_0,\ldots,p_{i-1},v,p_{i+1},\ldots,p_{k-1}.$$
We shall write $P'= \alpha \;\sharp\; P$. The homotopy class of $P'$
in $\cT_k$ does not depend on the specific choice of $P$ and
$\alpha$ within their homotopy classes. It will be denoted by
$[\alpha] \;\sharp\; [P]$. It is easy to check that
\begin{equation}\label{eq-Maslov-sharp}
m( [\alpha] \;\sharp\; P) = m([P])+m([\alpha])-n.
\end{equation}

\begin{prop}\label{prop-digon} Let $(L,K)$ be a finite type
collection of Lagrangian submanifolds. Let $\alpha$ be a digon with
boundaries on $(L,K)$ with the same vertices: $\alpha \in
\cT_2(a,a)$. Suppose that $m(\alpha)=n$. Then $\omega(\alpha)=0$.
\end{prop}

\medskip\noindent{\bf Proof.} Changing, if necessary, the orientation of
$\alpha$ we can assume that $\omega(\alpha) = c \geq 0$. Put
$\alpha_d:= \alpha \;\sharp\;\ldots\;\sharp\; \alpha$ taken $d$
times. Then, by \eqref{eq-Maslov-sharp}, we have that $m(\alpha_d)
=n$ while $\omega(\alpha_d)=dc$. Thus, by the finite type
condition, $dc$ is bounded as $d \to \infty$, and hence $c=0$.
\Qed

\subsection{Finite type collections of Lagrangian
submanifolds}\label{subsec-finitetype}

Here we discuss examples of finite type collections of Lagrangian
submanifolds.

\medskip
\noindent {\bf Proof of Proposition~\ref{prop-finite-type}.} Let
$\cL=(L_0,\ldots,L_{k-1})$ be a collection of monotone Lagrangian
submanifolds in general position with the same monotonicity
constant. Assume that for every $i$ the morphism $\pi_1(L_i) \to
\pi_1(M)$ has a finite image. We have to show that the collection
$\cL$ is of finite type. The latter assumption guarantees that
there exists only finite number of equivalence classes (in the
sense of Section~\ref{subsec-gluing} above) of  homotopy classes
of polygons with boundary on $\cL$. Suppose that $P \sim Q$ and
$m(P)=m(Q)$. Since all $L_i$ have the same monotonicity constant,
formula \eqref{eq-Maslov-sum} readily yields
$\omega(P)=\omega(Q)$. This, in turn, implies that $\cL$ is of
finite type. \Qed

\medskip
Consider now the cotangent bundle $T^*X$ of a closed manifold
$X^l$ equipped with the standard symplectic form $\sigma$. Let
$\cK=(K_0,\ldots,K_{k-1})$ be a collection of Lagrangian sections
of $T^* X$ of the form $K_i = \text{graph}\, dF_i$, where $F_i$ is
a smooth function on $X$. Suppose that $\cK$ is in general
position -- in particular, all functions $F_{i+1}-F_{i}$ are
Morse. Each intersection point $p \in K_i \cap K_{i+1}$ is a
critical point of $F_{i+1}-F_{i}$. Denote by $\nu(p)$ its Morse
index. One can readily check that for every polygon $P$ with
vertices $p_0,\ldots,p_{k-1}$ and boundary in $\cK$ one has $$
m(P) = \sum_i \nu(p_i),\ \sigma(P) = \sum_i (F_i(p_i) -
F_{i+1}(p_i)).$$ In particular,
\begin{equation}
\label{eq-sections} 0 \leq m(P) \leq kl, \;\;\;  |\sigma(P)| \leq
C(\cK):= 2k\cdot\max_i ||F_i||\;
\end{equation}
for {\it every} polygon $P$ with boundaries on $\cK$. This is a
considerable strengthening of the finite type property for the
collection $\cK$. In particular, it immediately yields the following
proposition.

\begin{prop}\label{prop-fintype-stab}
Let $(L_0,\ldots,L_{k-1})$ be any finite type collection of
Lagrangian submanifolds in a symplectic manifold $(M,\omega)$. Let
$K_i \subset T^*X$, $i=0,\ldots,k-1$ be sections as above. Then
the collection $(L_i \times K_i)_{i=0,\ldots,{k-1}}$ in $(M \times
T^*X, \omega \oplus \sigma)$ is of finite type with
\begin{multline}\label{eq-stab-precise}
A(L_0 \times K_0,\ldots,L_{k-1} \times K_{k-1};N) \leq\\
\leq \max \{A(L_0,\ldots,L_k;j) : j \in [N-kl;N]\}+C(\cK).
\end{multline}
\end{prop}

\medskip\noindent{\bf Proof.} Put $\Omega = \omega \oplus \sigma$.
Given a $k$-gon $P$ with boundary on $M \times T^*X$ and $m(P) =
N$, look at its projections $P_1$ and $P_2$ to $M$ and to $T^*X$
respectively. Then $m(P) = m(P_1) + m(P_2)$ and $\Omega (P) =
\omega(P_1) + \sigma(P_2)$. By \eqref{eq-sections},
$$-kl \leq m(P_1)-N \leq 0,\ |\omega(P_1)-\Omega(P)| \leq C(\cK).$$
 Thus
$$|\Omega(P)| \leq \max_j A(L_0,\ldots,L_k;j)+ C(\cK),$$
where $j$ runs over $[N-kl;N]$. Therefore the collection $(L_i
\times K_i)$ is of finite type and \eqref{eq-stab-precise} holds.
\Qed

\subsection{Preliminaries on Lagrangian Floer homology}\label{subsec-prelimLFH}

Here we sketch a definition of operations in Lagrangian Floer
homology (over $\Z_2$) -- the reader is referred to
\cite{DeSilva}, \cite{Seidel-book}, \cite{FOOO} for more details.

Let $(M^{2n},\omega)$ be a spherically monotone symplectic
manifold with a ``nice" behavior at infinity (e.g. geometrically
bounded \cite{ALP}). Let $\cL = (L_0,\ldots,L_{k-1})$ be a
collection of $k$ closed connected monotone Lagrangian
submanifolds, $k=2,3,4$. Our convention is that the indices of
$L_i$'s are taken modulo $k$, that is $L_k = L_0$, etc. Recall
that the minimal Maslov number $N_L$ of a Lagrangian submanifold
is the minimal positive generator of the image of $\pi_2(M,L)$
under the Maslov class. We put $N_L = +\infty$ if $\pi_2(M,L)=0$.

Throughout this section we shall assume that the following
conditions hold:
\begin{itemize}
\item[{({\bf F1})}] The whole collection $\cL$ is of finite type.
\item[{({\bf F2})}] Every pair $(L_i, L_{i+1})$ forms a collection
of finite type.
\item[{({\bf F3})}] The minimal Maslov number
$N_{L_i}$ of each $L_i$ is $\geq 2$.
\item[{({\bf F4})}] In case
$N_{L_i} = 2$, the number of pseudo-holomorphic discs of the
Maslov index $2$ passing through a generic point of $L_i$ is even.
In the terminology of \cite{FOOO} this means that the obstruction
class (over $\Z_2$) of each $L_i$ vanishes.
\end{itemize}
In addition we assume that $L_i$'s are in general position,
meaning that they intersect pairwise transversally and there are
no triple intersections, and if $k=4$, then
\begin{equation}
\label{eq-4inter-again} L_0 \cap L_2 = L_1 \cap L_3 = \emptyset.
\end{equation}

Consider the vector space $CF(L_i,L_{i+1}) := \text{Span}_{\Z_2}
(L_i \cap L_{i+1})$. Fix an $\omega$-compatible almost complex
structure $J$ on $M$. Given points $p_i \in L_{i-1} \cap L_{i}$,
$i=1,\ldots,k$, and a homotopy class $A \in \cT_k$ of $k$-gons
with boundary on $\cL$ and the vertices $p_1,\ldots,p_{k}$,
consider the moduli space $\cM_A(p_1,\ldots,p_k)$ of
$J$-holomorphic $k$-gons representing class $A$. A standard
transversality argument yields that for a generic $J$ this space
is a smooth manifold of the dimension
\begin{equation}
\label{exp-dim} \dim \cM_A(p_1,\ldots,p_{k})= m(A) + n(1-k) +
k-3.
\end{equation}

\medskip
\noindent \begin{rem} {\rm To make the transversality argument
actually work one needs to deal with a more involved version of
the $\bar{
\partial}$-equation (see \cite{Seidel-book}). We shall ignore this point in our sketch.
Furthermore, under certain assumptions there is a way to associate
an index, say $I(p)$, to each intersection point from $L_i \cap
L_{i+1}$ after equipping the Lagrangian submanifolds (and hence
the intersection points) with an additional structure of a
Lagrangian brane. In this case the dimension of the moduli space
$\cM_A(p_1,\ldots,p_{k})$ is given by a more standard expression
$$I(p_k) - \sum_{i=1}^{k-1}I(p_i)+ k-3\;$$ (see e.g.
\cite{Seidel-book}, formula (12.8) ). One can verify that
it coincides with \eqref{exp-dim}. We shall not enter the issue of
grading.}
\end{rem}

\medskip
We shall write $|Y|$ for the cardinality -- modulo 2 -- of a
finite set $Y$. Define a $\Z_2$-multi-linear map
$$\mu^{k-1}: CF(L_0,L_1) \otimes \ldots \otimes CF(L_{k-2},L_{k-1}) \to
CF(L_0,L_{k-1})$$ by
\begin{equation}\label{eq-Fukaya}
\mu^{k-1}(p_1,\ldots,p_{k-1}) = \sum_A \big{|}
\cM_A(p_1,\ldots,p_{k})\big{|} \cdot p_k,
\end{equation}
where the sum is taken over all $0$-dimensional moduli spaces.
Note that the moduli spaces $\cM_A(p_1,\ldots,p_{k})$ are
zero-dimensional (or empty) whenever $m(A)= n(k-1)-k+3$. Since our
collection is of finite type, the symplectic areas of all polygons
from such moduli spaces are bounded away from infinity. Thus a
compactness argument yields that the $0$-dimensional moduli spaces
are necessarily finite sets and that the sum in the right-hand
side of \eqref{eq-Fukaya} is finite.

The operation $\mu^1: CF(L_0,L_1)\to CF(L_0,L_1)$ is a
differential: $\mu^1 \circ \mu^1 =0$: this is guaranteed by Floer
gluing/compactness theorems and by the vanishing of the
obstruction class. For convenience we denote $\mu^1$ by $d$. The
corresponding homology $\text{Ker}\, d/\text{Im}\, d$ is called
the {\it Lagrangian Floer homology} $HF(L_0,L_1)$  of $L_0$ and
$L_1$. It is a $\Z_2$-module. In the same way we define Floer
homology $HF(L_i,L_j)$ for all $i,j$, and for the sake of brevity
use the same notation $d$ for the Floer differentials for all
$i,j$. Note that when $k=4$, the intersection condition
\eqref{eq-4inter-again} guarantees that
$HF(L_0,L_2)=HF(L_1,L_3)=0$.

Consider the  operation
$$\mu^2 : CF(L_0,L_1) \otimes CF(L_1,L_2) \to CF(L_0,L_2).$$
We shall abbreviate $\mu^2(a_1,a_2) = a_1 a_2$. This operation
satisfies the Leibnitz rule
$$d(a_1a_2) =a_1\cdot da_2 + da_1 \cdot a_2$$ and hence descends to
an operation in homology:
$$HF(L_0,L_1) \otimes HF(L_1,L_2) \to HF(L_0,L_2).$$
The latter is called the {\it triangle (or Donaldson) product} in
Lagrangian Floer homology. If $k=4$, we define the triangle
product for the triple $L_1,L_2,L_3$ in the same way and keep for
it the same notation. Note that for $k=4$ the intersection
condition \eqref{eq-4inter-again} guarantees that for pairwise
distinct $i,j,l$ the triangle product
 $$CF(L_i,L_j) \otimes CF(L_j,L_l) \to CF(L_i,L_l)$$
 vanishes already on the chain level.
The operation
$$\mu^3 : CF(L_0,L_1) \otimes CF(L_1,L_2) \otimes CF(L_2,L_3) \to CF(L_0,L_2)$$
satisfies the $A_{\infty}$-relation
\begin{multline}\label{eq-mu3-ainfty}
d\mu_3(a_1,a_2,a_3) = \mu^3(da_1,a_2,a_3)+ \mu^3(a_1,da_2,a_3) + \mu^3(a_1,a_2,da_3) \\
+ a_1(a_2a_3) + ( a_1a_2)a_3.
\end{multline}
This formula yields two useful facts. First, assume that $\mu^3
=0$. Then the triangle product is associative: $a_1(a_2a_3) =
(a_1a_2)a_3$. Second, we have the following proposition:

\begin{prop}\label{prop-mu3-correctness}
Assume that $L_0 \cap L_2 = L_1 \cap L_3 = \emptyset$. Then $\mu^3$
descends to an operation in Lagrangian Floer homology
$$ HF(L_0,L_1) \otimes HF(L_1,L_2) \otimes HF(L_2,L_3) \to HF(L_0,L_2).$$
\end{prop}

\medskip
By a slight abuse of notation, we shall still denote the
homological operation by $\mu^3$.

\medskip\noindent{\bf Proof.} The assumption on intersections yields that the
product $\mu^2$ vanishes for every triple $L_i,L_{i+1},L_{i+2}$.
Thus the terms $a_1(a_2a_3)$ and $( a_1a_2)a_3$ in
\eqref{eq-mu3-ainfty} vanish, which immediately yields the statement
of the proposition. \Qed

It is a folkloric fact that the Lagrangian Floer homology and the
operations introduced above remain invariant under exact
Lagrangian isotopies of the submanifolds $L_i$ (of course, in case
of the $\mu^3$-operation on homology one needs the intersection
assumption $L_0 \cap L_2 = L_1 \cap L_3 = \emptyset$ to remain
valid during the isotopies). We are going to discuss a particular
case of this statement in a slightly different language: instead
of deforming Lagrangian submanifolds we shall deform the
symplectic form on $M$. A crucial feature of this setting which
significantly simplifies the analysis is that the intersection
points from $L_i \cap L_{i+1}$ remain fixed and transversal in the
process of the deformation.

Consider a deformation $\omega_s$, $s \in [0;1]$, $\omega_0 =\omega$,
of the symplectic
form $\omega$ through symplectic forms on
$M$ which satisfy the following conditions:
\begin{itemize}
\item[{({\bf D1})}] $\omega_s=\omega$ near each $L_i$ for all $s$;
\item[{({\bf D2})}] $\omega_s$ is cohomologous to $\omega$ for all
$s$; \item[{({\bf D3})}] for any $s\in [0;1]$ and any
$i=0,\ldots,k-1$ the integrals of the forms $\omega_s$ and
$\omega$ over discs define the same functional $\pi_2(M,L_i\cup
L_{i+1})\to \R$.
\end{itemize}
Note that $L_i$ is a monotone Lagrangian submanifold of
$(M,\omega_s)$ and its monotonicity constant does not depend on $s$.
Furthermore, assumptions (F2)-(F4) hold automatically for all $s$.
We shall assume in addition that
\begin{itemize}
\item[{({\bf D4})}] the collection $\cL$ is of finite type with
respect to $\omega_s$ for all $s \in [0;1]$.
\end{itemize}

Choose a generic $1$-parametric family $J_s$, $s\in [0;1]$, of
$\omega_s$-compatible almost complex structures. Note that the
vector spaces $CF(L_i,L_{i+1})$ do not depend on $s$. Write $d_s$
for the Floer differential on $CF(L_i,L_j)$ with $i \neq j$.  Denote
by $HF_s(L_i,L_j)$ the Lagrangian Floer homology, and by
$\mu^{k-1}_s$ the operations associated to the collection $\cL$. We
shall write $\cM^s$ for the moduli space of $J_s$-holomorphic
polygons with boundaries on $\cL$ and $\cM_B^s
(p_1,\ldots,p_k)\subset \cM^s$ for the space of $\cJ_s$-holomorphic
polygons in a homotopy class $B$ with the vertices $p_1,\ldots,p_k$.

\begin{prop}\label{prop-LFHinvariance} Let $k=3$ or $4$.
There exist isomorphisms $$\phi_i: HF_0(L_i,L_{i+1}) \to
HF_1(L_i,L_{i+1}),\;\; i=0,\ldots,k-2,$$ and $\bar{\phi}_{k-1}:
HF_0(L_0,L_{k-1}) \to HF_1(L_0,L_{k-1})$ which send $\mu^{k-1}_0$
to $\mu^{k-1}_1$, i.e.
\begin{equation}
\label{eq-LFHinvariance}
\mu^{k-1}_1(\phi_0(x_0),\ldots,\phi_{k-2}(x_{k-2})) =
\bar{\phi}_{k-1}(\mu^{k-1}_0(x_0,\ldots,x_{k-2}))\end{equation}
for all $ x_i \in HF_0(L_i,L_{i+1}), \; i=0,\ldots,k-2$.
\end{prop}

\medskip
\noindent{\bf Proof.} Note that the differential $d_s$ and the
operations $\mu^{k-1}_s$ can change in the process of deformation
only due to bubbling-off. Since $L_i$'s are monotone with the
minimal Maslov number $\geq 2$, for a generic 1-parametric family
$J_s$ there is no bubbling-off of $J_s$-holomorphic discs and
spheres (and we assume that our 1-parametric family is chosen to
have this property). By the Gromov-Floer compactness result, other
possible degenerations of $J_s$-holomorphic polygons can be analyzed
by looking at possible degenerations of the disc $\D$ with the
marked points on the boundary into tree-like connected cusp-curves
with the marked points on them. Such an analysis, together with the
intersection assumptions $L_0 \cap L_1 \cap L_2 = \emptyset$ for
$k=3$ and \eqref{eq-4inter-again} for $k=4$, shows that the only
possible pattern of the bubbling-off is as follows: a
$J_s$-holomorphic digon $\beta$ with boundaries on some pair
$(L_i,L_j)$ of index $m(\beta)=n+1$ splits into the sum of two
digons $\beta = \beta' \;\sharp\; \alpha$ where $m(\alpha)=n$. This
splitting can take place for a finite set $T=\{0 < t_1 < \ldots <
t_N < 1\}$ of the parameter $s$ which we will call the {\it
critical} values. Thus on the intervals
$$[0;t_1),\ldots, (t_{l-1};t_l),\ldots, (t_N;1]$$
the Floer homology and the operations do not change and their
realizations for different values of the parameter (within such an
interval) will be identified. Without loss of generality, we can
assume that for every critical parameter $t \in T$ there is a
unique digon $\alpha$ with $m(\alpha)=n$. We shall call $\alpha$
{\it an exceptional digon}.

Fix a pair $(L,K)$ of distinct Lagrangian submanifolds from our
collection. Suppose that for $t \in T$ there exists an exceptional
digon $\alpha \in \cM^t_A(a,b)$, where $a,b \in L \cap K$. Note
that $\omega_t(\alpha)>0$ and hence Proposition~\ref{prop-digon}
above yields that $a \neq b$. Following Floer \cite[Lemma
3.5]{Floer-Morse-th-Lagr-int}, define an endomorphism $\psi^t$ of
$CF(L,K)$ by
\begin{equation}\label{eq-psi-def}
\psi^t(x) = x + (x,a)b,
\end{equation}
where $(x,a)$ is the coefficient at $a$ in the expansion of $x$ with
respect to the basis $L \cap K$ of $CF(L,K)$. Observe that $\psi^t
\circ \psi^t$ is the identity map (recall that we work over $\Z_2$)
and hence $\psi^t$ is an isomorphism. By using a gluing/compactness
argument Floer showed in \cite{Floer-Morse-th-Lagr-int} that for a
sufficiently small $\epsilon>0$
$$\psi^t \circ d_{t-\epsilon} = d_{t+\epsilon} \circ \psi^t.$$
Thus $\psi^t$ induces an isomorphism
$$\phi^t: HF_{t-\epsilon}(L,K) \to
HF_{t+\epsilon}(L,K).$$ Taking the composition of isomorphisms
$\phi^t$ over all critical parameters $t \in T$ we get an
isomorphism
\begin{equation}\label{eq-phikl}
\phi(L,K) :HF_0(L,K) \to HF_1(L,K).
\end{equation}
We claim that these isomorphisms send $\mu^{k-1}_0$ to
$\mu^{k-1}_1$. The proof is based on the very same Floer's argument.
Let us elaborate it in the case $k=4$ (the case $k=3$ is analogous).

Let us study what happens with  the operation $\mu^3_s$ when the
parameter $s$ passes a critical value $t \in T$. Let $\alpha \in
\cM^t(a,b)$ be the exceptional digon, and $A$ be its homotopy
class. We denote by $\cM^+$ and $\cM^-$ the moduli spaces of
$J_s$-holomorphic $k$-gons for  $s \in (t; t+ \epsilon)$ and $s
\in (t-\epsilon; t)$ respectively, and by $\mu^3_\pm$ the
corresponding $\mu^3$-operations.

\medskip
\noindent {\sc Case 1: $a,b \in L_0 \cap L_1$.} Consider a
0-dimensional moduli space of the form $\cM^s_B(b,p_2,p_3)$, where
$s \in (t-\epsilon,t+\epsilon)$, $p_2\in L_1\cap L_2$, $p_3\in
L_2\cap L_3$. It does not change when $s$ passes through the
critical value $t$. Take a $J_t$-holomorphic quadrilateral $P \in
\cM^t_B(b,p_2,p_3,q)$, $q\in L_0\cap L_3$, and look at the
quadrilateral $\alpha \;\sharp\; P$. A parametric version of the
standard compactness/gluing argument for pseudo-holomorphic
polygons yields that the following bifurcation takes place: there
exists a unique family of pseudo-holomorphic polygons from
$\cM^s_{A \;\sharp\; B}(a,p_2,p_3,q)$, where either $s \in
(t-\epsilon;t)$ or $s \in (t;t+\epsilon)$ but not both, which
bubbles off to  $\alpha \;\sharp\; P$ as $s=t$ and which
disappears as $s$ enters the other half of the interval
$(t-\epsilon; t+\epsilon)$. In other words, each $P $ contributes
$\pm 1$ to the difference
$$\big{|}\cM^+_{A\;\sharp\;B} (a,p_2,p_3,q)\big{|} -
\big{|}\cM^-_{A \;\sharp\; B}(a,p_2,p_3,q)\big{|}. $$ It follows
that
$$(\mu^3_+(a,p_2,p_3),q) - (\mu^3_-(a,p_2,p_3),q) =
(\mu^3_+(b,p_2,p_3),q),$$
(recall that we are counting modulo $2$), and hence
\begin{equation}
\label{eq-mu3-move1}
\mu^3_+(a,p_2,p_3)+\mu^3_+(b,p_2,p_3)=\mu^3_-(a,p_2,p_3).
\end{equation}
The cases when $a,b$ lie in $L_1 \cap L_2$ (respectively, in $L_2
\cap L_3$) yield similar equalities. The only difference with
\eqref{eq-mu3-move1} is that the points $a,b$ appear at the second
(respectively, at the third) position in $\mu^3_\pm$.

\medskip
\noindent {\sc Case 2: $a,b \in L_0 \cap L_3$.} Similarly, we look
at the broken quadrilateral $P \;\sharp\; \alpha$, where $P$ lies in
the $0$-dimensional moduli space $\cM^t_B(p_1,p_2,p_3,a)$, and
conclude that $P$ contributes $\pm 1$ to the difference
$$\big{|}\cM^+_{B\;\sharp\; A} (p_1,p_2,p_3,b)\big{|} - \big{|}\cM^-_{B
\;\sharp\; A}(p_1,p_2,p_3,b)\big{|}. $$ This yields (modulo $2$)
$$(\mu^3_+(p_1,p_2,p_3),b) - (\mu^3_-(p_1,p_2,p_3),b) =
(\mu^3_-(p_1,p_2,p_3),a).$$ For every $q \neq b$
$$(\mu^3_+(p_1,p_2,p_3),q)=(\mu^3_-(p_1,p_2,p_3),q)\;,$$
and hence (modulo 2)
\begin{equation}
\label{eq-mu3-move2} \mu^3_+(p_1,p_2,p_3)= \mu^3_-(p_1,p_2,p_3) +
(\mu^3_-(p_1,p_2,p_3),a)b.
\end{equation}

Suppose that the exceptional digon $\alpha$ is associated to the
pair $(L_u,L_v)$, where $(u,v)= (0,1),(1,2),(2,3)$ or $(0,3)$. Define
an  isomorphism $\psi^t_{ij}$ of $CF(L_i,L_j)$ by formula
\eqref{eq-psi-def} if $(i,j) = (u,v)$ and as the identity map
otherwise. Using formulas \eqref{eq-psi-def},\eqref{eq-mu3-move1}
and \eqref{eq-mu3-move2} we conclude that
\begin{equation}\label{eq-vspmu3}
\mu^3_+(\psi^t_{01}(x_0),\psi^t_{12}(x_1),\psi^t_{23}(x_2))=
\psi^t_{03}(\mu^3_-(x_0,x_1,x_2))
\end{equation}
for all
$$x_0 \in L_0 \cap L_1, x_1 \in L_1 \cap L_2, x_2 \in L_2 \cap L_3.$$
The composition of $\psi^t_{ij}$'s over all $t \in T$ is exactly the
isomorphism $\phi(L_i,L_j)$ introduced in \eqref{eq-phikl} above.
Put
$$\phi(L_0,L_1):=\phi_0,\; \phi(L_1,L_2):=\phi_1,\; \phi(L_2,L_3):=
\phi_2,\; \phi(L_0,L_3) = \bar{\phi}_3.$$  With this notation
formula \eqref{eq-vspmu3} readily yields \eqref{eq-LFHinvariance}.
This completes the proof of the proposition. \Qed

\subsection{The product formula}

Let $\cL= (L_0,\ldots,L_{k-1})$ be a generic collection of
Lagrangian submanifolds of a symplectic manifold $(M,\omega)$
satisfying assumptions (F1)-(F4) of Section~\ref{subsec-prelimLFH}
above and the intersection condition \eqref{eq-4inter-again}.
Choose a generic collection of $k$  sections
$\cK=(K_0,\ldots,K_{k-1})$ of $T^*S^1$. Assume that all $K_i$ are
exact, that is of the form $K_i = \text{graph}\, dF_i$ for some
functions $F_i : S^1 \to \R$. Consider a collection
$\widehat{\cL}:=(\widehat{L}_i := L_i \times K_i)$. It also
satisfies properties (F1) -(F4) and \eqref{eq-4inter-again}.
Indeed, (F1) and (F2) follow from
Proposition~\ref{prop-fintype-stab} and the remaining properties
readily follow from the definitions.

The K\"{u}nneth formula in Floer homology (which can be obtained
by considering
the Floer complexes for a split almost complex structure on
$M\times T^* S^1$)
yields
$$HF(\widehat{L}_i, \widehat{L}_{j}) = HF(L_i,L_j) \otimes HF(K_i,K_j).$$
With this  identification we have that
\begin{equation}\label{eq-prod-2}
(a \otimes A) \cdot (b \otimes B) = (ab) \otimes (AB)
\end{equation}
and
\begin{equation}\label{eq-prod-3}
 \mu^3(a \otimes A,b \otimes B,c \otimes C) = \mu^3(a,b,c) \otimes (ABC).
\end{equation}

It is well-known \cite{FO} that the $\Z_2$-module $HF(K_i,K_j)$ is
canonically identified with $H^1(S^1,\Z_2)$ so that the product
$\mu^2$ for $\cK$ corresponds to the cup-product and the
$\mu^3$-operation for $\cK$ vanishes. Let us mention that the
product  $\mu^2$ for $\cK$ is associative and hence the expression
$ABC$ is well-defined.

The conclusion of this discussion is that {\it the operations
$\mu^2$ and $\mu^3$ for $\widehat{\cL}$ do not vanish, provided
they do not vanish for $\cL$.}

The proof of \eqref{eq-prod-2} is straightforward and will be omitted.
The proof of \eqref{eq-prod-3} is a bit more delicate and will be
sketched below. For more information on the product formulae see
\cite{Amorim-PhD}.

\medskip
\noindent {\bf Sketch of the proof of formula \eqref{eq-prod-3}:}
Consider the space $\cP^4$ of all discs with four counterclockwise
cyclically ordered boundary points, and denote by
$\widetilde{\cP}^4$ its quotient by the natural action of $PU(1,1)$.
We shall denote by $\widetilde{P} \in \widetilde{\cP}^4$ the image
of $P \in \cP^4$ in $\widetilde{\cP}^4$.

\medskip
\noindent{\sc Step 1:} Fix the standard complex structure $I$ on
$T^*S^1$. Fix a generic almost complex structure $J$ on $M$. We
are studying $J \oplus I$-holomorphic maps $u$ from $P \in \cP$ to
$M \times T^*S^1$ with boundary on $\widehat{\cL}$. (In this
sketch we will not discuss the regularity of these almost complex
structures.)
%(One can check
%using the methods from \cite{Floer-unreg}, \cite{Oh-psh-discs},
%cf. \cite{DeSilva}, that $I$ is a regular almost complex
%structure, i.e. one for which the relevant moduli spaces are
%smooth manifolds of the expected dimension, and that so are $J$
%and $J\oplus I$ for a generic choice of $J$).

Each such map has
the form $u = (\phi,\psi)$, where $\phi: P \to M$ and $\psi: P \to
T^*S^1$. Using the dimension formula \eqref{exp-dim} and the fact that the
Maslov class is additive with respect to direct sums, we get that
the $0$-dimensional moduli space of such maps can
arise from two sources:
\begin{itemize}
\item[{(i)}] The map $\phi$ lies in the $0$-dimensional moduli
space of $J$-holomorphic quadrilaterals with boundary on $\cL$.
This picks a {\it finite subset}, say $Z$,  of possible classes
$\widetilde{P}$ in the space $\widetilde{\cP}^4$. To get a generic
existence of an $I$-holomorphic map $\psi:P \to T^*S^1 $ with
boundary on $\cK$ so that $[P] \in Z$, the map $\psi$ must lie in
a $1$-dimensional component of the moduli space of $I$-holomorphic
quadrilaterals with boundary on $\cK$ -- in this case by varying
$\psi$ we can ``tune in" its source to be in $Z$.

\item[{(ii)}] The same, but with $\phi$ lying in the
$1$-dimensional moduli space and $\psi$ lying in the
$0$-dimensional moduli space.
\end{itemize}
Note that the count of $0$-dimensional moduli spaces of
pseudo-holomorphic quadrilaterals yields the $\mu^3$-operation.
Since the latter vanishes for $\cK$, the scenario (ii) can be
disregarded. Thus we shall focus on (i) and study
$1$-dimensional components of the moduli space of $I$-holomorphic
quadrilaterals with boundary on $\cK$.

\begin{figure}
 \begin{center}
\scalebox{0.7}{\includegraphics*{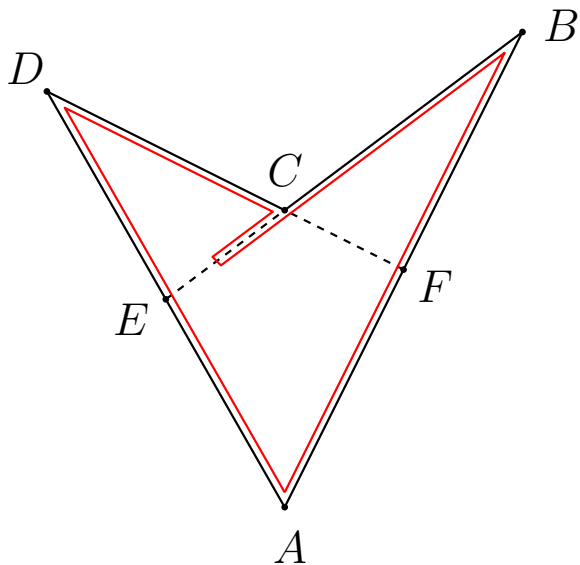}}
 \end{center}
  \caption{Non-convex quadrilateral $ABCD$}
   \label{fig3}
\end{figure}

\medskip
\noindent{\sc Step 2:} Pass to the universal cover $\R^2:= \C \to
T^*S^1$ and lift the sections $K_i$ (we keep the same notation for
the lifts). Look at the holomorphic quadrilaterals formed by
$K_0,K_1,K_2,K_3$. The holomorphic quadrilaterals of expected
dimension $1$ correspond to {\it embedded} quadrilaterals with
boundary on $\cK$ which have a unique interior angle $> \pi$. Fix
such a quadrilateral and, to make further analysis more
transparent, draw it as a non-convex Euclidean quadrilateral
$ABCD$ in $\R^2$, where the vertices are written in the
counter-clockwise order and the angle at $C$ is $>\pi$. Introduce
also the points $E$, which is the intersection of the edge $AD$
with the ray $[BC)$, and $F$, which is the intersection of the
edge $AB$ with the ray $[DC)$, see Figure \ref{fig3}.
 Suppose that
$$AB \subset K_1, \; BE \subset K_2,\; FD \subset K_3,\; DA \subset
K_0.$$ Introduce a parameter $t \in [0;1]$ on the broken line
$ECF$ so that $E$ corresponds to $t=0$, $C$ corresponds to
$t=1/2$, $F$ corresponds to $t=1$. Denote by $X_t$ the point on
the broken line $ECF$ corresponding to the value $t$ of the
parameter. Look at the following family of closed broken lines
which depends on a parameter $t\in (0,1)$:
\begin{itemize}
\item the line formed by the segments $AB \subset K_1$,
$BX_t\subset K_2$, $X_tC \subset K_2$, $CD \subset K_3$, $DA
\subset K_0$ for $t \in (0;1/2)$; \item the line formed by the
segments $AB \subset K_1$, $BC \subset K_2$, $CD \subset K_3$, $DA
\subset K_0$ for $t =1/2$; \item the line formed by the segments
$AB \subset K_1$, $BC \subset K_2$, $CX_t \subset K_3$, $X_tD
\subset K_3$, $DA \subset K_0$ for $t \in (1/2;1)$.
\end{itemize}

For every $t \in (0;1)$ this broken line bounds a holomorphic
polygon, say $\psi_t: P_t \to \C$, where $P_t =
\D(z_0,z_1,z_2,z_3)$, so that
$$\psi_t(z_0)=A,\; \psi_t(z_1) = B,\; \psi_t(z_2)=C,\; \psi_t(z_3) =
D.$$ Note that for $t \neq 1/2$ the map $\psi_t|_{S^1}$ hits $C$
twice, so $z_2$ corresponds to the second hit for $t < 1/2$ and to
the first hit for $t > 1/2$. By applying a M\"{o}bius
transformation we can assume that $z_0= -i, z_1 = 1, z_3=-1$ and
$z_2$ varies with $t \in (0;1)$ between $1$ and $-1$ (excluding
the endpoints themselves) in the upper half-circle.

Next, we wish to analyze the behavior of these holomorphic
quadrilaterals when $t \searrow 0$ and $t \nearrow 1$. For this
purpose let us recall (see \cite{FO}) that the Deligne-Mumford
compactification of $\widetilde{\cP}^4$ can be identified with
$[0;1]$, where the boundary point $0$ corresponds to the stable
curve $$\Sigma_0:= \D^-(z_0,z_1,z_*) \;\sharp\;
\D^+(z_*,z_2,z_3)$$ and the boundary point $1$ corresponds to the
stable curve $$\Sigma_1:= \D^-(z_0,z_*,z_3) \;\sharp\;
\D^+(z_*,z_1,z_2).$$ Here we denote by $\D^\pm$ two copies of the
unit disc.

When $t \searrow 0$, the bubbling-off happens: The map $\psi_t$
converges to a map $\psi_0= \psi_0^-\;\sharp\;\psi_0^+ $ from
$\Sigma_0$ to $\C$. Here $\psi_0^-: \D^- \to \C$ is a holomorphic
triangle with the three sides, respectively, on $K_0,K_1,K_2$ and
the vertices
$$\psi_0^-(z_0) = A, \psi_0^-(z_1) = B, \psi_0^-(z_*) =E,$$
and $\psi_0^+: \D^+ \to \C$ is a holomorphic triangle with the
three sides, respectively, on $K_0,K_2,K_3$ and the vertices
$$\psi_0^+(z_*) =E, \psi_0^+(z_2) = C, \psi_0^+(z_3) = D.$$

Similarly, when $t \nearrow 1$, the map $\psi_t$ converges to a map
$\psi_1= \psi_1^-\;\sharp\;\psi_1^+ $ from  $\Sigma_1$ to $\C$. Here
$\psi_1^-: \D^- \to \C$ is a holomorphic triangle with the three
sides, respectively, on $K_0,K_1,K_3$ and the vertices
$$\psi_1^-(z_0) = A, \psi_1^-(z_*) = F, \psi_1^-(z_3) =D,$$
and $\psi_1^+: \D^+ \to \C$ is a holomorphic triangle with the
three sides, respectively, on $K_1,K_2,K_3$ and the vertices
$$\psi_1^+(z_*) =F, \psi_1^+(z_1) = B, \psi_1^+(z_2) = C.$$
The bubbling pattern shows that $z_2 \to z_3=-1$ as $t\to 0$ and
$z_2 \to z_1 =1$ as $t \to 1$. This shows that the map
$$\Phi: [0;1] \to \text{Compactification}(\widetilde{\cP}^4), \; t \to \widetilde{P}_t,$$
has degree $1$ (modulo 2). Thus, generically, for every
$\widetilde{P} \in Z$ (where the finite set $Z$ was defined in Step
1) there exists an odd number of values of $t$ with
$\widetilde{P}_t=\widetilde{P}$.

\medskip{\sc Step 3:} We use the notations of Steps 1 and 2.
Consider all pairs $(P\in \cP^4, \phi: P\to M)$ such that the
image of $\phi$ is a polygon with vertices $a,b,c,d$. The moduli
space of such pairs consists of $(\mu^3(a,b,c),d)$ points. We have
seen that each such point and every quadrilateral $ABCD\subset T^*
S^1$ as above together contribute $1$ (modulo $2$) to the
coefficient
$$(\mu^3(a \otimes A, b \otimes B, c \otimes C), d \otimes D).$$
The analysis in Step 2 shows that the number of such
quadrilaterals $ABCD$ equals $(A\cdot B \cdot C, D)$. We conclude
that
$$(\mu^3(a \otimes A, b \otimes B, c \otimes C), d \otimes
D)=(\mu^3(a,b,c),d)\cdot (A\cdot B \cdot C, D),$$ which immediately
yields formula \eqref{eq-prod-3}. \Qed

\medskip
We refer  to \cite{Loday} and references therein for an algebraic
discussion on $A_\infty$-operations for a tensor product of
$A_\infty$-algebras.

\subsection{Application to Poisson brackets invariants}
\label{sec-mainlagthm}

The next result is a more precise version of
Theorems~\ref{thm-triang} and \ref{thm-quadrilat} stated in the
introduction. Let $\cL=(L_0,\ldots,L_{k-1})$, $k=3$ or $4$, be a
collection of Lagrangian submanifolds of a geometrically bounded
symplectic manifold. Assume that $L_i$'s are in general position,
satisfy conditions (F1)-(F4) and for $k=4$ satisfy the
intersection condition \eqref{eq-4inter-again}. Put
$$A_3 = A(L_0,L_1,L_2;2n), A_4= A(L_0, L_1, L_2, L_3;3n-1),$$
$$A'_3 = \max \{A(L_0,L_1,L_2;j) : j \in [2n-3;2n]\},$$
$$A'_4 = \max \{A(L_0,L_1,L_2,L_3;j) : j \in [3n-5;3n-1]\}.$$

\begin{thm}\label{thm-lag-main} Assume that the operation
$\mu^{k-1}$ in the Lagrangian Floer homology of $\cL$ does not
vanish. Then
\begin{itemize}
\item[{(i)}] If $k=3$, then
\begin{equation}\label{eq-pblag3-again}
pb_3(L_0,L_1,L_2) \geq \frac{1}{2A_3}
\end{equation}
and
\begin{equation}\label{eq-pblag3-again-stab}
pb_3(\s\, L_0,\s\, L_1,\s\, L_2) \geq \frac{1}{2A'_3}.
\end{equation}
\item[{(ii)}] If $k=4$, then
\begin{equation}\label{eq-pblag4-again}
pb_4(L_0,L_2,L_1,L_3) \geq 1/A_4
\end{equation}
and
\begin{equation}\label{eq-pblag4-again-stab}
pb_4(\s\, L_0,\s\, L_2,\s\, L_1, \s\, L_3 ) \geq 1/A'_4.
\end{equation}
\end{itemize}
\end{thm}

\medskip\noindent{\bf Proof.} We shall prove part (ii) (the proof of (i) is analogous).
By Proposition~\ref{prop-op}, $pb_4(L_0,L_2,L_1,L_3)=
\inf{||\{F,G\}||}$, where $F=0$ in a neighborhood of $L_0$, $F=1$
in a neighborhood of $L_2$, $G=0$ in a neighborhood of $L_1$,
$G=1$ in a neighborhood of $L_3$. Given such functions $F,G\in
C^\infty_c (M)$, consider the family of forms $$\omega_s := \omega
-sdF\wedge dG.$$ Note that $$dF \wedge dG \wedge \omega^{n-1} =
\frac{1}{n}\{F,G\}\cdot \omega^n.$$ Thus
$$\omega_s^n =(1-s\{F,G\})\omega^n.$$ Therefore
the form $\omega_s$ is symplectic for all
$$s \in I:= [0; 1/||\{F,G\}||).$$
A straightforward application of the Stokes formula shows that the
deformation $\omega_s$, $s \in I$, satisfies the assumptions
(D1)-(D4) of Section~\ref{subsec-prelimLFH} above. For instance,
in order to verify that the collection $\cL$ is of finite type for
every $s$, observe that $\int_\alpha dF\wedge dG = 1$ for every
quadrilateral $\alpha$ with the boundary on $\cL$ and hence
\begin{equation}\label{eq-quadr-area-bound}
\omega_s(\alpha) \leq \omega (\alpha) -s.
\end{equation}
At the same time the Maslov class $m(\alpha)$ does not change in the
process of deformation and hence the finite type condition for
$\omega_s$ follows from the one for $\omega$.

Choose a generic family of almost complex structures $J_s$
compatible with $\omega_s$. By
Proposition~\ref{prop-LFHinvariance}, the operation $\mu^3$ in the
Lagrangian Floer homology of $\cL$ with respect to $\omega_s$ does
not vanish. Thus for every $s \in I$ there exists a
$J_s$-holomorphic quadrilateral, say $\alpha$, with boundary on
$\cL$. The dimension of the moduli space of such quadrilaterals
equals $m(\alpha)-3n+1 = 0$ and thus the finite type condition
\eqref{eq-fintype} guarantees that $\omega(\alpha) \leq A_4$.
Thus, by \eqref{eq-quadr-area-bound}, $\omega_s (\alpha) \leq A_4
-s$. At the same time $\omega_s(\alpha)
> 0$ and hence $s \leq A_4$ for every $s \in I=[0;1/||\{F,G\}||)$. This yields
$||\{F,G\}|| \geq 1/A_4$, and inequality \eqref{eq-pblag4-again}
follows.

Let us pass to inequality \eqref{eq-pblag4-again-stab}. Choose a
collection $\cK= (K_0, K_1, K_2,K_3)$ of four generic sections of
$T^*S^1$ and take $\epsilon
>0$. Applying
Proposition~\ref{prop-fintype-stab} and using the product formula
\eqref{eq-prod-3} we get that the collection $(L_i \times \epsilon
K_i)_{i=0,1,2,3}$ also satisfies the assumptions of the theorem
and hence, by \eqref{eq-pblag4-again},
$$pb_4(L_0 \times \epsilon K_0, L_2\times \epsilon K_2, L_1\times \epsilon K_1, L_3\times \epsilon K_3)
\geq \frac{1}{A'_4-\epsilon C(\cK)}.$$ Note that $L_i \times
\epsilon K_i$ converges (in the sense of
Section~\ref{subsec-defnot2}) to $\s L_i$ as $\epsilon \to 0$. Thus
inequality \eqref{eq-pblag4-again-stab} immediately follows from
Corollary~\ref{cor-Hausdorff}. This completes the proof. \Qed

\medskip
\noindent
\begin{rem}\label{rm-stb-more}{\rm Let $K_0,\ldots,K_{k-1}$ be
arbitrary exact sections of $T^*S^1$ (not necessarily in general
position). Put $L'_i:= L_i \times K_i$. The same proof as above
shows that under the assumptions of the theorem the Poisson bracket
invariants $pb_3(L'_0,L'_1,L'_2)$ (when $k=3$) and
$pb_4(L'_0,L'_2,L'_1,L'_3)$ (when $k=4$) are positive. }\end{rem}

\medskip
\noindent
\begin{rem}\label{rem-nohandles} {\rm The results of Section \ref{sec-pb-pseudo-holo}
extend verbatim to the case when Lagrangian submanifolds from our
collections are not necessarily compact, but rather geometrically
bounded (see \cite[p.286]{ALP}), that is properly embedded with
``nice behavior" at infinity. In this case we should also assume
that the number of intersection points of each pair of
submanifolds is finite. Let us apply this remark to the quadruple
of circles $X_0,X_1,Y_0,Y_1$ on the torus $\T^2$ considered at the
end of Example~\ref{exam-pb4m3}. Fix a square $\Pi$ on $\T^2$
whose edges (in counter-clockwise cyclic order) lie on
$X_0,Y_0,X_1,Y_1$. Take any lift of the contractible curve
$\partial \Pi$ to the universal cover $\R^2 \to \T^2$. Its edges
lie on some lifts $\widetilde{X}_0,\widetilde{Y}_0,
\widetilde{X}_1,\widetilde{Y}_1$ of  $X_0,Y_0,X_1,Y_1$
respectively. The quadruple $\cL$ of lines
$\widetilde{X}_0,\widetilde{Y}_0, \widetilde{X}_1,\widetilde{Y}_1$
on $\R^2$ forms a collection of finite type. Take a generic
quadruple $\cK=(K_1,K_2,K_3,K_4)$ of exact Lagrangian sections of
$T^*S^1$ and put
$$\widehat{X}_0 = \widetilde{X}_0 \times K_1,  \widehat{X}_1 = \widetilde{X}_1 \times K_2,
\widehat{Y}_0 = \widetilde{Y}_0 \times K_3,\widehat{Y}_1 =
\widetilde{Y}_1 \times K_4.$$ Consider a pair $(F,G) \in
\cF'_4(X_0 \times K_1, X_1 \times K_2, Y_0 \times K_3, Y_1 \times
K_4)$. Consider the deformation $\omega_s = \omega - sdF \wedge
dG$ of the symplectic form $\omega$ on $\T^2 \times T^*S^1$. Let
$J_s$ be an $\omega_s$-compatible family of almost complex
structures on $\T^2 \times T^*S^1$. Denote by
$\widetilde{\omega}_s$ and $\widetilde{J}_s$ the lifts of
$\omega_s$ and $J_s$ to $\R^2 \times T^*S^1$. The periodicity of
$\widetilde{\omega}_s$ and $\widetilde{J}_s$ with respect to the
group $\Z^2$ acting on the $\R^2$-factor guarantees that
Lagrangian submanifolds
$\widehat{X}_0,\widehat{Y}_0,\widehat{X}_1,\widehat{Y}_1$ remain
geometrically bounded for every $s$ whenever $\omega_s$ is
symplectic. Moreover, the $\mu^3$ operation is well defined and
does not vanish: indeed, it does not vanish for the quadruple of
lines $\cL$ on $\R^2$ due to the contribution of the square $\Pi$,
and it survives the stabilization by the product formula
\eqref{eq-prod-3}. Thus we get a $\widetilde{J}_s$-holomorphic
quadrilateral $\widetilde{\Sigma}_s$ with the edges  on
$\widehat{X}_0,\widehat{Y}_0,\widehat{X}_1,\widehat{Y}_1$. Its
projection $\Sigma_s$ to $\T^2 \times T^*S^1$ satisfies
$$\int_{{\Sigma}_s} {\omega} \leq \text{Area}(\Pi) + \text{const}(\cK),\;\;
\int_{\partial \Sigma_s} FdG =1.$$ Applying
Proposition~\ref{prop-psh-curves-general} as in the proof of
Theorem~\ref{thm-lag-main} above we readily get that $pb_4(X_0
\times K_1, X_1 \times K_2, Y_0 \times K_3, Y_1 \times K_4) >0$.
This confirms the claim made in the end of
Example~\ref{exam-pb4m3}. }\end{rem}

\subsection {Lagrangian spheres and the triangle product}
\label{subsec-triangle-example}

Let $(M^{2n},\omega)$, $n\geq 2$, be
an exact convex symplectic manifold, meaning
there exists a 1-form $\theta$ on $M$ and an exhausting sequence of
compact manifolds with boundary $M_1\subset M_2\subset\ldots
\subset M$ such that $\omega = d\theta$ and for any $i$ the 1-form $\theta|_{\partial M_i}$
is contact.
%with contact type boundary, meaning that there exists a
%1-form $\theta$ on $M$ so that $\omega = d\theta$ and
%$\theta|_{\partial M}$ is a contact form.
Let $L_0$, $L_2$ be exact Lagrangian submanifolds of $(M,\omega)$
(meaning that the restrictions of $\theta$ on them are exact
1-forms). Let $L_1\subset M$ be a Lagrangian sphere. Assume $L_0\cap
L_1\cap L_2=\emptyset$ and all the $L_i$ intersect each other
transversally -- thus the collection $L_0, L_1,L_2\subset M$ is of
finite type.

Fix a diffeomorphism $f: S^n\to L_1$. This data allows to
associate to $L_1$ a compactly supported symplectomorphism
$\tau_{L_1}: M\to M$, called the {\it Dehn twist in $L_1$}. It
maps $L_1$ to itself. Therefore there is a canonical isomorphism
\begin{equation}
\label{eqn-Dehn-twist-isom} HF (\tau_{L_1}^{-1} (L_0), L_1)\cong HF
(L_0, L_1).
\end{equation}
Seidel showed \cite{Seidel-exact} that there is an
exact sequence:
$$\xymatrix{
HF(L_0,L_2)\ar[rr]& \ &  HF(\tau_{L_1}^{-1} (L_0),L_2)\ar[dl]\\
\ & HF (\tau_{L_1}^{-1} (L_0), L_1)\otimes HF(L_1,L_2) \ar[ul]_F &\
}$$
where the map
$$F: HF (\tau_{L_1}^{-1} (L_0), L_1)\otimes HF(L_1,L_2)\to HF(L_0,L_2)$$
is the composition of the isomorphism \eqref{eqn-Dehn-twist-isom}
and the triangle product
$$\mu^2: HF(L_0,L_1)\otimes HF(L_1,L_2)\to HF(L_0,L_2).$$
Therefore Seidel's exact
sequence implies that if
\begin{equation}
\label{eqn-HF-Dehn-twist} HF(L_0,L_2)\neq 0,\ HF(\tau_{L_1}^{-1}
(L_0), L_2)= 0,
\end{equation}
then the product $\mu^2$ is non-trivial.

We learned the following specific example of such a situation from
Ivan Smith \cite{Smith1} -- we thank him for explaining it to us.

Consider $\C^3$ with the complex coordinates $x,y,z$. Take a
smooth complex hypersurface $M$ in $\C^3$ given by the equation
$x^2 + y^2 + p(z) = 1$, where $p$ is a complex polynomial of
degree $5$ with $4$ non-degenerate critical points, say $z_0=0$,
$z_1=-1$, $z_2=i$, $z_3=1$. The symplectic structure on $\C^3$
induces the structure of an exact convex symplectic manifold on
$M$. The projection $\pi: M\to \C$ to the complex $z$-plane is a
Lefschetz fibration with the critical values $z_i$, $i=0,1,2,3$.

A smooth embedded path $\gamma: [0;1]\to \C$, which connects two
distinct critical values $z_i= \gamma (0)$ and $z_j=\gamma (1)$ and
does not pass through the other critical values, is called a {\it
matching path}. To a matching path $\gamma$ one can associate a
Lagrangian sphere $S\subset \pi^{-1} (\gamma) \subset M$, called a
{\it matching cycle}. (The construction is due to Donaldson, for
details see e.g. \cite{Seidel-book}, pp. 230-231. The sphere is
glued from two Lagrangian discs, called {\it Lefschetz thimbles},
coming out of the critical points $(0,0,z_i)$ and $(0,0,z_j)$ of
$\pi$ and having a common boundary which is a vanishing cycle in a
fiber of $\pi$).

Consider the matching paths $\gamma_{01}$, $\gamma_{02}$,
$\gamma_{03}$, $\gamma_{23}$ which are straight segments in $\C$
connecting, respectively, $z_0$ with $z_1$, $z_2$, $z_3$, and $z_2$
with $z_3$. Denote the corresponding matching cycles $S_{01}$,
$S_{02}$, $S_{03}$, $S_{23}$. An exact Lagrangian isotopy identifies
the matching cycle $S_{23}$ with $\tau_{S_{03}}^{-1} (S_{02})$
(under an appropriate identification of $S_{03}$ with $S^2$) -- see
\cite{Seidel-book}, p.232. We can perturb the matching cycles by
$C^\infty$-small exact Lagrangian isotopies so that those of them
that correspond to intersecting matching paths intersect each other
transversally at exactly one point and all the triple intersections
are empty.  (Obviously, matching cycles corresponding to disjoint
matching paths do not intersect each other).

Thus setting $L_0:= S_{02}$, $L_1:=S_{03}$, $L_2:=S_{01}$, we see
that $L_0$, $L_1$ and $L_2$ are Lagrangian spheres in $M$ such that
$L_0\cap L_1\cap L_2=\emptyset$ and which satisfy
\eqref{eqn-HF-Dehn-twist}. Therefore the triangle product
$$\mu^2: HF(L_0,L_1) \otimes HF(L_1,L_2) \to HF(L_0,L_2)$$
is non-trivial.

\section{Poisson bracket invariants and SFT}\label{sec-SFT}

In this section we prove Proposition \ref{prop-surfaces-sphere} by
using a method of Symplectic Field Theory \cite{EGH-SFT}.

\subsection{Lagrangian tori in $S^2 \times T^*S^1$}\label{subsec-Lagtorus}

Consider a symplectic manifold $V = S^2 \times T^*S^1$ equipped
with the split symplectic form $\omega_0$ so that the area of
$\gamma:= [S^2 \times \text{point}]$ equals $1$. Let $\Pi \subset
S^2$ be a disc with smooth boundary. Consider a Lagrangian torus
$L=\partial \Pi \times S^1$ in $V$.  The relative Hurewicz
morphism $\pi_2(V,L) \to H_2(V,L,\Z)$ is an isomorphism.  Denote
by $\alpha$ and $\beta$ the elements in $H_2(V,L,\Z)$ generated by
$\Pi \times \{\text{point}\}$ and $\overline{S^2 \setminus \Pi}
\times \{\text{point}\}$ respectively so that $\alpha + \beta =
\gamma$.

\begin{thm} \label{thm-SFT}
Let $\omega_\tau$, $\tau \in [0;1]$, be a smooth deformation of
$\omega_0$ through symplectic forms such that $L$ remains
$\omega_\tau$-Lagrangian for all $\tau$.   Then
\begin{equation}\label{eq-positivity}
 \omega_1 (\alpha) > 0\;\;\text{and}\;\; \omega_1(\beta) > 0.
\end{equation}
\end{thm}

\medskip
\noindent {\bf Proof: \footnote{ This proof is due to Richard Hind.
We thank him for his help and a considerable shortening of our
original argument.}}

 Given a Riemann surface with boundary, say $C$, attach
a punctured disc to each of its boundary components. The resulting
Riemann surface is denoted by $\widehat{C}$.

 Assume on the contrary that $\omega_1(\alpha) \leq 0.$
Since $\omega_0(\alpha) > 0$, there exists $t \in (0;1]$ such that
$\omega_t(\alpha) = 0$.

 We equip the torus $L$ with the Euclidean metric, make
an appropriate choice of an $\omega_t$-compatible almost complex
structure $J$ on $V$ and perform the stretching-the-neck procedure
near $L$ as in \cite{EGH-SFT}. As a result of the stretching, we
get an almost complex structure $J_b$ with a negative cylindrical
end on $V \setminus L$ and an almost complex structure $J_w$ on
$T^*L$ with a positive cylindrical end. Let us emphasize that the
structure $J_b$ is tamed by $\omega_t$.

The manifold $V$ is foliated by $J$-holomorphic spheres in the class
$\gamma$ \cite{Gromov-pshc}. The compactness theorem of \cite{BEHWE}
guarantees that after stretching the neck some of these spheres
split into a collection of multi-level pseudo-holomorphic curves
asymptotic to closed orbits of the Euclidean geodesic flow on $L$.
Without loss of generality we shall assume that there are just two
levels. Thus there exists
\begin{itemize}
\item  a partition of the sphere $S^2$ (equipped with the standard
complex structure) by $K>0$ boundary circles into blue and white
domains $B_1,\ldots,B_N$ and $W_1,\ldots,W_M$ so that any two
domains with a common boundary component have different colors;
\item pseudo-holomorphic maps $\phi_i: \widehat{B_i} \to (V
\setminus L, J_b)$ and $\psi_i: \widehat{W_i} \to (T^*L, J_w)$
whose negative (resp. positive) asymptotic ends are closed orbits
of the Euclidean geodesic flow on $L$.
\end{itemize}

By obvious topological reasons, there are at least two discs among
the domains of our partition. Since all Euclidean geodesics on the
two-torus are non-contractible, no white domain can be a disc.
Thus there are $N \geq 2$ blue domains.

Persistence of the fibration by $J$-holomorphic spheres in the
class $\gamma$ yields that $\omega_t(\gamma)>0$. Write the
relative homology class of $\phi_i(\widehat{B}_i)$ as $p_i \alpha +
q_i \gamma$. Since the relative homology class of each
$\psi_i(\widehat{W}_i)$ is zero (because $\pi_2 (T^* L,L)\cong
0$), the classes $\phi_i(\widehat{B}_i) = p_i \alpha + q_i \gamma$
add up to $\gamma$ and we get
\begin{equation}\label{eq-sumK}
\sum_{i=1}^N q_i =1.
\end{equation}
Since $\omega_t$ tames $J_b$, we have that $\omega_t(p_i \alpha +
q_i \gamma)>0$. Since  $\omega_t(\alpha)=0$ and
$\omega_t(\gamma)>0$, we necessarily have that $q_i >0$ for all
$i$. In view of \eqref{eq-sumK}, this contradicts $N \geq 2$.
Therefore $\omega_1(\alpha)>0$. Similarly, $\omega_1(\beta)>0$.
\qed

\subsection{An application to Poisson bracket invariants}

\noindent {\bf Proof of Proposition~\ref{prop-surfaces-sphere}:}
Arguing as in Section~\ref{exam-pb3-surfaces} we see that it
suffices to prove the lower bound
\begin{equation}
\label{eq-vsp-1000} pb_4(a_1 \times K, a_3 \times K, a_2 \times K,
a_4 \times K) \geq \max(1/A, 1/(B-A))\;\end{equation} assuming
that the boundary $\partial \Pi$ is smooth. Without loss of
generality, let $K=S^1$ be the zero section of $T^*S^1$ and $B =
\text{Area}(S^2)=1$.  Pick two functions $F,G \in \cF'_4(a_1
\times K, a_3 \times K, a_2 \times K, a_4 \times K)$ and consider
the deformation $\omega_s := \omega_0 -sdF \wedge dG$ of the split
symplectic form $\omega_0$ on $V= S^2 \times T^*S^1$. As we have
seen in Section~\ref{subsubsec-def-lowerbd}, $\omega_s$ is
symplectic for all $s < 1/||\{F,G\}$ and (in the notations of
Section~\ref{subsec-Lagtorus} above) $\omega_s(\alpha) = A-s$
(because $\int_\alpha dF\wedge dG = \int_{\partial\Pi} FdG =1$).
By Theorem \ref{thm-SFT}, $\omega_s(\alpha)>0$ and hence $s < A$.
Therefore
$$p:= pb_4(a_1 \times K, a_3 \times K, a_2 \times K, a_4 \times K)
\geq 1/A.$$

Applying the same argument to the quadrilateral $\Pi':=
\overline{S^2 \setminus \Pi}$ we get that
$$p':= pb_4(a_1 \times K,
a_4 \times K, a_2 \times K, a_1 \times K) \geq 1/(1-A).$$ By the
symmetry of the Poisson bracket invariants, $p=p'$ and hence we
get inequality \eqref{eq-vsp-1000}. \qed

\section{A vanishing result for $pb_4$} \label{sec-vanishing}

In this section we prove Proposition~\ref{prop-baby-example}. Let
us introduce the following terminology: Assume that $S$ is a
finite simplicial complex and $M$ is a manifold. Let $\phi: S \to
M$ be a homeomorphism from $S$ to its image $\phi(S)$ so that the
restriction of $\phi$ to every simplex is a smooth embedding. We
refer to the image $X:=\phi(S)$ as to an embedded simplicial
complex in $M$. We denote by $\dim X$ the maximal dimension of a
simplex from $S$.

\begin{prop}\label{prop-pbzero} Let $X_0,X_1$ be disjoint
embedded simplicial complexes in a symplectic manifold $(M,\omega)$.
Assume that
\begin{equation}\label{eq-dim}
\dim X_0 + \dim X_1 \leq 2n-2.
\end{equation}
Then for every pair of disjoint compact subsets $Y_0,Y_1$ of $M$
$$pb_4(X_0,X_1,Y_0,Y_1) = 0.$$
\end{prop}

\medskip\noindent{\bf Proof.} Assume on the contrary that $pb_4(X_0,X_1,Y_0,Y_1) = p
 >0$. Fix neighborhoods $U_i$ of $Y_i$, $i=0,1$.
Take any function $H \in C^\infty_c(M)$ so that $H=0$ on $U_0$ and
$H=1$ on $U_1$. Denote by $h_t$ the Hamiltonian flow generated by
$H$.  Put $T=2/p$ and set $Z:= \bigcup_{|t| \leq T} h_t(X_0)$. The
dimension formula \eqref{eq-dim} and a standard transversality
argument yield the existence of an arbitrary $C^0$-small
Hamiltonian diffeomorphism $f$ of $M$ supported near $X_1$ such
that
\begin{equation}\label{eq-displace}
f(X_1) \cap Z = \emptyset.
\end{equation}
Since $f$ is $C^0$-small, we can assume that $f(Y_i) \subset U_i$
for $i=0,1$. Moreover, since $f$ is supported near $X_1$ and $X_0
\cap X_1 =\emptyset$, we have $f(X_0)=X_0$. By symplectic
invariance of $pb_4$,
$$pb_4(f(X_0), f(X_1),f(Y_0),f(Y_1)) = p.$$
Theorem~\ref{thm-chords1} guarantees that there exists a Hamiltonian
chord of $h_t$ joining $f(X_0)=X_0$ and $f(X_1)$ of time-length
$\leq 1/p$. This contradicts property \eqref{eq-displace}, and hence
$p=0$. \Qed

Before proceeding further, let us introduce the following
notation. Consider the annulus $A = [-1;1] \times S^1$. Let
$x_1,\ldots,x_k, y_1,\ldots,y_l$ be pairwise distinct points in
$S^1$. Consider a 1-dimensional simplicial complex
$$A(x_1,\ldots,x_k, y_1,\ldots,y_l)\subset A$$ defined by
$$A(x_1,\ldots,x_k, y_1,\ldots,y_l) = (\{0\} \times S^1 ) \cup
\bigcup_{i=1}^k
([0;1] \times \{x_i\}) \cup \bigcup_{j=1}^l ([-1;0] \times \{y_j\}).$$

\medskip
\noindent {\bf Proof of Proposition~\ref{prop-baby-example}.}
Recall that
$M$ is a closed symplectic surface, $\Pi \subset M$
is a quadrilateral with edges denoted (in the cyclic order) by $a_1,a_2,a_3,a_4$ and
$K_1,\ldots,K_4$ is a generic quadruple of exact sections of $T^*S^1$.
Put $P_i=a_i \times K_i$. We have to show that
$$
p:=pb_4(P_1,P_3,P_2,P_4) = 0. $$

We will use the cyclic convention for the indices $i=1,2,3,4$ (that is $4+1=1$, $1-1=4$).

Choose in the obvious way parameterizations $\phi_i: A \to P_i$,
$i=1,2,3,4$,
such that $P_i \cap P_{i+1}$
consists of a finite number of points of the form
$$\phi_i(1,x_j) = \phi_{i+1}(-1,y_j),\ j =
1,\ldots,N(i),$$
for some $N(i)\in\N$.
Put
$$S_i = A(x_1, \ldots, x_{N(i)}, y_1,\ldots,y_{N(i-1)}).$$
Fix $\epsilon >0$ and observe that one can ``collapse" the annulus
$A$ to an $\epsilon$-neigh\-bor\-hood $S^\epsilon_i$ of $S_i$.
More precisely, there exists a family of embeddings $\psi_i^t: A
\to A$, $t \in [0;1]$, such that $\psi_0 = \id$, $\psi^1_i(A)=
S^\epsilon_i$ and $\psi_i^t=\id$ near $S_i$ for all $t$.

Put $Q_i = \phi_i(S_i)$ and $Q^\epsilon_i = \phi_i(S^\epsilon_i)$.
Observe that $P_i \cap P_{i+1}= Q_i \cap Q_{i+1}$ and therefore the
isotopies
$$\theta_i^t:= \phi_i \circ \psi_i^t: A \to P_i$$
have disjoint supports for distinct $i$. Since each $P_i$ is a
Lagrangian submanifold of $M \times T^*S^1$, the isotopies
$\theta_i^t$, $i=1,2,3,4$, extend {\bf simultaneously} to an ambient
Hamiltonian isotopy $\theta^t$ of $M \times T^*S^1$. By the
symplectic invariance of $pb_4$,
\begin{equation}\label{eq-p-epsilon}
pb_4(Q^\epsilon_1 ,Q^\epsilon_3 ,Q^\epsilon_2 ,Q^\epsilon_4)=p
\end{equation}
for all $\epsilon$.  Note that $Q^\epsilon_i \to Q_i$ as $\epsilon
\to 0$, where the convergence is understood in the sense of
Section~\ref{subsec-defnot2}. Thus, by
Proposition~\ref{prop-Hausdorff},
$$p \leq q:= pb_4(Q_1,Q_3,Q_2,Q_4).$$
Since each $Q_i$ is a one-dimensional embedded simplicial complex in
$M \times T^*S^1$, Proposition~\ref{prop-pbzero} yields $q=0$. Hence
$p=0$. This completes the proof. \Qed

\medskip
\begin{rem}\label{rem-pb3-vanishing} {\rm Let $X,Y,Z$ be
compact subsets of a symplectic manifold $(M,\omega)$ with $X \cap
Y \cap Z = \emptyset$. Assume that $X$ and $Z$ are embedded
simplicial complexes with $\dim X + \dim Z \leq 2n-2$. We claim
that $pb_3(X,Y,Z) = 0$. Indeed, represent $Z$ as the union $Z_1
\cup Z_2$ of two compact embedded complexes of the same dimension
so that $X \cap Z_1 = Y \cap Z_2 = \emptyset$. Combining
Proposition~\ref{prop-pbzero} with inequality \eqref{eq-pb3-pb4}
we get that
$$0 = pb_4(X,Z_1,Y,Z_2) \geq pb_3(X,Y,Z),$$
and the claim follows. }\end{rem}

\section{Discussion and further directions}
\label{sec-discussion}

\subsection{Hamiltonian chords and optimal control}

Hamiltonian chords joining two disjoint subsets of a symplectic
manifold appear in the mathematical theory of optimal control. For
instance, the shortest geodesic between two closed submanifolds of
a Riemannian manifold can be interpreted as a chord of the
geodesic flow joining their Lagrangian co-normals in the cotangent
bundle. As we have mentioned above (see discussion after
Theorem~\ref{thm-torus-chords}), in some situations such chords
can be captured by the Poisson brackets invariants.

A similar interpretation can be given to the extremals provided by
Pontryagin's maximum principle for an optimal-time control problem
with variable end-points \cite{Pontryagin et al}.  The Hamiltonian
functions appearing in this context are degree-one homogeneous in
the momenta and in general are not proper. It would be interesting
to understand whether methods of symplectic and contact topology
can detect Hamiltonian chords in this context.

The minimal time-length $T(X_0,X_1,G)$ introduced in
Section~\ref{subsec-intro-chords} has the following counterpart in
control theory. Let $(M,\omega)$ be a symplectic manifold (the
so-called state space), $U$ be the input space and $G: M \times U
\to \R$ be a controlled Hamiltonian (see e.g. \cite[Ch.
12]{Nijmeijer-vdSchaft} or \cite[Section 4.9.5]{Ivan-applied}).
For any path $u(t)$, $t\in\R$, in $U$ the function $(x,t)\mapsto
G(x,u(t))$ can and will be viewed as a time-dependent Hamiltonian
on $M$. The optimal-time control problem with the initial set
$X_0$ and the terminal set $X_1$ is to find the minimal possible
time $T:=T_{min}(X_0,X_1,G)$ and a sufficiently regular control $u
: [0;T] \to U$ so that  the Hamiltonian flow generated by
$G(x,u(t))$ admits a trajectory $x(t)$ with $x(0) \in X_0$ and
$x(T) \in X_1$.

In general, even if the minimal time $T_{min}(X_0,X_1,G)$ is finite,
there is no reason for it to remain uniformly bounded under
$C^0$-small perturbations of the controlled Hamiltonian $G$: such
perturbations may drastically change the dynamics. Suppose now that
$pb_4(X_0,X_1,Y_0,Y_1) = p >0$ for some subsets $Y_0,Y_1 \subset M$,
and in addition
$$\min_{Y_1}G(x,u_*) -\max_{Y_0}G(x,u_*) = a >0\;$$
for some input $u_* \in U$. Taking the constant control $u(t)
\equiv u_*$ and applying Theorem~\ref{thm-chords1} we get that
\begin{equation}\label{eq-control-bound}
T_{min}(X_0,X_1,G) \leq (ap)^{-1}.
\end{equation}
This upper bound for $T_{min}(X_0,X_1,G)$ is robust under
$C^0$-small perturbations of the controlled Hamiltonian $G$.

The methods of proving the bound \eqref{eq-control-bound}
developed in the present paper are very much disjoint from the
standard tools of control theory. It would be interesting to
explore their possible interrelations. As a starting point one may
consider the simplest case of an affine Hamiltonian control
system. Here the controlled Hamiltonian $G$ is of the form
$$G(x,u) = G_0(x) + \sum_{i=1}^k u_iG_i(x),$$
and the input space $U$ is the cube $$\{|u_i| \leq 1,
i=1,\ldots,k\} .$$ Suppose also that $X_0,X_1,Y_0,Y_1$ are
closed Lagrangian submanifolds in $M$ as in the setting of
Section~\ref{subsec-intro-pb-lagrangian} above. The maximum
principle with transversality conditions at $X_0$ and $X_1$
provides a wealth of information about time-optimal trajectories
joining $X_0$ and $X_1$ (note that these extremals may possess
switches of the control parameters which manifest the so-called
``bang-bang" control). It would be interesting to design specific
examples where the upper bound \eqref{eq-control-bound} can be
deduced from the maximum principle.

\subsection{Higher Poisson bracket invariants}

In this section we suggest a generalization of the Poisson brackets
invariants to ordered collections $\cX=(X_1,\ldots,X_N)$ of
compact subsets of a symplectic manifold $(M,\omega)$. We start
with the following data.

\begin{figure}
 \begin{center}
\scalebox{0.5}{\includegraphics*{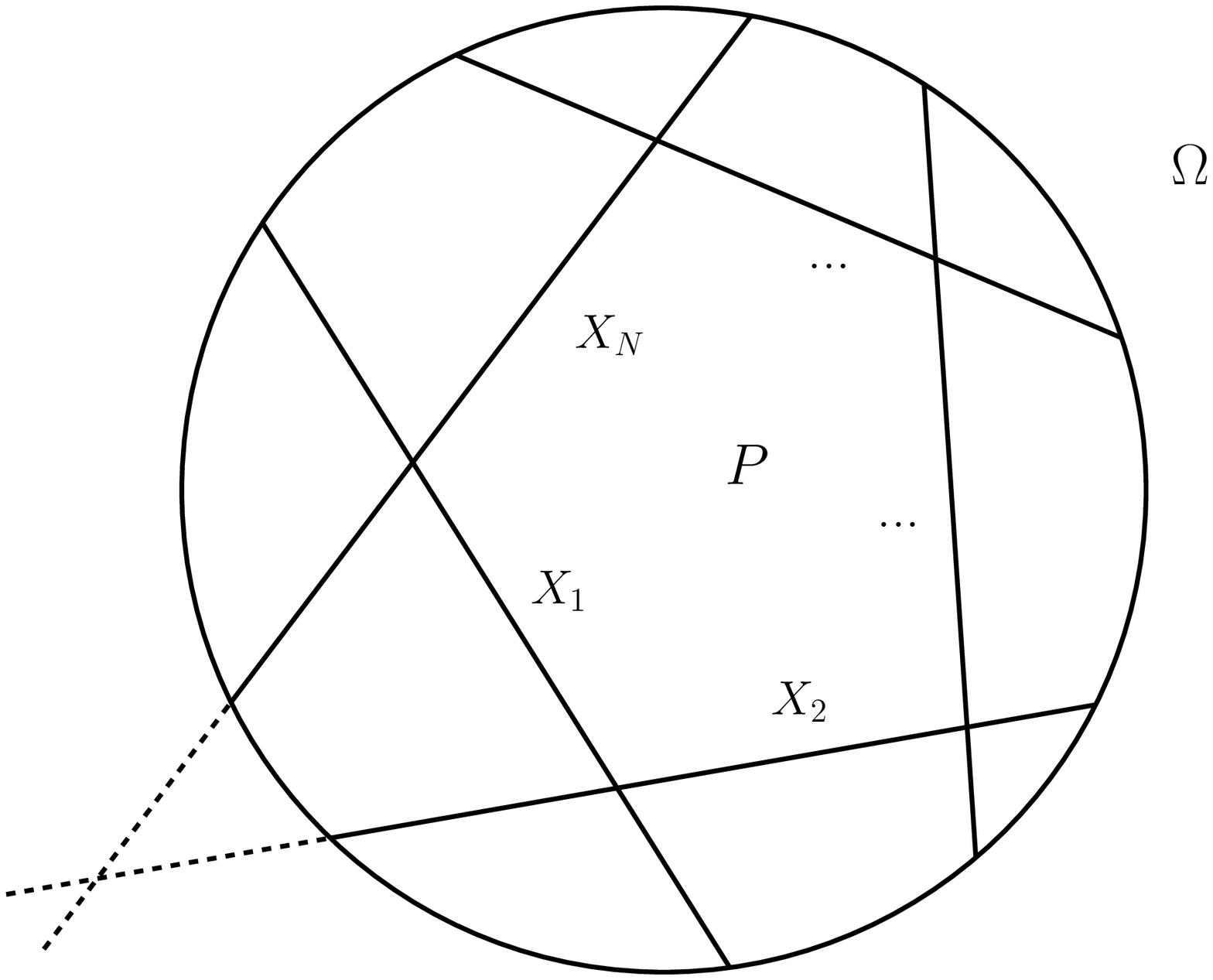}}
 \end{center}
  \caption{Condition $\diamondsuit$}
   \label{fig2}
\end{figure}

Let $a_i(s,t) = \alpha_i s + \beta_i t + \gamma_i$,
$i=1,\ldots,N$, be a collection of $N$  affine functions on the
plane $\R^2$ defining a convex  Euclidean polygon $P = \bigcup_i
\{a_i \geq 0\}$. Let $L_i = \{a_i = 0\}$ be the line containing
$i$-th edge of $P$. Fix a convex open domain $\Omega \subset \R^2$
containing $P$ with the following property:

\medskip
\noindent{\bf Condition $\diamondsuit$:}  For every $i\neq j$ the
set $L_i \cap L_j \cap \text{Closure}(\Omega)$ is either empty or
consists of a vertex of $P$ (i.e. the closure of $\Omega$ does not
contain the intersection points of the lines $L_i$ that are not
vertices of $P$), see Figure \ref{fig2}. In addition, we assume that $0 \in \Omega$ if $M$
is an open manifold.

\medskip
\noindent  A collection $\cX=\{X_1,\ldots,X_N\}$ of $N \geq 3$
compact subsets of $M$ is called {\it cyclic} if $X_i\cap X_j =
\emptyset$ unless $i$ and $j$ are equal or differ by $1$ (the
indices are taken modulo $N$). Given a cyclic collection $\cX$, denote
by $\cF_N(\cX,P,\Omega)$ the class of pairs of functions $(F,G)
\in \cF$ which satisfy
$$(F(x),G(x)) \in \Omega \;\; \forall x \in M$$
and
$$a_i(F(x),G(x)) \leq 0 \;\; \forall i = 1,\ldots,N, \;\; \forall x \in
X_i.$$ It is easy to see that  this class is non-empty: First, we
define $F$ and $G$ near $X_i \cap X_{i+1}$ as the $s-$ and
$t-$coordinates of the corresponding vertex of $P$, then we extend
$F,G$ to a neighborhood of $X_i$ so that $\alpha_i F(x) + \beta_i
G(x) +\gamma_i =0$ and $(F(x),G(x)) \in \Omega$, and finally we
cut off $(F,G)$ outside the union of $X_i$'s.

Define the Poisson bracket invariant $pb_N$ of a cyclic
collection $\cX$ by
$$pb_N(\cX,P,\Omega):= \inf ||\{F,G\}||,$$
where the infimum is taken over all $(F,G) \in
\cF_N(\cX,P,\Omega)$. The previously defined invariants $pb_3$ and
$pb_4$ are particular cases of this construction: the invariant
$pb_3(X_1,X_2,X_3)$ corresponds to the case when $N=3$, $\Omega =
\R^2$ and $P= \{s \geq 0, t \geq 0, s+t \leq 1\}$, while
$pb_4(X_1,X_3,X_2,X_4)$ (mind the order of the subsets)
corresponds to the case when $N=4$, $\Omega = \R^2$ and $P=\{0\leq
s\leq 1, 0 \leq t \leq 1\}$.

Higher Poisson bracket invariants can be studied along the lines of
the present paper. Denote by $\cF'_N(\cX,P,\Omega)$ the class of
pairs of functions $(F,G) \in \cF_N(\cX,P,\Omega)$ which satisfy
$$(F(x),G(x)) \in P \;\; \forall x \in M.$$ We claim that
$$pb_N(\cX,P,\Omega)=\inf ||\{F,G\}||,$$
where the infimum is taken over all $(F,G) \in
\cF'_N(\cX,P,\Omega)$. Indeed, condition $\diamondsuit$ yields (cf.
the proof of Lemma~\ref{lem-cutoff} above) the following fact: for
every $\kappa
>0$ there exists $\delta(\kappa)> 0$, with $\delta(\kappa) \to 0$ as
$\kappa \to 0$, and a map $T=(T_1,T_2) : \Omega \to P$ which takes
$\{a_i \leq \delta\}$ to the edge $ \{a_i =0\} \cap P$ for all
$i=1,\ldots,N$ and which satisfies
$$||\{T_1,T_2\}_{\R^2}|| \leq 1+\kappa.$$
The claim readily follows from this fact (cf. the proof of
Proposition~\ref{prop-op} above).

Suppose now that the symplectic manifold $(M,\omega)$ is closed
and admits a symplectic quasi-state satisfying the PB-inequality
(see \eqref{eq-frol} above).  We have then the following analogue
of Theorem~\ref{thm-qstates-pb4-dynamics}:  Assume that the sets
$Y_i:= \bigcup_{j=i}^{i+N-3}X_j$ are superheavy for all $i$ (we
use the cyclic convention for the indices with $N+1=1$). Then
$pb_N(\cX,P,\Omega) \geq c
>0$, where the constant $c$ depends only on the polygon $P$ and
on the constant $K$ entering the PB-inequality \eqref{eq-frol}.

Let us sketch a proof. Put $A_i =  \prod_{k=i}^{i+N-3} a_i$ and
$A= \sum_i A_i$. Observe that the function $A$ is strictly
positive on the polygon $P$. Indeed, the functions  $A_i$ are
non-negative, while the intersection of their zero sets is empty.
Put $c' = \min_P A>0$. Take a pair $(F,G) \in
\cF'_N(\cX,P,\Omega)$.  Define functions $H_i:=A_i(F,G)$ and $H:=
A(F,G)$ on $M$. It follows that $\zeta(H) \geq c'$. At the same
time $H_i$ vanishes on $Y_i$, and hence the superheaviness of
$Y_i$ yields $\zeta(H_i) = 0$. Thus if $||\{F,G\}||$ is
sufficiently small, the functions $H_i$ ``almost commute", and
hence, by the PB-inequality, $|\zeta(H)|$ should be strictly
smaller than $c'$ (which does not depend on $F,G$) yielding a
contradiction. Therefore $||\{F,G\}||$ cannot be small which
yields a lower bound on $pb_N$.

In case when $X_i$'s are Lagrangian submanifolds, one can pursue the
second approach to the positivity of $pb_N$:  deform the symplectic form
and study persistent pseudo-holomorphic polygons coming from Donaldson-Fukaya
category. Similarly to pseudo-holomorphic tri\-angles and
quadri\-la\-terals used for the study of $pb_3$ and $pb_4$,
pseu\-do-holomorphic polygons with a higher number of vertices can
be used to give a positive lower bound on $pb_N (\cX,P,\Omega)$. The
existence of persistent pseudo-holomorphic polygons can be extracted
from the higher (Massey-type) product
$$\mu^N: HF(X_1,X_2)\otimes\ldots\otimes HF(X_{N-1},X_{N})\to
HF(X_{1},X_N),$$ provided it is well-defined and non-trivial.

It would be interesting to explore applications of higher Poisson
bracket invariants beyond the applications of $pb_3$ or $pb_4$
described in this paper.

\subsection{Vanishing of Poisson bracket invariants}

According to the standard symplectic philosophy, the positivity of
the Poisson bracket invariant $pb_k(X_1,\ldots,X_k)$ should
manifest ``symplectic rigidity" of the collection of compact
subsets $X_1,\ldots,X_k$. Thus for a ``flexible" collection,
$pb_k$ should vanish. Proposition~\ref{prop-pbzero} and
Remark~\ref{rem-pb3-vanishing} above confirm this intuition:
$pb_4(X_0,X_1,Y_0,Y_1) = 0$ provided $\dim X_0 + \dim X_1 \leq
\dim M-2$ and $pb_3(X,Y,Z)=0$ provided $\dim X + \dim Z \leq \dim
M -2$.

The next natural test is the case when  $\dim M =4$ and our subsets
are two-dimensional surfaces.

\begin{prop}
Let $ (M^4,\omega) $ be a closed symplectic $4$-manifold and let $
X,Y \subset M $ be closed $2$-dimensional submanifolds such that
at any intersection point $ p \in X \cap Y $, the tangent spaces $
T_pX, T_pY \subset T_p M $ are transversal and symplectically orthogonal. Let $Z
\subset M$ be any compact set such that $ X \cap Y \cap Z =
\emptyset $. Then there exist smooth functions $ F,G : M
\rightarrow \R $ such that $ F = 0 $ on $ X $, $ G = 0 $ on $ Y $
and $ F+G \geqslant 1 $ on $ Z $, and, moreover, we have $ \{ F,G
\} = 0 $ on $ M $. As a consequence, we have $ pb_3(X,Y,Z) = 0 $.
\end{prop}

\medskip\noindent{\bf Proof.}
 For any point $ p \in X \cap Y $ there exists a neighborhood
 $W_p$ of $p$ with Darboux coordinates $ (x_1,y_1,x_2,y_2) $, where
$ p = (0,0,0,0) $, such that $ X \cap W_p $ coincides with $ x_1 =
y_1 = 0 $ and $ Y \cap W_p $ coincides with $ x_2 = y_2 = 0 $. Let
$ p_1,p_2,\ldots,p_k$ be the intersection points of $ X $ and $ Y
$. Replacing, if necessary, the neighborhoods $ W_{p_i} $, $
i=1,2,\ldots,k $, by smaller ones we may assume that $
W_{p_1},W_{p_2},\ldots,W_{p_k} $ are pairwise disjoint and $
W_{p_i} \cap Z = \emptyset $ for $ i=1,2,\ldots,k $. Take a small
$ a > 0 $ such that
$$ P_a := \{ (x_1,y_1,x_2,y_2) \, | \, x_1^2 + y_1^2 \leqslant
a^2,\;
 x_2^2 + y_2^2 \leqslant a^2 \} \subset W_{p_i} $$ for $ i=1,2,\ldots,k $.
Moreover, one can find tubular neighborhoods $ U_X $ of $ X $ and
$ U_Y $ of $ Y $ in $ M $ such that $ U_X \cap U_Y \subset
\cup_{i=1}^{k} W_{p_i} $, and $$ U_X \cap W_{p} =
 \{ (x_1,y_1,x_2,y_2) \in W_p \, | \,  x_1^2 + y_1^2 < a^2 \} ,$$
 $$ U_Y \cap W_{p} =
 \{ (x_1,y_1,x_2,y_2) \in W_p \, | \,  x_2^2 + y_2^2 < a^2 \} $$ for any $ p \in X \cap Y $.
Consider a smooth function $ u: [0;+\infty) \rightarrow [0;1] $
such that $ u(t) = 0 $ for $ t \in [0;a^2/3] $ and $ u(t) = 1 $
for $ t \in [2a^2/3;+\infty) $. Define functions
$$ f,g : \cup_{i=1}^{k} W_{p_i} \rightarrow \R $$ by $$
f(x_1,y_1,x_2,y_2) = u(x_1^2 + y_1^2),$$ $$ g(x_1,y_1,x_2,y_2) =
u(x_2^2 + y_2^2)$$ for $ (x_1,y_1,x_2,y_2) \in W_{p} $ for any $ p
\in X \cap Y $. One can easily find smooth functions $ F,G : M
\rightarrow [0;1] $, such that $ F(x) = 0 $ on $ X $, $ F(x) = 1 $
on $ M \setminus U_X $, $ G(x) = 0 $ on $ Y $, $ G(x) = 1 $ on $ M
\setminus U_Y $ and $ F(x) = f(x) $, $ G(x) = g(x) $ for $ x $
lying in a neighborhood of $ P_a \subset W_{p_i} $, $ i =
1,2,\ldots,k $. Then we will have $ \{ F,G \} = 0 $ on $ M $, $ F
= 0 $ on $ X $, $ G = 0 $ on $ Y $ and $ F+G \geqslant 1 $ on $ Z
$. \Qed

\medskip
A similar argument shows that if $X_0$ and $Y_0$ are closed
surfaces in a symplectic four-manifold which intersect
transversally and are symplectically orthogonal at each
intersection point, $pb_4(X_0,X_1,Y_0,Y_1) = 0$ for all $X_1,Y_1$
satisfying the intersection condition \eqref{eq-inter-4}.

We still do not know the answer to the following basic
question.

\begin{question}{\rm  Let $X,Y,Z \subset M^4$ be closed
$2$-dimensional {\bf non-Lagrangian} submanifolds with $X \cap Y
\cap Z =\emptyset.$ Is it true that $pb_3(X,Y,Z) =0$?}
\end{question}

\medskip
 The obvious analogue of this question for $pb_4$ is also
open.

\bigskip
\noindent {\bf Acknowledgements.} We thank Richard Hind for his
help with symplectic field theory and in particular for providing
us a short proof of Theorem~\ref{thm-SFT}. We are grateful to Paul
Seidel and Ivan Smith for useful consultations on Donaldson-Fukaya
category. In particular, Smith explained to us some interesting
examples where operations $\mu^2$ and $\mu^3$ do not vanish.  We
thank Paul Biran for a number of stimulating discussions.
We thank Strom Borman and an anonymous referee for comments and corrections --
in particular, Borman suggested a considerable improvement of the numerical constants in Section~\ref{subseq-pb-qs}.
A part
of this paper was written during our stay in MSRI, Berkeley. We
thank MSRI for the hospitality. The third-named author thanks the
Simons Foundation for sponsoring this stay.

\bibliographystyle{alpha}

\medskip

\begin{tabular}{l}
Lev Buhovsky\\
Department of Mathematics\\
University of Chicago\\
Chicago, IL 60637, USA\\
levbuh@gmail.com\\
\end{tabular}

\medskip

\begin{tabular}{l}
Michael Entov \\
Department of Mathematics\\
Technion - Israel Institute of Technology\\
Haifa 32000, Israel \\
entov@math.technion.ac.il \\
\end{tabular}

\medskip

\begin{tabular}{l}
Leonid Polterovich\\
Department of Mathematics\\
University of Chicago\\
Chicago, IL 60637, USA\\
and\\
School of Mathematical Sciences\\
Tel Aviv University\\
Tel Aviv 69978, Israel\\
polterov@runbox.com\\
\end{tabular}

\end{document}